# Functional penalised basis pursuit on spheres ☆

Matthieu Simeoni [1]

*École Polytechnique Fédérale de Lausanne (EPFL), Lausanne, Switzerland*

A B S T R A C T

In this paper, we propose a unified theoretical and practical spherical approximation framework for functional inverse problems on the hypersphere $\mathbb{S}^{d-1}$. More specifically, we consider recovering spherical fields directly in the continuous domain using functional penalised basis pursuit problems with generalised total variation (gTV) regularisation terms. Our framework is compatible with various measurement types as well as non-differentiable convex cost functionals. Via a novel representer theorem, we characterise their solution sets in terms of spherical splines with sparse innovations. We use this result to derive an approximate canonical spline-based discretisation scheme, with vanishing approximation error. To solve the resulting finite-dimensional optimisation problem, we propose an efficient and provably convergent primal-dual splitting algorithm. We illustrate the versatility of our framework on real-life examples from the field of environmental sciences.

## 1. Introduction

### 1.1. Spherical approximation: an overview

Many scientific inquiries in natural sciences, such as environmental and planetary sciences [1–3], acoustics [4] or astronomy [5–7], involve approximating a *spherical field* – a scalar quantity such as a function or measure defined over a *continuum* of directions, from a finite number of measurements acquired by probing sensors. During the reconstruction task, the physical evidence is compared to some prior model of the unknown spherical field, reflecting the analyst's a priori beliefs about the latter. In practice, a trade-off

☆ This work was in part supported by the Swiss National Science Foundation grant number 200021 181978/1, "SESAM - Sensing and Sampling: Theory and Algorithms".

*E-mail address:* matthieu.simeoni@epfl.ch.

[1] For part of this work, the author was also with the Foundations of Cognitive Solutions group in the IBM Research Laboratory of Zurich.

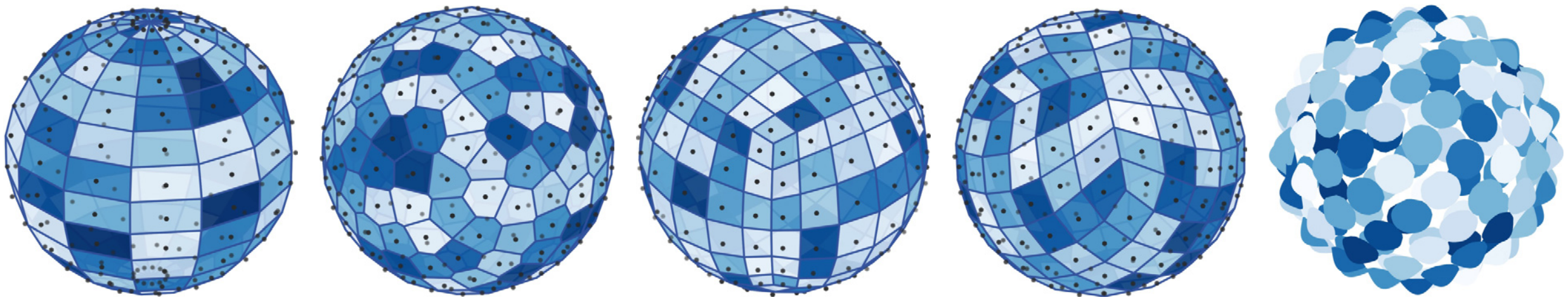

(a) (b) (c) (d) (e)

**Fig. 1.** Examples of discretisation schemes on the sphere, with an approximate resolution of 200 for each scheme. Figs. 1a to 1d show examples of non-uniform (Fig. 1a) and quasi-uniform (Figs. 1b to 1d) spherical point sets (marked by black dots). Fig. 1e on the other hand, shows an example of parametric discretisation by means of bell-shaped zonal basis functions. The equal-angle point set in Fig. 1a is obtained by gridding the azimuth-elevation domain. The point set in Fig. 1b is called the Fibonacci lattice, and can be generated as explained in [13]. The point sets in Figs. 1c and 1d are obtained from the cell centroids of the cubic and HEALPix spherical tessellations respectively. The cubic tessellation is obtained by projecting the pixelated faces of a cube onto the sphere. The HEALPix tessellation, very popular in cosmology and astronomy, is constructed by hierarchically subdividing the faces of a dodecahedron [10].

between fidelity to the data and compliance with this prior is assessed via a composite convex optimisation problem, linear combination of a *cost functional* and a *regularisation* term. Popular regularisation strategies include generalised Tikhonov (gTikhonov) or generalised total variation (gTV) [8], which favour physically admissible spherical fields with *smooth* and *sharp* variations respectively. Since spherical fields encountered in nature are continuous and hence have infinitely many degrees of freedom, scientists often constrain the approximation problem using *domain discretisation schemes*, which help reducing the number of degrees of freedom to something more manageable, ideally comparable to the size of the available data. For Euclidean domains, it is for example common practice to approximate the continuum by means of discrete *uniformly distributed point sets*, typically forming regular rectangular grids.[2] The popularity of such domain discretisation schemes can be primarily explained by their simplicity and computational conveniency. Indeed, signals defined over rectangular grids admit a natural representation as multi-dimensional arrays, a data structure commonly used in computer science for computation, storage and visualisation purposes.

Unfortunately, the sphere manifold structure makes it much more difficult to discretise by means of uniformly distributed point sets. For example, points gridded regularly on the azimuth-elevation domain $[0, 2\pi] \times [-\pi/2, \pi/2]$ are highly non-uniformly distributed at the surface of the sphere, with a much higher concentration of points near the poles (see Fig. 1a). As a matter of fact, uniform spherical point sets are only known [4, Chapter 3] for fixed numbers of points: 4, 6, 8, 12 and 20. They are respectively obtained from the vertices of the five Platonic solids: the *tetrahedron*, *cube*, *octahedron*, *dodecahedron* and *icosahedron*. For arbitrary numbers of points, spherical point sets with *quasi-uniform* distribution have been proposed [9,4,10]. The spherical Voronoi diagrams of the latter typically tile the sphere with near-regular polygonal tiles,[3] see for example Fig. 1. Unfortunately, quasi-uniform spherical point sets are significantly more complicated to work with as they are not easily represented by array-like data-structures. Moreover, derivatives and more generally pseudo-differential operators are difficult to approximate on quasi-uniform spherical point sets [11,12], making it cumbersome to work with gTikhonov and gTV priors.

The difficulty in designing domain discretisation schemes for the sphere has led scholars to consider alternative *parametric* discretisation schemes, where the unknown spherical field is constrained to a finite dimensional functional space, typically spanned by *zonal basis functions*[4] [15–19], i.e. functions with rotational invariance around a particular central direction on the sphere. The majority of zonal basis functions used in practice take the form of positive and *smooth* bell-shaped functions, sharply decaying to zero as the angular distance from their central direction increases (see Fig. 1e). They possess moreover many useful properties, particularly convenient for practical purposes:

- They are identical, spatially localised and highly symmetric, and thus easy to evaluate and amenable to sparse, parallel computations.
- Their overlapping supports and strong regularity make them well-suited for approximating smooth natural phenomena.
- Their centres can be positioned arbitrarily at the surface of the sphere, permitting for example the concentration of more zonal basis functions in regions more susceptible of welcoming the signal.
- *Strictly positive definite* zonal basis functions [18] are all linearly independent, irrespective of the chosen centres [20]. This guarantees a non-redundant representation and limits the risk of numerical instability.
- They are particularly well-suited for *scattered data interpolation problems* [15–17] where the spatial samples to interpolate may be non-uniformly distributed.

[2] This discretisation scheme is sometimes called *pixelisation* in visual computing.
[3] In practice, quasi-uniform spherical point sets are actually often constructed from centroids of spherical tessellation cells.
[4] Zonal basis functions are the spherical analogues of *radial basis functions* [14], used for scattered data interpolation in Euclidean settings.

This last fact is probably the main reason for their wide adoption in the literature. As a matter of fact, some zonal basis functions are not merely well-suited but *canonical* to spherical scattered data interpolation. This is notably the case for a specific type of zonal functions, called *spherical splines* [21], which arise naturally as solutions of interpolation problems on the sphere. As an illustration, consider the simplest interpolation problem where one wishes to find all *maximally smooth* functions with prescribed values $\{y_1, \dots, y_L\} \subset \mathbb{R}$ at directions $\{\boldsymbol{r}_1, \dots, \boldsymbol{r}_L\} \subset \mathbb{S}^2$. The relevant notion of smoothness is of course application dependent, but is generally enforced by seeking an interpolant with minimal generalised Tikhonov (gTikhonov) norm, induced by some linear, self-adjoint and strictly positive definite *pseudo-differential operator* $\mathscr{D}$. In mathematical terms, the prototypical interpolation problem can be formulated as:

$$\min_{f \in \mathscr{H}_{\mathscr{D}}} \|\mathscr{D} f\|_2 \quad \text{such that} \quad f(\boldsymbol{r}_i) = y_i,\ i = 1, \dots, L, \tag{1}$$

where the search space $\mathscr{H}_{\mathscr{D}}$ is an appropriately chosen reproducing kernel Hilbert space (RKHS) so that all the quantities involved in (1) are well-defined. It is then possible to show that there exists a *unique* maximally smooth interpolant $f^\star$ and that the latter has exactly $L$ degrees of freedom. Moreover, the maximally smooth interpolant $f^\star$ can be expressed [21, Section 6.3] as a spherical spline[5] with *knots* coinciding with the sampling directions $\{\boldsymbol{r}_1, \dots, \boldsymbol{r}_L\} \subset \mathbb{S}^2$. The $L$ spline weights can moreover be recovered by solving a square linear system [21, Section 6.3]. This result is quite remarkable, since it provides us with a *canonical* discretisation scheme operating in a *lossless* fashion: the infinite-dimensional optimisation problem (1) is transformed into an *equivalent* finite-dimensional optimisation problem, amenable to numerical optimisation. Theorem 6.40 in [21, Section 6.4.2] generalises[6] this result to the *smoothing spline approximation* problem:

$$\min_{f \in \mathscr{H}_{\mathscr{D}}} \sum_{i=1}^{L} |y_i - f(\boldsymbol{r}_i)|^2 + \lambda \|\mathscr{D} f\|_2^2, \tag{2}$$

used as an alternative to (1) in the context of noisy spatial samples, since it is less prone to overfitting.

Unfortunately, both problems (1) and (2) are too restrictive for most spherical approximation tasks encountered in practice. This is for example the case when the measurements are corrupted by non Gaussian noise – hence requiring a more general cost functional than a quadratic one – or do not consist in directional samples of the spherical field, but rather in local averages or more generally filtrations of the latter. In addition, gTikhonov-regularised optimisation problems *à la* (1) and (2) suffer from two main drawbacks: they tend to produce *overly smooth* interpolants and are too sensitive to the sampling locations $\{\boldsymbol{r}_1, \ldots, \boldsymbol{r}_L\} \subset \mathbb{S}^2$. This is a general behaviour of smooth spline approximation even in the Euclidean setting [8], which can be explained in part by the fact that the gTikhonov regularising norm is a weighted $\mathscr{L}^2$ norm. The latter favours indeed functions with relatively *smooth* variations, rendering the spline interpolant incapable of adapting to rapid changes in the data. To overcome this limitation scholars have, motivated by empirical studies [23], advocated the use of generalised total variation (gTV) regularisation norms, promoting functions with *relatively sparse* but potentially *sharp* variations, as often encountered in natural phenomena. However, representer theorems [24,8,25–27] characterising the form of the solutions yielded by the use of such regularisation strategies are, to date, unadapted to spherical geometries (see Section 1.3 for a discussion on prior art). For this reason, and in contrast with most fields of signal processing, total variation based penalties are still very much unexplored in spherical setups. As explained in greater detail in Section 1.2, one of the goals of this work is to close this theoretical gap and promote the wider use of such recovery methods across the community of geomathematicians.

[5] With respect to the operator $\mathscr{D}^2$.
[6] See also [22] for similar results in the specific case of the circle $\mathbb{S}^1$.

### *1.2. Contributions and organisation of the manuscript*

A primary goal of this paper is to offer a unified theoretical and practical approximation framework for gTV regularised infinite dimensional inference on the hypersphere $\mathbb{S}^{d-1}$ of arbitrary dimension $d \geq 2$. Out of concern for making the content of this work accessible to a wider audience, care has been taken in thoroughly interpreting and analysing the stated results and their assumptions, with a particular focus on their practical implications. Multiple real-life examples from environmental sciences and radio astronomy are moreover considered. The main contributions of this manuscript are listed hereafter and classified into three categories: *theory*, *algorithms* and *applications*.

*Theory.* Our main theoretical contribution is a *representer theorem* for functional penalised basis pursuit (FPBP) problems on the hypersphere, established and proved in Section 4. FPBP problems essentially seek spherical functions, measures or distributions minimising an optimal trade-off between a *convex* cost functional and a gTV penalty term. They are well suited for solving *functional linear inverse problems* on the hypersphere. To accommodate various measurements types, this class of optimisation problems is built upon a convenient *generalised sampling framework*, modelling the measurement process in terms of sampling linear functionals. Our representer Theorem 2 shows that the solution sets of FPBP problems can be characterised geometrically as the (weak*) *closed convex hull* of their extreme points. These extreme points take moreover the form of spherical splines with *sparse* innovations – i.e. fewer degrees of freedom than the total number of measurements – and unknown knots. The splines are canonically associated to the pseudo-differential operator $\mathscr{D}$ intervening in the gTV regularisation term. A key assumption of Theorem 2 is that the pseudo-differential operator $\mathscr{D}$ is *invertible*. With this assumption, Theorem 2 can be obtained as corollary of Lemma 1, an abstract representer theorem pertaining to penalised convex optimisation in abstract Banach spaces. This result, which generalises [28, Theorem 6] to the case of *non-strictly convex* cost functionals, is based on [28, Theorem 5], [8, Proposition 8] and [27, Theorem 3.1]. Note that assuming the regularising pseudo-differential operator $\mathscr{D}$ to be invertible does not restrict too much the applicability of Theorem 2 in practice. Indeed, we show that non-invertible pseudo-differential operators such as the *Laplace-Beltrami operator* $\Delta_{\mathbb{S}^{d-1}}$ can be brought into the scope of Theorem 2 if properly regularised on their nullspaces.

Finally, we provide a sufficient condition for sampling linear functionals to be *compatible* with a given regularising operator $\mathscr{D}$. More precisely, we show that choosing the sampling linear functionals in the predual of the search space yields well-defined FPBP problems. This allows us, in particular, to show that most practically interesting sampling functionals are compatible with the class of pseudo-differential operators considered in this paper, including atomic Dirac measures – convenient for spatial sampling – or square-integrable functions.

*Practical aspects & algorithms.* Adopting a more practical point of view, we outline in Section 5 the profound practical implications of the representer theorem for FPBP problems. First, we leverage Theorem 2 to derive a *canonical search space discretisation scheme* based on *quasi-uniform* $\mathscr{D}$*-splines*, i.e. splines whose knot sets form *quasi-uniform spherical point sets* [9]. This discretisation scheme is motivated by Proposition 10, which shows that, under mild conditions on the pseudo-differential operator $\mathscr{D}$, *quasi-uniform* $\mathscr{D}$-splines can approximate arbitrarily well most solutions of FPBP problems[7] when the number of knots tends to infinity. In Theorem 3 finally, we show that the finite dimensional problem resulting from this search space discretisation is a classical penalised basis pursuit (PBP) problem [29], which we propose to solve by means of provably convergent *fully-split proximal iterative methods* [30]. The latter involve only simple matrix-vector multiplications and *proximal steps*. We treat the most general case where the cost function is proximable but not necessarily differentiable with the *primal-dual splitting method (PDS)* introduced by Condat in his seminal work [31]. In the simpler (yet prevailing in practice) case where the cost functional is also *differentiable* and with *Lipschitz continuous* derivative, we leverage an optimal first-order method called *accelerated proximal gradient descent (APGD)* [32,30], with faster convergence rate than the PDS method. We conclude by discussing suitable choices of pseudo-differential and spherical splines for practical purposes. We introduce *Wendland* and *Matérn* pseudo-differential operators, whose Green kernels have simple closed-form expressions and good localisation in space. The latter have moreover spectra equivalent to those of *Sobolev operators* $\mathscr{D}_\beta = [\mathrm{Id} - \Delta_{\mathbb{S}^{d-1}}]^\beta$, $\beta > (d-1)/2$, often used in functional analysis.

*Applications.* To demonstrate the versatility of our proposed sparse spline approximation framework, we put it to the test in Section 6 on two datasets originating from real-life spherical approximation problems encountered in environmental sciences. First, we reconstruct a global map of sea surface temperature anomalies from recordings collected by drifting floats of the ARGO fleet [33,34]. Such maps are used in environmental sciences to monitor global climate change as well as manage the population of marine species and ecosystems particularly sensitive to fluctuations in the water temperature. We compare our reconstruction strategy to the state-of-the-art smoothing spline method with Tikhonov regularisation and observe that total variation penalties yield indeed higher resolution maps. In a second example, we build global density maps of wildfires across the globe for the year 2016, using fire counts recorded by NASA's Aqua and Terra satellites. Wildfires maps allow scientists to better understand the chemistry of the atmosphere and its impact on climate. Because of the Poisson-like distribution of count data, we investigate the use of a *KL-divergence* as data-fidelity term. In both cases, interactive versions of the estimated spherical maps are available online at the links https://matthieumeo.github.io/temperature_anomalies.html and https://matthieumeo.github.io/fire_density.html respectively.

### 1.3. Representer theorems in the literature

In this section, we review the most notable representer theorems proposed in the literature for gTV regularisation, and discuss their limitations in the context of spherical approximation. A summary of this section is provided in Table 1.

Inspired by the pioneering work of Fisher and Jerome in [35], Unser et al. have investigated in a series of papers [24,8,36–38] gTV-regularised functional inverse problems over $\mathbb{R}^d$. Like Theorem 2, they consider a generalised sampling framework, compatible with a great variety of linear measurements. Unlike Theorem 2 however, they do not assume that the regularising pseudo-differential operator is invertible, and allow it to have a finite dimensional nullspace. The conclusions of their representer theorem in [24] – derived for Euclidean domains only – are analogous to the ones of Theorem 2, proposed in this paper for spherical domains. In subsequent publications [8,36,38], the authors proposed canonical discretisation schemes as well as numerical algorithms for approximating extreme point solutions of gTV-penalised convex optimisation problems. In [8], they consider a discretisation based on cardinal splines of the gTV pseudo-differential operator, with uniform knots chosen over a dense grid. In [36], they propose a numerically stabler discretisation scheme, based this time on *multi-resolution B-splines* with refinable grid sizes. In both cases they solve the resulting discrete optimisation problem with a two-stage procedure leveraging *proximal gradient descent* and the *simplex* algorithm. In the specific case where the domain is $\mathbb{R}$ and the regularising operator is the second derivative, Debarre et al. describe in [38] an even simpler reconstruction algorithm. In contrast with the proximal algorithms proposed in this paper, the optimisation pipelines proposed in these papers are however limited to differentiable cost functionals with Lipschitz-continuous derivatives. While remarkably generic, their spline approximation framework is only valid for functions defined over $\mathbb{R}^d$. An extension of the theoretical framework to the hypertorus $\mathbb{T}^d = \mathbb{R}^d/2\pi\mathbb{Z}^d$ has been recently proposed [39], but a similar work for the spherical setting is still missing.

[7] Of course with gTV regularisation norms induced by $\mathscr{D}$.

**Table 1**
Comparison between existing representer theorems and the one established in this paper.

| | Domain | Existence | Generalised sampling | Cost functional | Charac. of solutions | Search space | Charac. of predual |
|---|---|---|---|---|---|---|---|
| Boyer et al. [27] Bredies et al. [26] | Agnostic | ✗ | ✓ | (quasi-) Convex | Convex-hull of Extreme Points | Banach | ✗ |
| Unser [28] | Agnostic | ✓ | ✓ | Strictly Convex | Duality Map | Banach | ✗ |
| Unser, Fageot et al. [24,8,37,39] | $\mathbb{R}^d, \mathbb{T}^d$ | ✓ | ✓ | Convex | Convex-hull of Extreme Points | Maximal | ✓ |
| This paper | $\mathbb{S}^{d-1}$ | ✓ | ✓ | Convex | Convex-hull of Extreme Points | Maximal | ✓ |

In a subsequent work [25], Flinth et al. proposed an alternative proof of the representer theorem proposed in [24]. Their proof is based on a limit argument, considering nested finite dimensional discretisations of the domain $\Omega \subset \mathbb{R}^d$ based on finer and finer uniform rectangular grids. They claim that such an approach, presented in the Euclidean case for the sake of simplicity, could easily be adapted to domains more general than $\mathbb{R}^d$, such as the torus or any separable, locally compact topological space. They specify however that such an extension would require modifying adequately the discretisation scheme to the specific geometry of the domain, without giving additional details on how this could be achieved canonically. Unfortunately, such a task may be very complex if even possible at all for geometries such as the sphere. Indeed, discretising the sphere by means of nested quasi-uniform point sets with finer and finer resolution is, as previously discussed, a nontrivial problem.

More recently, Boyer et al. [27] and Bredies et al. [26] have independently shown that the solutions to infinite dimensional optimisation problems with convex regularisers are convex combinations of extreme points of the regulariser level sets. This result applies notably to gTV regularisers with not only scalar but also *vector* pseudo-differential operators such as the gradient. This is in contrast with the previously cited works which were all limited to scalar pseudo-differential operators such as the Laplace-Beltrami operator. While theoretically applicable to spherical geometries, their result neither addresses existence conditions nor characterises the minimal search space (and its corresponding predual) associated to a certain gTV norm. This is problematic for practical purposes, where it is crucial to know if a given optimisation problem admits a solution or understand which sampling linear functionals are compatible with a specific choice of gTV penalty.

Finally, Unser [28] established a Banach representer theorem with very broad applicability. Unlike the previously cited results, this representer theorem relies on the notion of *duality map*, which generalises the Hilbert notion of *Riesz map* to Banach spaces. More precisely, it shows that the solutions of convex regularised inverse problems are contained in the image by a certain duality map of a linear combination of the sensing linear functionals. As acknowledged by the author, this result is however of limited use in the context of gTV regularisation, since the duality map is unknown, nonlinear and set-valued.

### *1.4. Notations and terminology*

Throughout the manuscript, we adopt the following conventions:

- We use the term *spherical field* to refer, depending on the context, to *functions*, *measures* or *generalised functions* [40] defined over the sphere $\mathbb{S}^{d-1}$ for any dimension[8] $d \geq 2$. In full generality, one shall think at a spherical field as an element of some *infinite-dimensional* Banach space $f \in \mathscr{B}$.
- It is traditional to call the 1-sphere $\mathbb{S}^1 \subset \mathbb{R}^2$ a *circle*, the 2-sphere $\mathbb{S}^2 \subset \mathbb{R}^3$ a *sphere* and the $(d-1)$-sphere $\mathbb{S}^{d-1} \subset \mathbb{R}^d$, $d \geq 2$ a *hypersphere*. For the sake of simplicity, we break with tradition and use the appellation "sphere" agnostic to the underlying dimension. Moreover, we denote by $\mathfrak{a}_d$ the *area* of the unit sphere $\mathbb{S}^{d-1}$, $d \geq 2$, given in general by: $\mathfrak{a}_d = 2\pi^{d/2}/\Gamma(d/2)$. We have notably $\mathfrak{a}_2 = 2\pi$ and $\mathfrak{a}_3 = 4\pi$.
- Vectors and matrices are written in bold face, in an attempt to make finite-dimensional quantities more apparent. The *adjoint*, *Moore-Penrose pseudo-inverse*, *range* and *nullspace* of a linear operator $\Phi$ are denoted by $\Phi^*$, $\Phi^\dagger$, $\mathcal{R}(\Phi)$ and $\mathcal{N}(\Phi)$ respectively. In the specific case of real and complex matrices $\boldsymbol{A}$, we denote the adjoint by $\boldsymbol{A}^T$ and $\boldsymbol{A}^H$, respectively. For scalars $z \in \mathbb{C}$ finally, we denote by $\bar{z}, |z|, \Re(z), \Im(z)$ the *conjugate*, *modulus*, *real part* and *imaginary part* of $z$ respectively.
- For a real functional $F : \Theta \to \mathbb{R}$, $\operatorname{Arg\,min}_{\theta\in\Theta} F(\theta)$ denotes the set $\{\hat{\theta} \in \Theta : F(\hat{\theta}) = \min_{\theta\in\Theta} F(\theta)\}$ of all minimisers of $F$, while $\arg\min_{\theta\in\Theta} F(\theta)$ denotes an element $\hat{\theta} \in \operatorname{Arg\,min}_{\theta\in\Theta} F(\theta)$ of this set. When there is a unique minimiser, we have in particular: $\operatorname{Arg\,min}_{\theta\in\Theta} F(\theta) = \{\arg\min_{\theta\in\Theta} F(\theta)\} = \{\hat{\theta}\}$.

## 2. Preliminaries

This section briefly reviews the key ingredients of our approximation framework, borrowed from the fields of functional analysis, measure theory and Fourier analysis on the hypersphere respectively. We begin our discussion in Section 2.1 with the concept of *duality* in topological vector spaces, central to the process of generalised sampling introduced in Section 4.1. In Section 2.1.3, we define and provide a dual characterisation of the *total variation norm* for regular Borel measures. The latter generalises the discrete $\ell_1$ norm to continuous setups, and will be used in our approximation framework as a sparsity-promoting regularisation norm. Finally, we conclude in Section 2.2 with important notions from Fourier analysis on the hypersphere, namely *spherical harmonics* and *spherical zonal functions*, which will help us define the class of *spline-admissible pseudo-differential operators* in Section 3. Readers already familiar with all these notions can skip this section and continue their read from Section 3. The main reference for the topics discussed in Sections 2.1 and 2.1.3 is [40, Chapter 3]. For Section 2.2, we refer to the thesis dissertation [15] and the book chapter [21, Chapter 5]. Additional results and theorems relating to functional analysis and Fourier analysis on the hypersphere are also available in [41, Chapter 1, 3 and 4].

### *2.1. Duality in topological vector spaces*

In order to approximate a spherical field, modelled here as a generic element $f$ of a vector space $\mathscr{B}$, one must first collect evidence of the latter, often by sensing it via a *linear acquisition device*. As we shall see in Section 4.1, these linear measurements can in general be modelled as the outcomes of a collection of device-specific *linear functionals* acting on the object $f$ of interest. In this section, we investigate the structure of the space of all linear functionals associated with a given vector space $\mathscr{B}$, with a special focus on those that yield *well-defined* measurements when acting on any element $f \in \mathscr{B}$.

[8] Of course the cases $d = 2, 3$ will be particularly prevailing in real-life applications.

### 2.1.1. Schwartz duality product

For a vector space $\mathscr{B}$ over a scalar field $\mathbb{C}$, the space of all linear functionals $f : \mathscr{B} \to \mathbb{C}$ is a vector space called the *algebraic dual* and is denoted by $\mathscr{B}^*$. It is customary to write the action of a linear functional $f \in \mathscr{B}^*$ onto an element $h \in \mathscr{B}$ by means of a *bilinear map* $\langle \cdot | \cdot \rangle : \mathscr{B}^* \times \mathscr{B} \to \mathbb{C}$ called the *Schwartz duality product* defined as

$$\langle \cdot | \cdot \rangle : \begin{cases} \mathscr{B}^* \times \mathscr{B} \to \mathbb{C}, \\ (f, h) \mapsto \langle f | h \rangle := f(h). \end{cases} \tag{3}$$

The *bra-ket* notation used in (3) to denote the Schwartz duality product is common in quantum mechanics and was introduced by Paul Dirac in 1939. Its resemblance to the inner product is not fortuitous, and is motivated by Hilbert space theory. Indeed, the *Riesz-Fréchet representation theorem* [42] states that, for a Hilbert space $\mathscr{H}$, every linear functional in $\mathscr{H}^*$ can be written as an inner product with some unique element $g$ of $\mathscr{H}$, i.e. $\forall f \in \mathscr{H}^*$, $\exists! g \in \mathscr{H}$ such that $\langle f | h \rangle = \langle h, g \rangle_{\mathscr{H}}$.

### 2.1.2. Topological dual

Any sensible acquisition system should react continuously to variations in its input. It seems hence reasonable to require that the linear functionals modelling it be *continuous* as well. The subset of continuous linear functionals in the algebraic dual is a linear subspace, called the *topological dual*[9] and denoted by $\mathscr{B}'$. In infinite dimensions, not all linear functionals are guaranteed to be continuous so we have in general $\mathscr{B}' \subset \mathscr{B}^*$. Moreover, since continuity is a topological notion, there can exist multiple topological duals of a given space $\mathscr{B}$ depending on the topology chosen on the latter. In the special case where $(\mathscr{B}, \|\cdot\|)$ is a Banach space equipped with its canonical normed topology, continuous linear functionals can be characterised as:

$$\mathscr{B}' = \left\{ f \in \mathscr{B}^* : |||f||| := \sup_{h \in \mathscr{B}, \|h\|=1} |\langle f | h \rangle| < \infty \right\}. \tag{4}$$

The norm $|||\cdot||| : \mathscr{B}' \to \mathbb{R}_+$ used in (4) makes $\mathscr{B}'$ a *Banach space* and is called the *dual norm* induced by the norm $\|\cdot\|$ on $\mathscr{B}$. In plain words, (4) states that continuous linear functionals are *bounded*[10] linear functionals in $\mathscr{B}^*$. The latter will hence produce bounded measurements, well-defined for any input $f \in \mathscr{B}$. The Banach topology induced by the dual norm defined in (4) is called the *strong topology*. The topological dual $\mathscr{B}'$ can also be equipped with the weak* topology, which is the *coarsest* topology on $\mathscr{B}'$ such that elements $h \in \mathscr{B}$ are continuous functionals on $\mathscr{B}'$. This is the topology of *pointwise* convergence: $\{f_n,\, n \in \mathbb{N}\} \subset \mathscr{B}'$ converges towards a limit functional $f^* \in \mathscr{B}'$ with respect to the weak* topology i.f.f. $\lim_{n\to\infty} |\langle f^* - f_n | h \rangle| = 0$, $\forall h \in \mathscr{B}$. Throughout, we will employ expressions such as "weak* compact", "weak* closed" or "weak* convergent" when it is important to make obvious the underlying topology with respect to which the topological notions should be understood.

**Vocabulary 1** *(Predual and duality pair).* A topological vector space $\mathscr{B}$ and its topological dual $\mathscr{B}'$ are said to form a *duality pair*. Moreover, $\mathscr{B}$ is called the *predual* of $\mathscr{B}'$.

---

[9] In the manuscript, the shorthand expression "dual space" is sometimes used to refer to the topological dual space.

[10] Note that in finite dimensions things are much less complicated since every linear functional is bounded and all norms are topologically equivalent, hence $\mathscr{B}^* = \mathscr{B}'$.

#### *2.1.3. Duality pairs for common functional spaces*

In this section we provide well-known duality pairs for functional spaces of interest.

*Schwartz functions and generalised functions.* The space of *generalised functions* or *distributions* is almost the largest functional space that can be defined on $\mathbb{S}^{d-1}$. It contains as subspaces the Lebesgue spaces as well as the spaces of continuous functions or regular Borel measures. It is denoted $\mathscr{S}'(\mathbb{S}^{d-1})$ and is defined as the topological dual of the space of *Schwartz functions*[11] $\mathscr{S}(\mathbb{S}^{d-1}) = \mathscr{C}^\infty(\mathbb{S}^{d-1})$, equipped with the *metric* topology generated by the family of norms:

$$\|\cdot\|_{n,\infty} : \begin{cases} \mathscr{C}^\infty(\mathbb{S}^{d-1}) \to \mathbb{R}_+ \\ h \mapsto \|h\|_{n,\infty} := \|(\mathrm{Id} - \Delta_{\mathbb{S}^{d-1}})^n h\|_\infty \end{cases} \qquad \forall n \in \mathbb{N},$$

where $\Delta_{\mathbb{S}^{d-1}}$ is the *Laplace-Beltrami* operator on $\mathbb{S}^{d-1}$. The Schwartz space is a *locally convex Fréchet space*.[12] It is *not* normable and in particular not complete with the supremum norm (indeed, it is dense in the space of continuous functions). Note that since the metric topology on $\mathscr{S}(\mathbb{S}^{d-1})$ is not induced by a norm,[13] it is not possible to define a strong topology on $\mathscr{S}'(\mathbb{S}^{d-1})$. We will hence always assume the weak* topology as canonical topology on $\mathscr{S}'(\mathbb{S}^{d-1})$.

*Continuous functions, measures and total variation norm.* The *Riesz-Markov representation theorem* (see [41, Theorem 2.5] for a statement of the latter in the spherical setup) establishes the duality pair

$$(\mathscr{C}(\mathbb{S}^{d-1}), \|\cdot\|_\infty)' \cong (\mathcal{M}(\mathbb{S}^{d-1}), \|\cdot\|_{TV})$$

between the space of continuous functions equipped with the supremum norm and the space of regular Borel measures equipped with the total variation (TV) norm [41, Definition 2.2]. The TV norm can be thought as an $\mathscr{L}^1$ norm for measures,[14] which motivates its use as sparsity-inducing regularisation norm in penalised convex optimisation problem involving measures as considered in this paper. From (4) and the Riesz-Markov representation theorem, we moreover get the following dual characterisation of $\mathcal{M}(\mathbb{S}^{d-1})$:

$$\mathcal{M}(\mathbb{S}^{d-1}) \cong \left\{ f \in \mathscr{S}'(\mathbb{S}^{d-1}) : \ \|f\|_{TV} = \sup_{h \in \mathscr{S}(\mathbb{S}^{d-1}),\, \|h\|_\infty = 1} |\langle f | h \rangle| < +\infty \right\}, \tag{5}$$

where we have used the density of Schwartz functions in the space of bounded continuous functions of [24, Section 3]. Equation (5) permits us to see $\mathcal{M}(\mathbb{S}^{d-1})$ as the subspace of generalised functions with finite dual norm.

### *2.2. Fourier analysis on the hypersphere*

The class of spherical pseudo-differential operators introduced in Section 3 are defined implicitly in the Fourier domain. In this section, we hence introduce the basic mathematical machinery needed for performing Fourier analysis on the hypersphere. The material presented in this chapter is based on the formalism adopted in [21,4,15]. Useful additional results are also available in [41, Chapter 3].

[11] The hypersphere being compact and bounded, Schwartz functions simply reduce to infinitely smooth functions.
[12] Fréchet spaces generalise Banach spaces with metrics that do not originate from norms.
[13] It is possible to show that $\mathscr{S}(\mathbb{S}^{d-1})$ is not normable [40].
[14] As a matter of fact, the TV norm of a measure absolutely continuous w.r.t. the Lebesgue measure is given by the $\mathscr{L}^1$ norm of its density (or Radon-Nikodym derivative) [41, Remark 2.1].

#### *2.2.1. Spherical harmonics*

One possible route towards defining the Fourier basis on the hypersphere is to proceed analogously to Fourier and study fundamental solutions of the *heat differential equation* on $\mathbb{S}^{d-1}$. The separation of variables technique reveals that the spherical component of such fundamental solutions are eigenfunctions of the *Laplace-Beltrami operator*[15] on $\mathbb{S}^{d-1}$ [21, Chapter 5]. They are called *spherical harmonics*.

**Definition 1** *(Spherical harmonics).* Let $\Delta_{\mathbb{S}^{d-1}}$ be the *Laplace-Beltrami operator* on $\mathbb{S}^{d-1}$ with *spectrum* $\{\lambda_n = -n(n+d-2),\, n \in \mathbb{N}\}$. We call *spherical harmonic* of order $n$ any *eigenfunction* $Y$ in the *eigenspace* $\mathrm{Harm}_n(\mathbb{S}^{d-1})$ associated to the *eigenvalue* $\lambda_n$:

$$\mathrm{Harm}_n(\mathbb{S}^{d-1}) := \left\{ Y : \mathbb{S}^{d-1} \to \mathbb{C} \,|\, \Delta_{\mathbb{S}^{d-1}} Y = -n(n+d-2)Y \right\}.$$

Moreover, we denote by $\mathfrak{B}_n := \{Y_n^m,\ m = 1, \cdots, N_d(n)\}$ any *orthonormal basis* of $\mathrm{Harm}_n(\mathbb{S}^{d-1})$, where $N_d(n)$ is the *geometric multiplicity* of the eigenvalue $\lambda_n$.

**Remark 1** *(Geometric multiplicity).* The geometric multiplicity $N_d(n)$ of each eigenspace $\mathrm{Harm}_n(\mathbb{S}^{d-1})$ can be computed explicitly [15, Chapter 2]. It is given in general by:

$$N_d(0) = 1, \qquad \& \qquad N_d(n) = \frac{2n+d-2}{n}\binom{n+d-3}{n-1}, \qquad n \geq 1.$$

In particular, for $d = 2, 3$, we get $N_2(n) = 2$ and $N_3(n) = 2n + 1$, $n \geq 1$. We have moreover the asymptotic behaviour [15, Chapter 2]:

$$N_d(n) = \mathcal{O}\left(n^{d-2}\right). \tag{6}$$

It is possible to show [21, Chapter 5] that the direct summation of the eigenspaces $\mathrm{Harm}_n(\mathbb{S}^{d-1})$ is *asymptotically dense* in $\mathscr{L}^2(\mathbb{S}^{d-1})$: $\mathscr{L}^2(\mathbb{S}^{d-1}) = \bigoplus_{n=0}^{+\infty} \mathrm{Harm}_n(\mathbb{S}^{d-1})$. This yields the Fourier expansion theorem on the sphere:

**Theorem 1** *(Spherical Fourier expansion [21]). Let $d \geq 2$, $n \in \mathbb{N}$ and $\mathfrak{B}_n = \{Y_n^m, m = 1, \ldots, N_d(n)\}$ be an* orthonormal basis *of $\mathrm{Harm}_n(\mathbb{S}^{d-1})$. Then, every function $f \in \mathscr{L}^2(\mathbb{S}^{d-1})$ admits a* spherical Fourier expansion *given by*

$$f \overset{\mathscr{L}^2}{=} \sum_{n=0}^{+\infty} \sum_{m=1}^{N_d(n)} \hat{f}_n^m \, Y_n^m,$$

*where the* spherical Fourier coefficients $\{\hat{f}_n^m\} \subset \mathbb{C}$ *of $f$ are given by the spherical harmonic transform (SHT):*

$$\hat{f}_n^m = \langle f, Y_n^m \rangle_{\mathbb{S}^{d-1}} = \int\limits_{\mathbb{S}^{d-1}} f(\boldsymbol{r}) \overline{Y_n^m(\boldsymbol{r})} \, d\boldsymbol{r}, \qquad n \in \mathbb{N},\, m = 0, \ldots, N_d(n).$$

**Remark 2** *(Fully normalised spherical harmonics).* Note that the spherical harmonics $Y_n^m$ in Definition 1 and Theorem 1 are not uniquely specified, since there exist infinitely many orthonormal bases $\mathfrak{B}_n$ of $\mathrm{Harm}_n(\mathbb{S}^{d-1})$. In practice, the convention is to work with the system of so-called fully normalised spherical harmonics (FNSH),[16] obtained inductively by the method of separation of variables applied to the eigenvalue problem for the Laplace-Beltrami operator $\Delta_{\mathbb{S}^{d-1}}$. This construction is detailed in [21, Section 5.2] for the case $d = 2$ and in [43] for the general case. In all that follows, we will always assume the fully normalised spherical harmonics as canonical Fourier basis. Closed-form analytical expressions of the fully normalised spherical harmonics for $d = 2, 3$ are provided in [41, Example 3.1]. Formulae for the more general case $d > 3$ are available in [43].

[15] The *Laplace-Beltrami operator* on $\mathbb{S}^{d-1}$ generalises the *Laplace operator* $\Delta_{\mathbb{R}^d}$ in $\mathbb{R}^d$ to the manifold setting. Both operators are linked by the relationship: $\Delta_{\mathbb{R}^d} = \frac{\partial^2}{\partial\rho^2} + \frac{d-1}{\rho}\frac{\partial}{\partial\rho} + \frac{1}{\rho^2}\Delta_{\mathbb{S}^{d-1}}$, where, for every $\boldsymbol{x} \in \mathbb{R}^d \backslash \{\boldsymbol{0}\}$, we define $\boldsymbol{x} := \rho\boldsymbol{r}$ with $\rho := \|\boldsymbol{x}\| \in \mathbb{R}_+$ and $\boldsymbol{r} \in \mathbb{S}^{d-1}$.

Using the bilinearity of the Schwartz duality product, it is possible to extend the SHT to generalised functions [41, Remark 3.3]. The generalised spherical harmonic transform (gSHT) of a generalised function $f \in \mathscr{S}'(\mathbb{S}^{d-1})$ is defined as

$$f \overset{\text{weak}^*}{=} \sum_{n=0}^{+\infty} \sum_{m=1}^{N_d(n)} \hat{f}_n^m Y_n^m, \quad \text{where} \quad \hat{f}_n^m := \langle f | Y_n^m \rangle, \quad n \in \mathbb{N}, m = 0, \ldots, N_d(n). \tag{7}$$

Notice that the convergence of the infinite series in (7) is w.r.t. the weak* topology (see Section 2.1.2). Since the spherical harmonics are infinitely differentiable and hence in the predual $\mathscr{S}(\mathbb{S}^{d-1})$ of $\mathscr{S}'(\mathbb{S}^{d-1})$, the Fourier coefficients $\{\hat{f}_n^m\}$ are moreover well-defined.

### *2.2.2. Spherical zonal kernels*

Another route towards defining the Fourier basis consists of looking at eigenfunctions of linear shift invariant systems, or *convolution operators* [44, Chapters 3 and 4]. As the hypersphere is a manifold, there is no intrinsically defined notion of convolution.[17] It is however possible to define a class of linear integral operators which "behave" as traditional convolution operators. In the Euclidean setting, convolution operators have *shift-invariant* kernels, whose value at a pair $(\boldsymbol{r}, \boldsymbol{s}) \in \mathbb{R}^d \times \mathbb{R}^d$ depends only on the distance $\|\boldsymbol{r} - \boldsymbol{s}\|_{\mathbb{R}^d}$ between the two points. We can extend this notion to the hypersphere by noticing that the *chord distance* between two input directions $\boldsymbol{r}, \boldsymbol{s} \in \mathbb{S}^{d-1}$ is given by: $\|\boldsymbol{r} - \boldsymbol{s}\|_{\mathbb{R}^d} = \sqrt{\|\boldsymbol{r}\|^2 + \|\boldsymbol{s}\|^2 - 2\langle \boldsymbol{r}, \boldsymbol{s}\rangle} = \sqrt{2 - 2\langle \boldsymbol{r}, \boldsymbol{s}\rangle}$. Notice that this quantity depends only on the *inner product* between the two directions $\boldsymbol{r}, \boldsymbol{s}$. This observation naturally leads to the notion of *zonal kernel*:

**Definition 2** *(Spherical zonal kernel).* A kernel $\Psi : \mathbb{S}^{d-1} \times \mathbb{S}^{d-1} \to \mathbb{C}$ is called a *spherical zonal kernel* if there exists a function $\psi : [-1, 1] \to \mathbb{C}$ such that: $\Psi(\boldsymbol{r}, \boldsymbol{s}) = \psi(\langle \boldsymbol{r}, \boldsymbol{s}\rangle)$, $\forall (\boldsymbol{r}, \boldsymbol{s}) \in \mathbb{S}^{d-1} \times \mathbb{S}^{d-1}$. For brevity, the function $\psi$ is often abusively referred to as the *zonal kernel* and no reference is made to $\Psi$. Moreover, for every $\boldsymbol{s} \in \mathbb{S}^{d-1}$, we call the trace $\psi(\langle \cdot, \boldsymbol{s}\rangle) : \mathbb{S}^{d-1} \to \mathbb{C}$ of a zonal kernel a *zonal function*.

Zonal kernels can then be used to construct *spherical convolution operators* [21]:

**Definition 3** *(Spherical convolution operator).* Let $\psi \in \mathscr{L}^2([-1, 1])$ be a *zonal kernel*. The *spherical convolution operator* $\mathscr{I}_\psi : \mathscr{L}^2(\mathbb{S}^{d-1}) \to \mathscr{L}^2(\mathbb{S}^{d-1})$ is defined as

$$\mathscr{I}_\psi : \begin{cases} \mathscr{L}^2(\mathbb{S}^{d-1}) \to \mathscr{L}^2(\mathbb{S}^{d-1}) \\ f \mapsto \{\psi * f\}(\boldsymbol{r}) = \int\limits_{\mathbb{S}^{d-1}} \psi(\langle \boldsymbol{r}, \boldsymbol{s}\rangle) f(\boldsymbol{s})\, d\boldsymbol{s}, \quad \forall \boldsymbol{r} \in \mathbb{S}^{d-1}. \end{cases} \tag{8}$$

**Remark 3.** It is shown in [21, Theorem 7.2] that, under the assumptions of Definition 3, the image of the convolution operator (8) is indeed $\mathscr{L}^2(\mathbb{S}^{d-1})$.

---

[16] The terminology "fully normalised" is slightly deceptive. Indeed, all $Y_n^m$ are normalised, independently of the orthonormal system $\mathfrak{B}_n$ chosen on $\text{Harm}_n(\mathbb{S}^{d-1})$.

[17] See Chapter 1 of [4] for an in-depth discussion on the topic.

Using two important results from spherical Fourier analysis, namely the *addition theorem* [21, Theorem 5.11] and the *Funk-Hecke formula* [21, Theorem 7.3], it can be shown [41, Proposition 3.8] that spherical convolution operators defined from zonal kernels are indeed diagonalised by spherical harmonics.

## 3. Hyperspherical splines

In this section we introduce *hyperspherical splines* – or *spherical splines* for short, which play a central role in spherical approximation theory [21, Chapter 6]. To this end, we extend the approach of [40, Chapter 6] to the spherical setting and construct spherical splines as "primitives" of finite Dirac streams w.r.t. a certain class of *pseudo-differential operators*, called *spline-admissible*. In short, spline-admissible operators are such that their fundamental solutions, called *Green functions*, are *ordinary* functions.[18] We derive a sufficient condition for spline-admissibility and provide examples of spline-admissible operators among the pseudo-differential operators most commonly used in practice.

### *3.1. Spherical pseudo-differential operators*

By analogy with the Euclidean case, we define spherical pseudo-differential operators as *Fourier multipliers* with *slowly growing* spectra.

**Definition 4** *(Spherical pseudo-differential operator).* We call *spherical pseudo-differential operator* any *linear* operator of the form

$$\mathscr{D}: \begin{cases} \mathscr{S}(\mathbb{S}^{d-1}) \to \mathscr{S}(\mathbb{S}^{d-1}) \\ h \mapsto \mathscr{D}h := \sum_{n=0}^{+\infty} \hat{D}_n \left[ \sum_{m=1}^{N_d(n)} \hat{h}_n^m Y_n^m \right], \end{cases} \tag{9}$$

where $\{\hat{h}_n^m,\, n \in \mathbb{N},\, m = 1, \dots, N_d(n)\}$ are the *spherical Fourier coefficients* of $h$ and $\{\hat{D}_n\}_{n\in\mathbb{N}} \in \mathbb{R}^{\mathbb{N}}$ is a sequence of *real numbers* – called the *Fourier symbol* of $\mathscr{D}$ – such that the set

$$\mathfrak{K}_{\mathscr{D}} := \left\{ n \in \mathbb{N} : \ |\hat{D}_n| = 0 \right\}, \tag{10}$$

is *finite*, i.e. $\#\mathfrak{K}_{\mathscr{D}} := N_0 < +\infty$, and

$$|\hat{D}_n| = \Theta\left(n^p\right), \tag{11}$$

for some real number $p \geq 0$, called the *spectral growth order* of $\mathscr{D}$.

**Remark 4** *(Theta notation).* The condition $|\hat{D}_n| = \Theta\left(n^p\right)$ for some $p \geq 0$ means that $|\hat{D}_n| = \mathcal{O}\left(n^p\right)$ and $|\hat{D}_n| = \Omega\left(n^p\right)$, i.e. there exists $n_0 \in \mathbb{N}$ such that $\forall n \geq n_0$ we have $C_1 n^p \leq |\hat{D}_n| \leq C_2 n^p$, for some positive constants $C_1, C_2 \in \mathbb{R}_+$. In other words, the sequence $\{|\hat{D}_n|\}_{n\in\mathbb{N}}$ is asymptotically comparable to the polynomial $n^p$.

Notice that $\mathscr{D}$ multiplies the Fourier coefficients of its argument $h$ by a sequence $\{\hat{D}_n\}_{n\in\mathbb{N}}$ with polynomial growth order. This filtering operation effectively boosts the high frequency content of $h$, hence making it "rougher" (less regular). This behaviour is reminiscent of the one of the Laplace-Beltrami operator, which can indeed be shown to verify Definition 4. Examples of common pseudo-differential operators are provided in the following example.

[18] Ordinary functions are functions which are everywhere defined pointwise.

**Example 1** *(Common pseudo-differential operators).* Consider the Laplace-Beltrami operator $\Delta_{\mathbb{S}^{d-1}}$.

- **Iterated Laplace-Beltrami operators:** these operators are obtained as *integer powers* (i.e. successive compositions) of the Laplace-Beltrami operator $\mathscr{D} := \Delta^k_{\mathbb{S}^{d-1}}$, with $k \in \mathbb{N}$. They are indeed pseudo-differential operators since they can be written[19] as in (9) with $\hat{D}_n = (-n(n+d-2))^k$, $n \in \mathbb{N}$. We have indeed $\hat{D}_n \in \mathbb{R}$ (since $k \in \mathbb{N}$), $|\hat{D}_n| = \Theta(n^{2k})$, and $\mathfrak{K}_{\mathscr{D}} = \{0\}$ is finite. Notice that $\Delta^k_{\mathbb{S}^{d-1}}$ is positive semi-definite for $k$ even and negative semi-definite for $k$ odd.
- **Fractional Laplace-Beltrami operators:** these operators are obtained as $p$-th *roots* of the negative Laplace-Beltrami operator $\mathscr{D} := (-\Delta_{\mathbb{S}^{d-1}})^{1/p}$, with $p \in \mathbb{N}^*$. They are indeed pseudo-differential operators since they can be written as in (9) with $\hat{D}_n = \sqrt[p]{n(n+d-2)}$, $n \in \mathbb{N}$. Since $n(n+d-2) > 0$, $\forall n \geq 1$, $d \geq 2$, we have indeed $\hat{D}_n \in \mathbb{R}$. Moreover we also have $|\hat{D}_n| = \Theta(n^{2/p})$, and $\mathfrak{K}_{\mathscr{D}} = \{0\}$ is finite. Notice that $(-\Delta_{\mathbb{S}^{d-1}})^{1/p}$ is always positive semi-definite. The case $p = 2$ yields the square-root of the Laplace-Beltrami operator, which is intimately linked to the spherical *gradient* differential operators $\boldsymbol{\nabla}_{\mathbb{S}^{d-1}}$. The latter is however *vector-valued*, and hence does not belong to the class of pseudo-differential operators considered in Definition 4.
- **Sobolev operators:** these operators are defined as $\mathscr{D} := (\mathrm{Id} - \Delta_{\mathbb{S}^{d-1}})^{\beta}$, with $\beta > 0$. They are indeed pseudo-differential operators since their Fourier symbols are given by $\hat{D}_n = (1 + n(n+d-2))^{\beta}$, $n \in \mathbb{N}$, and hence $\hat{D}_n \in \mathbb{R}$, $|\hat{D}_n| = \Theta(n^{2\beta})$, and $\mathfrak{K}_{\mathscr{D}} = \emptyset$. Notice that $(\mathrm{Id} - \Delta_{\mathbb{S}^{d-1}})^{\beta}$ is always positive definite.

In order to gain further insight on Definition 4 and the motivations behind it, it is helpful to look at some key properties of spherical pseudo-differential operators:

**Proposition 1** *(Properties of pseudo-differential operators). Let $\mathscr{D}$ be a* spherical pseudo-differential operator *as in Definition 4. Then the following holds:*

1. *$\mathscr{D}$ is* self-adjoint*, i.e. $\mathscr{D}^* = \mathscr{D}$.*
2. *$\mathscr{D}$ is* isotropic*, i.e. any $Y \in Harm_n(\mathbb{S}^{d-1})$ is an* eigenfunction *of $\mathscr{D}$, with common eigenvalue $\lambda_n = \hat{D}_n$.*
3. *$\mathscr{D}$ has* finite-dimensional *nullspace, given by*

$$\begin{aligned}\mathcal{N}(\mathscr{D}) &= \left\{h \in \mathscr{S}(\mathbb{S}^{d-1}) : \langle h, Y_n^m\rangle_{\mathbb{S}^{d-1}} = 0,\, n \in \mathbb{N}\backslash\mathfrak{K}_{\mathscr{D}},\, m = 1,\dots,N_d(n)\right\}\\ &= span\left\{Y_n^m,\ n \in \mathfrak{K}_{\mathscr{D}},\ m = 1,\dots,N_d(n)\right\}.\end{aligned}$$

4. *$\mathscr{D}$ is an* endomorphism *on $\mathscr{S}(\mathbb{S}^{d-1})$, i.e. it maps infinitely differentiable functions onto infinitely differentiable functions.*

The proof of Proposition 1 is not detailed here due to space constraints but is relatively easy to obtain from the assumptions of Definition 4 [41, Proposition 4.1]:

- Properties 1 and 3 are direct consequences of the fact that the sequence $\{\hat{D}_n\}_{n\in\mathbb{N}}$ is respectively real and null for at most finitely many integers.
- The isotropy property 2 is implicitly assumed in Definition 4 since the spectrum of $\mathscr{D}$ was chosen in (9) to be constant for all $m = 1,\dots,N_d(n)$ in a given frequency level $n \in \mathbb{N}$. As shall be seen in

[19] Recall that the spherical harmonics were defined as eigenfunctions of the Laplace-Beltrami operator: $\forall Y \in \mathrm{Harm}_n(\mathbb{S}^{d-1})$, $\Delta_{\mathbb{S}^{d-1}} Y = -n(n+d-2)Y$.

Section 3.2, this construction guarantees that – when they exist – the spherical splines associated to a given pseudo-differential operator are sums of zonal functions, and hence fast to evaluate.

- Property 4 finally, results from the polynomial growth of the sequence $\{\hat{D}_n\}_{n\in\mathbb{N}}$. More specifically, it results from $\{\hat{D}_n\}_{n\in\mathbb{N}}$ being asymptotically bounded *from above* by a polynomial sequence $|\hat{D}_n| = \mathcal{O}(n^p)$, implied by $|\hat{D}_n| = \Theta(n^p)$.

Note that Definition 4 has an additional requirement that $\{\hat{D}_n\}_{n\in\mathbb{N}}$ be asymptotically bounded *from below* by a polynomial sequence – i.e. $|\hat{D}_n| = \Omega(n^p)$. This assumption comes into play when considering *primitives* w.r.t. a particular pseudo-differential operator $\mathscr{D}$, obtained via the *Moore-Penrose pseudo-inverse* $\mathscr{D}^\dagger$ of $\mathscr{D}$.

**Proposition 2** *(Moore-Penrose pseudo-inverse of $\mathscr{D}$). Let $\mathscr{D}$ be a pseudo-differential operator as in Definition 4. The* Moore-Penrose pseudo-inverse *$\mathscr{D}^\dagger$ of $\mathscr{D}$ is given by*

$$\mathscr{D}^\dagger : \begin{cases} \mathscr{S}(\mathbb{S}^{d-1}) \to \mathscr{S}(\mathbb{S}^{d-1}) \\ h \mapsto \mathscr{D}^\dagger h := \sum_{n\notin\mathfrak{K}_{\mathscr{D}}} \frac{1}{\hat{D}_n} \left[\sum_{m=1}^{N_d(n)} \hat{h}_n^m Y_n^m\right], \end{cases} \tag{12}$$

*where $\{\hat{h}_n^m,\, n\in\mathbb{N},\, m=1,\dots,N_d(n)\}$ are the* spherical Fourier coefficients *of $h$.*

**Proof.** The proof of Proposition 2 is available in Appendix A. □

Note that the primitive operator $\mathscr{D}^\dagger$ acts as an *integral operator* and smooths out high frequency content with a polynomially *decaying* sequence. In what follows, we will sometimes need to extend by duality the action of $\mathscr{D}$ (respectively $\mathscr{D}^\dagger$) to generalised functions $f \in \mathscr{S}'(\mathbb{S}^{d-1})$. Since $\mathscr{D}$ (respectively $\mathscr{D}^\dagger$) is self-adjoint, this can easily be achieved by understanding $\mathscr{D}f$ as the element of $\mathscr{S}'(\mathbb{S}^{d-1})$ with point-wise definition:

$$\langle \mathscr{D}f|\varphi\rangle := \langle f|\mathscr{D}\varphi\rangle, \quad \forall\varphi\in\mathscr{S}(\mathbb{S}^{d-1}). \tag{13}$$

Equation (13) is indeed well-defined since, from Item 4 of Proposition 1, $\mathscr{D}$ is an endomorphism on $\mathscr{S}(\mathbb{S}^{d-1})$.

### *3.2. Green functions and spline-admissibility*

The next important ingredient for the definition of spherical splines is the notion of *Green function* of a pseudo-differential operator $\mathscr{D}$. A Green function is a fundamental solution of $\mathscr{D}$, obtained by taking the primitive of some Dirac measure.

**Definition 5** *(Green function).* Let $\mathscr{D}$ be a *pseudo-differential operator* as in Definition 4. Consider moreover the *Moore-Penrose pseudo-inverse* $\mathscr{D}^\dagger$ of $\mathscr{D}$, extended into an endomorphism on $\mathscr{S}'(\mathbb{S}^{d-1})$ with (13). Then, a generalised function $\Psi_{\boldsymbol{s}}^{\mathscr{D}} \in \mathscr{S}'(\mathbb{S}^{d-1})$ is said to be a *Green function* for $\mathscr{D}$ if:

$$\Psi_{\boldsymbol{s}}^{\mathscr{D}} = \mathscr{D}^\dagger\delta_{\boldsymbol{s}}, \tag{14}$$

where $\delta_{\boldsymbol{s}} \in \mathcal{M}(\mathbb{S}^{d-1}) \subset \mathscr{S}'(\mathbb{S}^{d-1})$ is the *Dirac measure* for some direction $\boldsymbol{s}\in\mathbb{S}^{d-1}$.

The Green functions of an operator $\mathscr{D}$ can be expressed as traces of a certain zonal kernel, called the *zonal Green kernel*:

**Proposition 3** *(Zonal Green kernel). Let $\{\Psi_{\boldsymbol{s}}^{\mathscr{D}},\, \boldsymbol{s} \in \mathbb{S}^{d-1}\} \subset \mathscr{S}'(\mathbb{S}^{d-1})$ be Green functions for a pseudo-differential operator $\mathscr{D}$. We have then, for each $\boldsymbol{s} \in \mathbb{S}^{d-1}$:*

$$\langle \Psi_{\boldsymbol{s}}^{\mathscr{D}} | \varphi \rangle = \int_{\mathbb{S}^{d-1}} \psi_{\mathscr{D}}(\langle \boldsymbol{r}, \boldsymbol{s}\rangle)\varphi(\boldsymbol{r})\, d\boldsymbol{r}, \quad \forall \varphi \in \mathscr{S}(\mathbb{S}^{d-1}). \tag{15}$$

*The* zonal Green kernel $\psi_{\mathscr{D}}$ *is moreover such that $\{\psi_{\mathscr{D}}(\langle \cdot, \boldsymbol{s}\rangle),\, \boldsymbol{s} \in \mathbb{S}^{d-1}\} \subset \mathscr{S}'(\mathbb{S}^{d-1})$, and is defined as*

$$\psi_{\mathscr{D}}(\langle \boldsymbol{r}, \boldsymbol{s}\rangle) := \sum_{n \in \mathbb{N} \setminus \mathfrak{K}_{\mathscr{D}}} \frac{N_d(n)}{\mathfrak{a}_d \hat{D}_n} P_{n,d}(\langle \boldsymbol{r}, \boldsymbol{s}\rangle), \quad \boldsymbol{r} \in \mathbb{S}^{d-1}, \tag{16}$$

*where $\mathfrak{a}_d$ is the area of the unit sphere $\mathbb{S}^{d-1}$ and $P_{n,d} : [-1, 1] \to \mathbb{R}$ denotes the $d$-dimensional* ultraspherical polynomial *of degree $n \in \mathbb{N}$ [41, Definition 3.4].*

**Proof.** The proof of Proposition 3 is available in Appendix A. □

Observe that the traces of the zonal Green kernel (16) are generalised functions, which make sense when integrated against a Schwartz function but which may not admit a pointwise interpretation. When they do admit a pointwise interpretation, we say that the operator $\mathscr{D}$ is *spline-admissible*:

**Definition 6** *(Spline-admissible pseudo-differential operator).* Let $\mathscr{D}$ be a *pseudo-differential operator* with *zonal Green kernel* $\psi_{\mathscr{D}}$. We say that $\mathscr{D}$ is *spline admissible* if all traces $\{\psi_{\mathscr{D}}(\langle \cdot, \boldsymbol{s}\rangle),\, \boldsymbol{s} \in \mathbb{S}^{d-1}\} \subset \mathscr{S}'(\mathbb{S}^{d-1})$ of $\psi_{\mathscr{D}}$ are *ordinary functions*, i.e. they are *pointwise defined.*

The following result provides us with a sufficient condition for a pseudo-differential operator to be spline-admissible:

**Proposition 4** *(Sufficient condition for spline-admissibility). Let $\mathscr{D}$ be a pseudo-differential operator, with* spectral growth order *$p > d - 1$ and* zonal Green kernel *$\psi_{\mathscr{D}}$. Then we have*

$$\{\psi_{\mathscr{D}}(\langle \cdot, \boldsymbol{s}\rangle) : \mathbb{S}^{d-1} \to \mathbb{R},\ \boldsymbol{s} \in \mathbb{S}^{d-1}\} \subset \mathscr{C}(\mathbb{S}^{d-1}),$$

*and hence $\mathscr{D}$ is* spline-admissible.

**Proof.** The proof of Proposition 4 is available in Appendix A. □

We conclude this section by providing, for the specific case of $\mathbb{S}^2$, some examples (and non-examples) of spline-admissible pseudo-differential operators.

**Example 2** *(Common spline-admissible operators on $\mathbb{S}^2$).* Consider the specific case $d = 3$.

- **Laplace-Beltrami operator:** $\Delta_{\mathbb{S}^2}$ is *not* spline-admissible. Indeed, its zonal Green kernel is given by [45, Lemma 4.3],

$$\psi_{\Delta_{\mathbb{S}^2}}(\langle \boldsymbol{r}, \boldsymbol{s}\rangle) = \frac{1}{4\pi}\ln(1 - \langle \boldsymbol{r}, \boldsymbol{s}\rangle) + \frac{1}{4\pi} - \frac{1}{4\pi}\ln 2, \quad \forall \boldsymbol{r}, \boldsymbol{s} \in \mathbb{S}^{d-1},$$

  which is not defined for $\boldsymbol{r} = \boldsymbol{s}$. Note that we have $p = 2 = d - 1$ which shows that the bound on the spectral growth order in Proposition 4 is tight. $\Delta_{\mathbb{S}^2}^2$ on the other hand is, from Proposition 4, spline-admissible. Indeed, its spectral order $p$ is such that $p = 4 > 2 = d - 1$. Moreover, its zonal Green kernel admits a closed-form expression [45, Corollary 4.24].

- **Sobolev operators:** from Proposition 4, the Sobolev operators $(\mathrm{Id}-\Delta_{\mathbb{S}^2})^\beta$ are spline-admissible whenever $\beta > (d-1)/2 = 1$. There exists however no known closed-form expression for their zonal Green kernel, which is given by

$$\psi_\beta(\langle \boldsymbol{r}, \boldsymbol{s}\rangle) = \sum_{n=0}^{+\infty} \frac{2n+1}{4\pi\,(1+n(n+1))^\beta} P_n(\langle \boldsymbol{r}, \boldsymbol{s}\rangle), \qquad \forall \boldsymbol{r}, \boldsymbol{s} \in \mathbb{S}^{d-1}. \tag{17}$$

### 3.3. Spherical splines

We are now in a position to introduce spherical splines. Roughly speaking, spherical splines are primitives (w.r.t. a particular spline-admissible pseudo-differential operator) of Dirac streams with finite innovations:

**Definition 7** *($\mathscr{D}$-spline).* Let $\Xi_M = \{\boldsymbol{r}_1, \ldots, \boldsymbol{r}_M\} \subset \mathbb{S}^{d-1}$ be a *set of points* on the hypersphere and $\mathscr{D}$ a *spline-admissible pseudo-differential operator*. Then, a *$\mathscr{D}$-spline* is a generalised function $\mathfrak{s} \in \mathscr{S}'(\mathbb{S}^{d-1})$ such that

$$\mathscr{D}\mathfrak{s} = \sum_{i=1}^{M} \alpha_i\, \delta_{\boldsymbol{r}_i}, \tag{18}$$

where $\{\alpha_i,\, i = 1, \ldots, M\} \subset \mathbb{C}$ are called the *amplitudes* of the spline, while the directions $\boldsymbol{r}_i$ in the *knot set* $\Xi_M$ are called the *knots* of the spline. The pairs $(\alpha_i, \boldsymbol{r}_i)$ of amplitudes and knots are called the *innovations* of the spline, and their collection $\mathfrak{X}(\nu) = \{(\alpha_i, \boldsymbol{r}_i),\, i = 1, \ldots, M\}$ is called the *innovation set* of the spline.

Finally, we denote by

$$\mathfrak{S}_{\mathscr{D}}(\mathbb{S}^{d-1}, \Xi_M) := \left\{ \mathfrak{s} \in \mathscr{S}'(\mathbb{S}^{d-1}) : \mathscr{D}\mathfrak{s} = \sum_{i=1}^{M} \alpha_i \delta_{\boldsymbol{r}_i},\, \alpha_i \in \mathbb{C},\, \boldsymbol{r}_i \in \Xi_M \right\}$$

the linear subspace of $\mathscr{D}$-splines associated with the knot set $\Xi_M$.

**Remark 5** *(Non-trivial nullspace and constrained amplitudes).* Notice that (18) implicitly constrains the spline amplitudes $\alpha_i$ when $\mathscr{D}$ has a nontrivial nullspace. Indeed, for a Schwartz function $\varphi \in \mathscr{N}(\mathscr{D})$, we have from the definition of $\mathscr{D}$ for generalised functions $\langle \mathscr{D}\mathfrak{s} | \varphi \rangle = \langle \mathfrak{s} | \mathscr{D}\varphi \rangle = 0$. However, we also have from the right-hand side of (18): $\langle \mathscr{D}\mathfrak{s} | \varphi \rangle = \sum_{i=1}^{M} \alpha_i \varphi(\boldsymbol{r}_i)$. We have hence necessarily $\sum_{i=1}^{M} \alpha_i \varphi(\boldsymbol{r}_i) = 0$ for all $\varphi \in \mathscr{N}(\mathscr{D})$, which holds if and only if:

$$\sum_{i=1}^{M} \alpha_i Y_n^m(\boldsymbol{r}_i) = 0, \qquad \forall n \in \mathfrak{K}_{\mathscr{D}},\, m = 1, \ldots, N_d(n). \tag{19}$$

The following result characterises the splines associated to a spline-admissible operator in terms of its zonal Green kernel:

**Proposition 5** *(Characterisation of $\mathscr{D}$-splines). Let $\mathscr{D}$ be a* spline-admissible operator*, with* zonal Green kernel $\psi_{\mathscr{D}}$. *Let further $\mathfrak{s} \in \mathscr{S}'(\mathbb{S}^{d-1})$ be a $\mathscr{D}$-spline with knot set $\Xi_M = \{\boldsymbol{r}_1, \ldots, \boldsymbol{r}_M\} \subset \mathbb{S}^{d-1}$ and valid coefficients $\{\alpha_i\}_{i=1,\ldots,M} \subset \mathbb{C}$. Then, we have*

$$\mathfrak{s}(\boldsymbol{r}) = \sum_{i=1}^{M} \alpha_i \psi_{\mathscr{D}}(\langle \boldsymbol{r}, \boldsymbol{r}_i \rangle) \quad + \quad \sum_{n \in \mathfrak{K}_{\mathscr{D}}} \sum_{m=1}^{N_d(n)} \hat{\beta}_n^m Y_n^m(\boldsymbol{r}), \quad \forall \boldsymbol{r} \in \mathbb{S}^{d-1}, \tag{20}$$

*where* $\hat{\beta}_n^m := \langle \mathfrak{s} | Y_n^m \rangle$, $\forall n \in \mathfrak{K}_{\mathscr{D}}, m = 1, \ldots, N_d(n)$. *In particular, when* $\mathfrak{K}_{\mathscr{D}} = \emptyset$ *and* $p > d-1$, *we have*

$$\mathfrak{s}(\boldsymbol{r}) = \sum_{i=1}^{M} \alpha_i \psi_{\mathscr{D}}(\langle \boldsymbol{r}, \boldsymbol{r}_i \rangle), \qquad \forall \boldsymbol{r} \in \mathbb{S}^{d-1}, \tag{21}$$

*and*

$$\mathfrak{S}_{\mathscr{D}}(\mathbb{S}^{d-1}, \Xi_M) := \mathit{span}\{\psi_{\mathscr{D}}(\langle \cdot, \boldsymbol{r}_i \rangle), \boldsymbol{r}_i \in \Xi_M, i = 1, \ldots, M\} \subset \mathscr{L}^2(\mathbb{S}^{d-1}). \tag{22}$$

**Proof.** The proof of Proposition 5 is available in Appendix A. □

**Remark 6.** Observe from (21) that when $\mathfrak{K}_{\mathscr{D}} = \emptyset$ the $\mathscr{D}$-splines are linear combinations of zonal functions, and hence very easy to evaluate. This nice feature is due to the fact that we restricted ourselves to *isotropic* pseudo-differential operators in Definition 4.

## 4. Functional penalised basis pursuit on the hypersphere

In this section, we leverage the tools introduced in Sections 2 and 3 to formulate functional inverse problems on the hypersphere in a common *generalised sampling* framework. Our formulation allows us to see most spherical approximation problems as specific instances of a generic penalised optimisation problem. We investigate the latter in the specific case of generalised total variation penalties, yielding a functional penalised basis pursuit (FPBP) problem. We define the search space of a given FPBP problem, and characterise its predual in which the sampling linear functionals must be chosen. We moreover provide a *representer theorem* characterising the form of the solutions to FPBP problems in terms of sparse spherical splines.

### *4.1. Generalised sampling & functional inverse problems*

Most real-life spherical approximation problems take the form of *functional inverse problems*. In a typical inverse problem formulation, an unknown *spherical field*[20] $f \in \mathscr{B}'$ is probed by some *sensing device*, resulting in a data vector $\boldsymbol{y} = [y_1, \ldots, y_L] \in \mathbb{C}^L$ of $L$ *measurements*. To account for potential inaccuracies in the measurement process, the data vector $\boldsymbol{y}$ is often modelled as the outcome of a *random vector* $\boldsymbol{Y} = [Y_1, \ldots, Y_L] : \Omega \to \mathbb{C}^L$, fluctuating according to some application-dependent *noise distribution*.[21] When the measurement process is unbiased, entries of the expectation of $\boldsymbol{Y}$ can be thought of as the *ideal measurements* which would be obtained in a noise-free environment. In most cases, the ideal measurements are linked to the unknown spherical field by some *linear relationship*, called *generalised sampling* [24]:

$$\mathbb{E}[Y_i] = \langle f | \varphi_i \rangle, \qquad i = 1, \ldots, L, \tag{23}$$

where $\langle \cdot | \cdot \rangle : \mathscr{B}' \times \mathscr{B} \to \mathbb{C}$ denotes the *Schwartz duality product* for some duality pair $(\mathscr{B}, \mathscr{B}')$ and $\{\varphi_1, \ldots, \varphi_L\} \subset \mathscr{B}$ are *linear sampling functionals* modelling the action of the sensing device on the spherical field $f \in \mathscr{B}'$. Since most real-life acquisition systems react continuously to variations in their inputs, the dual space $\mathscr{B}'$ is generally equipped with the weak* topology, so that the linear functionals $\{\varphi_1, \ldots, \varphi_L\}$ modelling the instrument are all continuous (see Section 2.1.2). In such a formalism, the ideal measurements in (23) are often referred to as *generalised samples* [24] of $f$. This is because the Schwartz duality product is a generalised *evaluation* map $\langle f|\varphi_i\rangle = f(\varphi_i)$, allowing us to interpret the ideal measurements as samples of $f$ evaluated at "points" $\varphi_1, \ldots, \varphi_L \in \mathscr{B}$. For convenience, it is moreover customary to write the generalised sampling equations (23) in terms of a *sampling operator*[22] $\mathbf{\Phi} : \mathscr{B}' \to \mathbb{C}^L$ (see [44, Chapter 5]) defined as:

$$\mathbf{\Phi} : \begin{cases} \mathscr{B}' \to \mathbb{C}^L \\ f \mapsto [\langle f|\varphi_1\rangle, \cdots, \langle f|\varphi_L\rangle]. \end{cases} \tag{24}$$

[20] The generic appellation "spherical field" is used here to designate any element of $\mathscr{S}'(\mathbb{S}^{d-1})$, such as a function or a measure defined over the sphere.
[21] In the absence of noise, $\boldsymbol{Y}$ can simply be chosen as a deterministic random vector.

Reformulating (23) in terms of $\mathbf{\Phi}$ yields:

$$\mathbb{E}[\boldsymbol{Y}] = \begin{bmatrix} \mathbb{E}[Y_1] \\ \vdots \\ \mathbb{E}[Y_L] \end{bmatrix} = \begin{bmatrix} \langle f|\varphi_1\rangle \\ \vdots \\ \langle f|\varphi_L\rangle \end{bmatrix} = \mathbf{\Phi}(f). \tag{25}$$

The goal of a functional inverse problem is then to recover a spherical field $f \in \mathscr{B}'$ which best explains the observed generalised samples $\boldsymbol{y}$, given a particular noise and functional data model (25). Since the search space $\mathscr{B}'$ is infinite-dimensional and the data finite-dimensional, this task is fundamentally *ill-posed* and will in general elicit *infinitely many* candidate solutions. To discriminate among such solutions, it is customary to resort to *regularisation*, which can be seen as implementing *Occam's razor principle*[23] by favouring solutions with *simple* behaviours. This is typically achieved by means of *penalised convex optimisation problems* of the form:

$$\min_{f\in\mathscr{B}'} \left\{F\left(\boldsymbol{y}, \mathbf{\Phi}(f)\right) \;+\; \Lambda\left(|||f|||\right)\right\}, \tag{26}$$

where

- $F : \mathbb{C}^L \times \mathbb{C}^L \to \mathbb{R}_+ \cup \{+\infty\}$ is a *cost functional*, measuring the discrepancy between the *observed* and *predicted* generalised samples $\boldsymbol{y}$ and $\mathbf{\Phi}(f)$ respectively. Common choices of discrepancy measures are discussed in Remark 7. In what follows, we will assume that $F$ is such that for all $\boldsymbol{y} \in \mathbb{C}^L$,

$$F(\boldsymbol{y}, \cdot) : \begin{cases} \mathbb{C}^L \to \mathbb{R}_+ \cup \{+\infty\} \\ \boldsymbol{z} \mapsto F(\boldsymbol{y}, \boldsymbol{z}) \end{cases} \tag{27}$$

  is *proper*, *convex* and lower semi-continuous.
- $|||\cdot||| : \mathscr{B}' \to \mathbb{R}_+$ is the *dual norm* on $\mathscr{B}'$, called *regularisation norm*, which implements Occam's razor principle. Intuitively, elements $f \in \mathscr{B}'$ with small regularisation norm are *simple* and *well-behaved*, typically with a finite number of degrees of freedom (df).
- $\Lambda : \mathbb{R} \to \mathbb{R}_+$ is some *convex regularisation function*, *strictly increasing* on $\mathbb{R}_+$. In practice, $\Lambda$ often takes the form of a *monomial* $t \mapsto \lambda t^p$, where $p \geq 1$ and $\lambda > 0$. The parameter $\lambda$ is called *regularisation parameter* and controls the amount of regularisation by putting the regularisation norm and the cost functional on a similar scale.

**Remark 7** *(Choosing the cost functional).* In practice, the cost functional $F$ is often chosen in one of the following two ways:

[22] In finite dimensions, the sampling operator is generally called *forward*, *design* or *sensing* matrix.

[23] Occam's razor principle is a philosophical principle also known as the "law of briefness" or in Latin *lex parsimoniae*. It was supposedly formulated by William of Ockham in the 14th century, who wrote in Latin *"Entia non sunt multiplicanda praeter necessitatem"*. In English, this translates to *"More things should not be used than are necessary"*. In essence, this principle states that when two equally good explanations for a given phenomenon are available, one should always favour the simplest, i.e. the one that introduces the least explanatory variables.

- *Noiseless case:* In a noiseless setup, one has full trust in the generalised samples. It is therefore natural to require that any solution of (26) be *consistent* [44, Chapter 5] with the samples at hand, i.e. $\boldsymbol{y} = \boldsymbol{\Phi}(f)$, $\forall f \in \mathcal{V}$. This can be achieved by choosing the cost functional as $F(\boldsymbol{y}, \boldsymbol{\Phi}(f)) = \iota(\boldsymbol{y} - \boldsymbol{\Phi}(f))$, where $\iota : \mathbb{C}^L \to \{0, +\infty\}$ is the *indicator function* with value 0 if $\boldsymbol{z} = \boldsymbol{0}$ and $+\infty$ otherwise. Such cost functionals are for example used in *interpolation problems*.
- *Noisy case:* In a noisy setup, consistency is not desired anymore, as it almost always leads to *overfitting* the noisy data. In this case, one can use general $\ell_p$ cost functionals $F(\boldsymbol{y}, \boldsymbol{\Phi}(f)) = \|\boldsymbol{y} - \boldsymbol{\Phi}(f)\|_p^p$, where $p \in [1, +\infty]$ is typically chosen according to the tail behaviour[24] of the noise distribution [46]. Another approach consists in using the *negative log-likelihood* of the data $\boldsymbol{y}$ as a measure of discrepancy, i.e. $F(\boldsymbol{y}, \boldsymbol{\Phi}(f)) = -\ell(\boldsymbol{y}|\boldsymbol{\Phi}(f))$. This choice makes (26) resemble a *maximum a posteriori* problem with *improper prior*. In the case of centred Gaussian white noise, both discrepancy measures coincide, yielding the classical quadratic cost functional $F(\boldsymbol{y}, \boldsymbol{\Phi}(f)) = \|\boldsymbol{y} - \boldsymbol{\Phi}(f)\|_2^2$.

### 4.2. Generalised total variation regularisation

Notice that the regularisation norm $|||\cdot|||$ in (26) entirely determines the search space $\mathscr{B}'$ (see (4) page 8). Candidate regularisation norms for spherical approximation problems can hence be constructed as follows:

1. Identify *interesting* functional spaces $\mathscr{B}' \subset \mathscr{S}'(\mathbb{S}^{d-1})$, whose elements are *regular* enough;
2. Find a norm $|||\cdot|||$ on $\mathscr{B}'$ such that $\mathscr{B}'$ admits a predual $\mathscr{B}$ and characterise this predual.

For example, one could consider choosing $\mathscr{B}'$ as a *generalised Sobolev space* of the form:

$$\mathscr{H}_{\mathscr{D}}(\mathbb{S}^{d-1}) = \left\{ f \in \mathscr{S}'(\mathbb{S}^{d-1}) : \; \mathscr{D}f \in \mathscr{L}^2(\mathbb{S}^{d-1}) \right\}, \tag{28}$$

where $\mathscr{D} : \mathscr{S}'(\mathbb{S}^{d-1}) \to \mathscr{S}'(\mathbb{S}^{d-1})$ is some pseudo-differential operator as in Definition 4, which we will assume here to have trivial nullspace for simplicity. This is the space of generalised functions regular enough so that their generalised derivatives w.r.t. $\mathscr{D}$ are *square-integrable*. The associated regularisation norm can be shown to be the generalised Tikhonov (gTikhonov) norm $\|\mathscr{D}f\|_2$ [41, Section 2.1]. While extensively used in the literature, this notion of regularity may however be considered too restrictive, since the Sobolev space (28) is notably not large enough[25] to contain $\mathscr{D}$-splines in cases where $\mathscr{D}$ is *spline-admissible*. This is particularly cumbersome, since the latter are, by definition of spline-admissible operators, *ordinary functions* and hence relatively *well-behaved*. To include $\mathscr{D}$-splines, one must consider the larger space

$$\mathcal{M}_{\mathscr{D}}(\mathbb{S}^{d-1}) = \left\{ f \in \mathscr{S}'(\mathbb{S}^{d-1}) : \; \mathscr{D}f \in \mathcal{M}(\mathbb{S}^{d-1}) \right\}, \tag{29}$$

where $\mathcal{M}(\mathbb{S}^{d-1})$ denotes the space of *spherical regular Borel measures* introduced in (5). This is the space of generalised functions regular enough so that their generalised derivatives w.r.t. $\mathscr{D}$ are *Borel measures*. When $\mathscr{D}$ has trivial nullspace, $\mathcal{M}_{\mathscr{D}}(\mathbb{S}^{d-1})$ can be equipped with the generalised total variation (gTV) norm:

$$\|f\|_{\mathscr{D},TV} := \|\mathscr{D}f\|_{TV} = \sup_{\varphi \in \mathscr{S}(\mathbb{S}^{d-1}),\, \|\varphi\|_\infty = 1} |\langle \mathscr{D}f | \varphi \rangle|, \quad \forall f \in \mathcal{M}_{\mathscr{D}}(\mathbb{S}^{d-1}). \tag{30}$$

The gTV norm can be interpreted as measuring the *variations* of the generalised derivative $\mathscr{D}f$. Since $\mathscr{D}$ is bijective we can consider its inverse which can be shown to define an *isometric isomorphism* between the spaces $(\mathcal{M}(\mathbb{S}^{d-1}), \|\cdot\|_{TV})$ and $(\mathcal{M}_{\mathscr{D}}(\mathbb{S}^{d-1}), \|\cdot\|_{\mathscr{D},TV})$. Indeed, we can uniquely write any element $f$ in $\mathcal{M}_{\mathscr{D}}(\mathbb{S}^{d-1})$ as:

$$f = \mathscr{D}^{-1}\mu, \quad \mu \in \mathcal{M}(\mathbb{S}^{d-1}), \quad \text{with} \quad \|f\|_{\mathscr{D},TV} = \|\mathscr{D}f\|_{TV} = \|\mu\|_{TV}. \tag{31}$$

[24] The following rule of thumb is proposed in [46]: $p$ should be close to 1 for heavy-tailed distributions, close to 2 for Gaussian-like distributions, and close to $+\infty$ for compactly supported distributions.

[25] $\mathscr{D}$-splines are indeed defined as $\mathscr{D}$-primitives of Dirac streams, i.e. $\mathscr{D}\mathfrak{s} = \sum_{i=1}^M \alpha_i \delta_{\boldsymbol{r}_i}$. Their $\mathscr{D}$-derivatives are hence not in $\mathscr{L}^2(\mathbb{S}^{d-1})$.

This isometry implies that the metric space $(\mathcal{M}_{\mathscr{D}}(\mathbb{S}^{d-1}), \|\cdot\|_{\mathscr{D},TV})$ is actually a *Banach space*, and allows us to characterise its predual:

**Proposition 6** *(Predual of $\mathcal{M}_{\mathscr{D}}(\mathbb{S}^{d-1})$). The Banach space*

$$\mathscr{C}_{\mathscr{D}}(\mathbb{S}^{d-1}) = \left\{h \in \mathscr{S}'(\mathbb{S}^{d-1}) : \; h = \mathscr{D}\eta, \, \eta \in \mathscr{C}(\mathbb{S}^{d-1})\right\}, \tag{32}$$

*equipped with the norm $\|h\|_{\mathscr{D},\infty} = \|\mathscr{D}^{-1}h\|_\infty = \|\eta\|_\infty$, is the* predual *of the Banach space $(\mathcal{M}_{\mathscr{D}}(\mathbb{S}^{d-1}), \|\cdot\|_{\mathscr{D},TV})$,*

$$(\mathscr{C}_{\mathscr{D}}(\mathbb{S}^{d-1}), \|\cdot\|_{\mathscr{D},\infty})' \cong (\mathcal{M}_{\mathscr{D}}(\mathbb{S}^{d-1}), \|\cdot\|_{\mathscr{D},TV}).$$

**Proof.** The proof of Proposition 6 is available in Appendix B. □

We have hence established the duality pair

$$(\mathscr{C}_{\mathscr{D}}(\mathbb{S}^{d-1}), \|\cdot\|_{\mathscr{D},\infty})' \cong (\mathcal{M}_{\mathscr{D}}(\mathbb{S}^{d-1}), \|\cdot\|_{\mathscr{D},TV}),$$

showing that the gTV norm $\|\cdot\|_{\mathscr{D},TV}$ is actually a dual norm which can hence be used as regularisation norm in (26). For such a choice of regularisation norm, it is customary to set the regularisation function to $\Lambda(t) = \lambda t$ with $\lambda > 0$. This yields the following optimisation problem:

$$\min_{f \in \mathcal{M}_{\mathscr{D}}(\mathbb{S}^{d-1})} \left\{F\left(\boldsymbol{y}, \boldsymbol{\Phi}(f)\right) \; + \; \lambda \|\mathscr{D}f\|_{TV}\right\}, \tag{33}$$

where the sampling operator $\boldsymbol{\Phi} : \mathcal{M}_{\mathscr{D}}(\mathbb{S}^{d-1}) \to \mathbb{C}^L, f \mapsto [\langle f|\varphi_1\rangle, \ldots, \langle f|\varphi_L\rangle]$ is such that $\{\varphi_1, \ldots, \varphi_L\} \subset \mathscr{C}_{\mathscr{D}}(\mathbb{S}^{d-1})$. We call (33) a functional penalised basis pursuit (FPBP) problem. Because of the gTV regularisation norm, solutions to FPBP problems will tend to have *few* variations in their generalised derivatives. When $\mathscr{D}$ is *spline-admissible* and invertible, such functions are templated by the $\mathscr{D}$-splines, which, from Proposition 5, take the form

$$\mathfrak{s} = \sum_{i=1}^{M} \alpha_i \psi_{\mathscr{D}}(\langle \cdot, \boldsymbol{r}_i\rangle), \tag{34}$$

where $\{\boldsymbol{r}_1, \ldots, \boldsymbol{r}_M\} \subset \mathbb{S}^{d-1}$ and $\psi_{\mathscr{D}}$ is the zonal Green kernel of $\mathscr{D}$. For such functions, we have indeed

$$\|\mathscr{D}\mathfrak{s}\|_{TV} = \left\|\sum_{i=1}^{M} \alpha_i \mathscr{D}\psi_{\mathscr{D}}(\langle \cdot, \boldsymbol{r}_i\rangle)\right\|_{TV} = \left\|\sum_{i=1}^{M} \alpha_i \delta_{\boldsymbol{r}_i}\right\|_{TV} = \sum_{i=1}^{M} |\alpha_i| \underbrace{\|\delta_{\boldsymbol{r}_i}\|_{TV}}_{=1} = \|\boldsymbol{\alpha}\|_1.$$

Hence $\mathscr{D}$-splines with small $\ell_1$ norm in their coefficients will also have small gTV norm. It is then expected for solutions of (33) to take the form of $\mathscr{D}$-splines (34) with *few* innovations $M$. In Theorem 2 we will show that extreme points of the solution set of (33) indeed take such a form, with $M < L$.

*4.3. Representer theorem*

We now establish a representer theorem characterising the solution set of the FPBP problem (33). For simplicity, we state this theorem in the case where the pseudo-differential operator $\mathscr{D}$ used to define the gTV regularisation norm has a *trivial null space* ($\mathfrak{K}_{\mathscr{D}} = \emptyset$). While this may seem like a limiting assumption for practical purposes, we show that pseudo-differential operators with nontrivial nullspaces (such as the Laplace-Beltrami operator $\Delta_{\mathbb{S}^{d-1}}$) can be brought into the scope of the representer theorem if properly regularised on their nullspace. Our theorem shows that the solution set of an FPBP problem is *nonempty* and the *weak* closed convex-hull* of extreme points taking the form of $\mathscr{D}$-primitives of Dirac streams – i.e. $\mathscr{D}$-splines when $\mathscr{D}$ is spline-admissible – with less innovations than available measurements. This result can be seen as an extension to the spherical setup of [8, Theorem 4].

**Theorem 2** *(Representer theorem). Consider the following assumptions:*

*A1* $\mathscr{D} : \mathscr{S}'(\mathbb{S}^{d-1}) \to \mathscr{S}'(\mathbb{S}^{d-1})$ *is some* pseudo-differential operator *with* trivial nullspace *and* Green functions $\{\Psi^{\mathscr{D}}_{\boldsymbol{r}}, \boldsymbol{r} \subset \mathbb{S}^{d-1}\}$*;*
*A2* $(\mathcal{M}_{\mathscr{D}}(\mathbb{S}^{d-1}), \|\mathscr{D} \cdot\|_{TV})$ *is the space defined in* (29)*, with* topological dual $\mathscr{C}_{\mathscr{D}}(\mathbb{S}^{d-1})$ *given by* (32)*;*
*A3* *span*$\{\varphi_i, i = 1, \dots, L\} \subset \mathscr{C}_{\mathscr{D}}(\mathbb{S}^{d-1})$*, with the* $\varphi_i$ *being* linearly independent*;*
*A4* $\boldsymbol{\Phi} : \mathcal{M}_{\mathscr{D}}(\mathbb{S}^{d-1}) \to \mathbb{C}^L$ *is a* sampling *operator, defined as in* (24)*;*
*A5* $F : \mathbb{C}^L \times \mathbb{C}^L \to \mathbb{R}_+ \cup \{+\infty\}$ *is a cost functional as in* (27)*;*
*A6* $\lambda$ *is a* positive *regularisation constant.*

*Then, for any* $\boldsymbol{y} \in \mathbb{C}^L$*, the solution set of the FPBP problem*

$$\mathcal{V} = \operatorname*{Arg\,min}_{f \in \mathcal{M}_{\mathscr{D}}(\mathbb{S}^{d-1})} \{F(\boldsymbol{y}, \boldsymbol{\Phi}(f)) + \lambda \|\mathscr{D} f\|_{TV}\}, \tag{35}$$

*is* nonempty*, and the* weak* closed convex hull *of its* extreme points.[26] *The latter are moreover necessarily of the form:*

$$f^\star = \sum_{i=1}^M \alpha_i \, \mathscr{D}^{-1} \delta_{\boldsymbol{r}_i} = \sum_{i=1}^M \alpha_i \, \Psi^{\mathscr{D}}_{\boldsymbol{r}_i} = \sum_{i=1}^M \alpha_i \left[ \sum_{n \in \mathbb{N}} \sum_{m=1}^{N_d(n)} \frac{Y_n^m(\boldsymbol{r}_i)}{\hat{D}_n} Y_n^m \right], \tag{36}$$

*for some weights* $\{\alpha_1, \dots, \alpha_M\} \subset \mathbb{C}$*, directions* $\{\boldsymbol{r}_1, \dots, \boldsymbol{r}_M\} \subset \mathbb{S}^{d-1}$*, and where* $1 \le M \le L$*. In particular when* $\mathscr{D}$ *is* spline-admissible*, the extreme points of* $\mathcal{V}$ *are* ordinary functions *and take the form of* $\mathscr{D}$-splines

$$f^\star(\boldsymbol{r}) = \sum_{i=1}^M \alpha_i \, \psi_{\mathscr{D}}(\langle \boldsymbol{r}, \boldsymbol{r}_i \rangle), \qquad \forall \boldsymbol{r} \in \mathbb{S}^{d-1}, \tag{37}$$

*where* $\psi_{\mathscr{D}}$ *is the* zonal Green kernel *of* $\mathscr{D}$*.*

**Proof.** The proof of Theorem 2 is available in Appendix B. □

**Remark 8.** Theorem 2 allows us to write the solution set $\mathcal{V}$ of (35) as the weak* closed convex-hull of a (potentially infinite) set of extreme points $\delta\mathcal{V} \subset \mathcal{V}$:

---

[26] An extreme point $v \in \mathcal{V}$ is a point which does not lie on any open line segment joining two distinct points of $\mathcal{V}$ [41, Definition 2.4].

$$\mathcal{V} = \overline{\mathrm{Hull}\,(\delta\mathcal{V})}^{\mathrm{weak}^*}$$
$$= \overline{\left\{ \sum_{k=1}^{n} \alpha_{i_k} f^\star_{i_k} \,\middle|\, n \in \mathbb{N},\, \{i_1,\ldots,i_n\} \subset \mathbb{N},\, \sum_{k=1}^{n} \alpha_{i_k} = 1, \text{and } \ 0 \le \alpha_{i_k} \le 1, f^\star_{i_k} \in \delta\mathcal{V} \right\}}^{\mathrm{weak}^*}, \tag{38}$$

where the extreme points $f^\star \in \delta\mathcal{V}$ are, when $\mathscr{D}$ is spline-admissible, $\mathscr{D}$-splines of the form:

$$f^\star(\boldsymbol{r}) = \sum_{(\alpha_m,\boldsymbol{r}_m)\in\Xi(f^\star)} \alpha_m \psi_{\mathscr{D}}(\langle \boldsymbol{r}, \boldsymbol{r}_m\rangle), \quad \boldsymbol{r} \in \mathbb{S}^{d-1}, \quad \#\Xi(f^\star) = M(f^\star) \le L.$$

Note that the spline innovation sets $\{\Xi(f^\star),\, f^\star \in \delta\mathcal{V}\} \subset \mathscr{P}(\mathbb{C} \times \mathbb{S}^{d-1})$ are a priori *unknown* and *different* for each extreme point. However, they all have *bounded* cardinality $\#\Xi(f^\star) \le L$, where $L$ corresponds to the dimension of the data vector $\boldsymbol{y}$. Extreme points are hence $\mathscr{D}$-splines with *sparse innovations*: *they have at most as many degrees of freedom as available data.* This remarkable result is reminiscent of a similar property of basis pursuit problems in discrete setups (see [29, Theorem 19]). Unfortunately, it is only valid for extreme points and does not hold for arbitrary interior points in $\mathcal{V}^\circ = \mathcal{V}\backslash\delta\mathcal{V}$. Indeed, $\mathcal{V}^\circ$ consists in general of all *finite convex combinations of extreme points* in $\delta\mathcal{V}$, as well as *limits* of sequences of the latter under the weak* topology. Such limits may take the form of *infinite* summations and are hence not $\mathscr{D}$-splines anymore.

**Remark 9** *(gTV regularisation with non-injective operators).* Consider a pseudo-differential operator $\mathscr{D}$ with nontrivial nullspace $\mathcal{N}(\mathscr{D})$ and let $\Pi_{\mathcal{N}(\mathscr{D})} : \mathscr{S}'(\mathbb{S}^{d-1}) \to \mathcal{N}(\mathscr{D})$ be the *orthogonal projection operator*[27] onto $\mathcal{N}(\mathscr{D})$. Then, the regularised operator $\mathcal{L}_\gamma = \mathscr{D} + \gamma\Pi_{\mathcal{N}(\mathscr{D})}$ with $\gamma > 0$ is an injective pseudo-differential operator. The latter can hence be used to define a gTV regularisation norm. The FPBP problem associated to this choice of gTV regularisation norm is then:

$$\min_{f\in\mathcal{M}_{\mathcal{L}_\gamma}(\mathbb{S}^{d-1})} \left\{ F\left(\boldsymbol{y}, \boldsymbol{\Phi}(f)\right) \; + \; \lambda \|(\mathscr{D} + \gamma\Pi_{\mathcal{N}(\mathscr{D})}) f\|_{TV} \right\}, \tag{39}$$

where the sampling operator $\boldsymbol{\Phi} : \mathcal{M}_{\mathcal{L}_\gamma}(\mathbb{S}^{d-1}) \to \mathbb{C}^L, f \mapsto [\langle f|\varphi_1\rangle, \ldots, \langle f|\varphi_L\rangle]$ is such that $\{\varphi_1, \ldots, \varphi_L\} \subset \mathscr{C}_{\mathcal{L}_\gamma}(\mathbb{S}^{d-1})$. It is then possible to invoke Theorem 2 to show that the solution set of (39) is *nonempty* and the *weak* closed convex hull* of extreme points of the form

$$f^\star = \sum_{i=1}^{M} \alpha_i \left[\mathscr{D} + \gamma\Pi_{\mathcal{N}(\mathscr{D})}\right]^{-1} \delta_{\boldsymbol{r}_i} = \sum_{i=1}^{M} \alpha_i \left[\mathscr{D}^\dagger + \frac{1}{\gamma}\Pi_{\mathcal{N}(\mathscr{D})}\right] \delta_{\boldsymbol{r}_i}$$
$$= \sum_{i=1}^{M} \alpha_i \left[ \sum_{n\notin\mathfrak{K}_{\mathscr{D}}} \sum_{m=1}^{N_d(n)} \frac{Y_n^m(\boldsymbol{r}_i)}{\hat{D}_n} Y_n^m + \sum_{n\in\mathfrak{K}_{\mathscr{D}}} \sum_{m=1}^{N_d(n)} \frac{Y_n^m(\boldsymbol{r}_i)}{\gamma} Y_n^m \right], \tag{40}$$

with $1 \le M \le L$, $\{\alpha_1, \ldots, \alpha_M\} \in \mathbb{C}$ and $\{\boldsymbol{r}_1, \ldots, \boldsymbol{r}_M\} \subset \mathbb{S}^{d-1}$. Note that if $\gamma$ is close to zero, the solution is mainly contained in the nullspace $\mathcal{N}(\mathscr{D})$, while if $\gamma$ is very large, the solution is mainly contained in the orthogonal complement $\mathcal{N}(\mathscr{D})^\perp$ of the nullspace. Finally, when $\mathscr{D}$ is *spline-admissible*, it is moreover possible, using the addition theorem [21, Theorem 5.11], to rewrite (40) as an ordinary function, given by

$$f^\star(\boldsymbol{r}) = \sum_{i=1}^{M} \alpha_i \left[ \psi_{\mathscr{D}}(\langle \boldsymbol{r}, \boldsymbol{r}_i\rangle) + \frac{1}{\gamma} \sum_{n\in\mathfrak{K}_{\mathscr{D}}} \frac{N_d(n)}{\mathfrak{a}_d} P_{n,d}(\langle \boldsymbol{r}, \boldsymbol{r}_i\rangle) \right], \quad \boldsymbol{r} \in \mathbb{S}^{d-1}, \tag{41}$$

[27] The latter is easily shown to be given by $\Pi_{\mathcal{N}(\mathscr{D})} = \mathrm{Id} - \mathscr{D}\mathscr{D}^\dagger$, where $\mathscr{D}^\dagger$ is the *Moore-Penrose pseudo-inverse* of $\mathscr{D}$ in (12).

where $P_{n,d} : [-1,1] \to \mathbb{R}$ is the *ultraspherical polynomial* of degree $n$ [41, Definition 3.4], and $\mathfrak{a}_d > 0$ denotes the surface area of $\mathbb{S}^{d-1}$. Observe that although resembling it, (41) is *not* a $\mathscr{D}$-spline as the coefficients $\{\alpha_1, \ldots, \alpha_M\}$ have a priori no reason to verify $\sum_{i=1}^{M} \alpha_i Y_n^m(\boldsymbol{r}_i) = 0$, $\forall n \in \mathfrak{K}_{\mathscr{D}}$, $m = 1, \ldots, N_d(n)$, as requested for spherical splines associated to a pseudo-differential operator with nontrivial nullspace (see discussion in Remark 5).

### 4.4. Compatible sampling functionals

Note that assumption A3 of Theorem 2 requires the sampling functionals $\{\varphi_1, \ldots, \varphi_L\}$ to be *compatible* with the regularising operator $\mathscr{D}$, in the sense that they must be included in the predual $\mathscr{C}_{\mathscr{D}}(\mathbb{S}^{d-1})$ of the search space $\mathcal{M}_{\mathscr{D}}(\mathbb{S}^{d-1})$. In practice, this assumption may not always be easy to verify. In this section, we therefore provide sufficient conditions on the spectral growth order of $\mathscr{D}$ for Dirac measures and square-integrable functions – which are the two most common types of sampling functionals encountered in practice – to be compatible with a given pseudo-differential operator $\mathscr{D}$.

**Proposition 7** *(Compatibility of Dirac measures). Let $\mathscr{D}$ be a* pseudo-differential operator *as in Definition 4 with* spectral growth order *$p > d - 1$, and $(\mathcal{M}_{\mathscr{D}}(\mathbb{S}^{d-1}), \|\mathscr{D} \cdot\|_{TV})$ the space defined in (29) equipped with the gTV norm. Then, all Dirac measures are included in the predual $\mathscr{C}_{\mathscr{D}}(\mathbb{S}^{d-1})$ of $\mathcal{M}_{\mathscr{D}}(\mathbb{S}^{d-1})$, i.e. $\{\delta_{\boldsymbol{r}}, \ \boldsymbol{r} \in \mathbb{S}^{d-1}\} \subset \mathscr{C}_{\mathscr{D}}(\mathbb{S}^{d-1})$.*

**Proof.** The proof of Proposition 7 is available in Appendix B. □

**Proposition 8** *(Compatibility of square-integrable functions). Let $\mathscr{D}$ be a* pseudo-differential operator *as in Definition 4 with trivial nullspace and* spectral growth order *$p > (d-1)/2$, and $(\mathcal{M}_{\mathscr{D}}(\mathbb{S}^{d-1}), \|\mathscr{D} \cdot\|_{TV})$ the space defined in Eq. (28) equipped with the gTV norm. Then, all square-integrable functions are included in the predual $\mathscr{C}_{\mathscr{D}}(\mathbb{S}^{d-1})$ of $\mathcal{M}_{\mathscr{D}}(\mathbb{S}^{d-1})$, i.e. $\mathscr{L}^2(\mathbb{S}^{d-1}) \subset \mathscr{C}_{\mathscr{D}}(\mathbb{S}^{d-1})$.*

**Proof.** The proof of Proposition 8 is available in Appendix B. □

## 5. Practical implications

Theorem 2 has profound practical implications, which we now outline. First, we use (37) to derive a discretisation scheme allowing us, under mild assumptions on $\mathscr{D}$, to approximate with arbitrary precision functions in the solution set $\mathcal{V}$ of a given FPBP problem. The obtained discretisation scheme is moreover canonical to the FPBP problem at hand, as it transforms it into a traditional *discrete* penalised basis pursuit (PBP) problem [29]. Depending on the nature of the convex cost functional, we propose to solve this discrete PBP problem with two state-of-the-art first-order proximal algorithms [30], namely the *accelerated proximal gradient descent method* [32] or the *primal-dual splitting method* [31]. We conclude by discussing suitable choices of regularising pseudo-differential for practical purposes. We introduce *Wendland* and *Matérn* operators, whose Green functions have simple closed-form expressions and good localisation in space. The latter are moreover closely related to the popular *Sobolev operators* $\mathscr{D}_\beta = [\mathrm{Id} - \Delta_{\mathbb{S}^{d-1}}]^\beta$, $\beta > (d-1)/2$.

### 5.1. Canonical discretisation of FPBP problems

In this section, we use Theorem 2 to derive a *canonical search space discretisation scheme* for the FPBP problem (35). The idea of search space discretisation schemes is to restrict the search space $\mathscr{B}'$ of an optimisation problem to a well-chosen finite-dimensional subspace of the form $\mathrm{span}\{\psi_1, \ldots, \psi_N\} \subset \mathscr{B}'$. In the particular case of (35), Theorem 2 tells us that the solution set $\mathcal{V}$ is the *closed convex-hull* of $\mathscr{D}$-splines[28] with *sparse*[29] innovation sets (see Remark 8). This essentially means that any non limit point[30] of $\mathcal{V}$ is itself a $\mathscr{D}$-spline, as *finite* convex combination of $\mathscr{D}$-splines. Our goal is therefore to find a finite family of functions capable of approximating well enough any arbitrary $\mathscr{D}$-spline. To this end, it will help to characterise the functional space in which $\mathscr{D}$-splines naturally live, called the *native space* [21, Chapter 6]. When $\mathscr{D}$ is a *positive-definite operator* with spectral growth order $p > d-1$, the native space is the generalised Sobolev space $\mathscr{H}_{\mathscr{D}^{1/2}}(\mathbb{S}^{d-1})$ associated to the Hermitian square-root $\mathscr{D}^{1/2}$ of $\mathscr{D}$:

**Proposition 9** *(Native space for $\mathscr{D}$-splines). Let $\mathscr{D}$ be a* positive-definite pseudo-differential operator *with* spectral growth order $p > d-1$. *Then the Hermitian square-root $\mathscr{D}^{1/2}$ of $\mathscr{D}$ is a pseudo-differential operator with Fourier coefficients $\{\sqrt{\hat{D}_n},\ n \in \mathbb{N}\} \subset \mathbb{R}_+$. Moreover, the generalised Sobolev space $\mathscr{H}_{\mathscr{D}^{1/2}}(\mathbb{S}^{d-1})$ is a reproducing kernel Hilbert space (RKHS) which contains all $\mathscr{D}$-splines.*

**Proof.** The proof of Proposition 9 is available in Appendix C. □

The next proposition, adapted from [21, Theorem 6.36], shows the error incurred by approximating elements of $\mathscr{H}_{\mathscr{D}^{1/2}}(\mathbb{S}^{d-1})$ – and hence in particular arbitrary $\mathscr{D}$-splines of interest here – with $\mathscr{D}$-splines with fixed knot set $\Xi_N \subset \mathbb{S}^{d-1}$ of size $N$.

**Proposition 10** *(Approximation error analysis). Consider a* knot set $\Xi_N = \{\boldsymbol{r}_1, \dots, \boldsymbol{r}_N\} \subset \mathbb{S}^{d-1}$ *with* nodal width

$$\Theta_{\Xi_N} := \max_{\boldsymbol{r} \in \mathbb{S}^{d-1}} \min_{\boldsymbol{s} \in \Xi_N} \|\boldsymbol{r} - \boldsymbol{s}\|_{\mathbb{R}^d}. \tag{42}$$

*Let further $\mathscr{D}$ denote a* positive-definite, spline-admissible pseudo-differential operator *with* spectral growth order $p > \frac{d+1}{2}$ *and $\mathfrak{S}_{\mathscr{D}}(\mathbb{S}^{d-1}, \Xi_N)$ be the space of* spherical $\mathscr{D}$-splines *associated to the knot set $\Xi_N$. Then, for every function $h \in \mathscr{H}_{\mathscr{D}^{1/2}}(\mathbb{S}^{d-1})$ we have*

$$\frac{\|h - \mathfrak{s}_N^{\perp}\|_{\infty}}{\|h\|_{\mathscr{D}^{1/2}}} \leq 2^{3/2} L_{\mathscr{D}} \sqrt{\Theta_{\Xi_N}}, \tag{43}$$

*where $\|h\|_{\mathscr{D}^{1/2}} := \sqrt{\langle \mathscr{D}^{1/2} h, \mathscr{D}^{1/2} h \rangle_{\mathbb{S}^{d-1}}}$, $L_{\mathscr{D}} > 0$ is a known*[31] *positive constant depending only on $\mathscr{D}$ and $\mathfrak{s}_N^{\perp} \in \mathfrak{S}_{\mathscr{D}}(\mathbb{S}^{d-1}, \Xi_N)$ is a $\mathscr{D}$-spline verifying $\mathfrak{s}_N^{\perp} := \arg\min_{\mathfrak{s} \in \mathfrak{S}_{\mathscr{D}}(\mathbb{S}^{d-1}, \Xi_N)} \|h - \mathfrak{s}\|_{\mathscr{D}^{1/2}}$, i.e. $\mathfrak{s}_N^{\perp}$ is the* orthogonal projection *of $h$ onto $\mathfrak{S}_{\mathscr{D}}(\mathbb{S}^{d-1}, \Xi_N)$.*

**Proof.** The proof of Proposition 10 is available in Appendix C. □

Notice that the approximation error in Proposition 10 is bounded by the nodal width (42) of the knot set, which can be interpreted geometrically as the largest distance from an arbitrary point on the sphere to the knot set $\Xi_N$ (see Fig. 2a). $\mathscr{D}$-splines with knot sets minimising this quantity for a fixed number of knots $N$ will hence yield the smallest approximation error. From the geometric interpretation of the nodal width, it is easy to see that knot sets with minimal nodal width distribute their knots *uniformly* over the hypersphere. Unfortunately, distributing points uniformly over $\mathbb{S}^{d-1}$ is a notoriously hard problem for $d > 2$

[28] Of course, under the assumption that the operator $\mathscr{D}$ used to define the gTV regularisation norm is spline-admissible.
[29] I.e. with cardinality bounded by the number of available measurements.
[30] As explained in Remark 8, the limit points are not necessarily splines. We are however not interested in approximating such limit points, since they may have infinitely many df.
[31] As shown in the proof, $L_{\mathscr{D}}$ is the *uniform Lipschitz constant* of the zonal Green kernel $\psi_{\mathscr{D}}$.

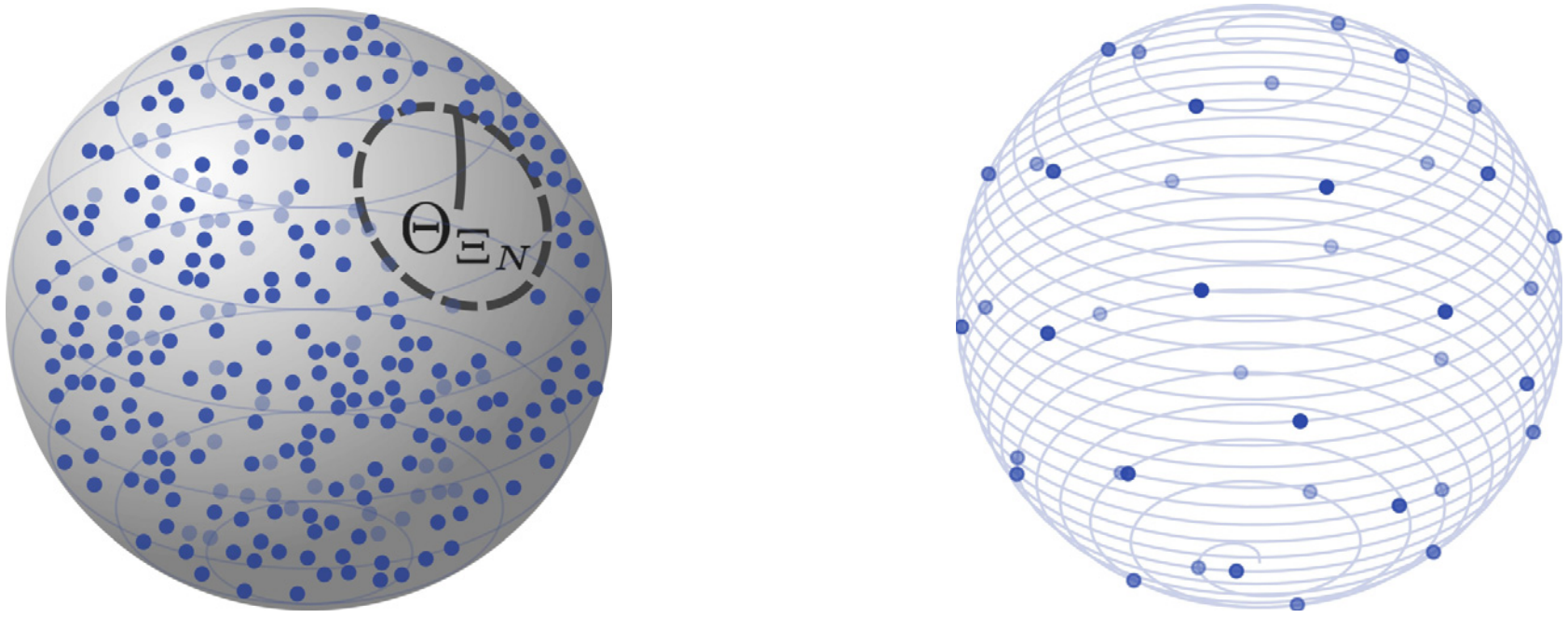

(a) The nodal width of a knot set represents the largest distance from an arbitrary point on the sphere to the knot set.

(b) Fibonacci lattice with $N = 41$ points.

**Fig. 2.** Visual representation of the nodal width (a) and the quasi-uniform Fibonacci lattice (b).

[9], making uniform knot sets impractical. For $d = 3$ however, it is possible to obtain *quasi-uniform* knot sets with *quasi-optimal* nodal widths [9]. An example of quasi-uniform knot set is the *Fibonacci lattice* [13,9] described in the subsequent example. In [9], the authors provide a comprehensive list of quasi-uniform knot sets easy to generate in practice. For each knot set, the asymptotic behaviour of the nodal width is assessed, either numerically or theoretically.

**Example 3** *(Fibonacci lattice).* In nature, many plant leaves are arranged according to *phyllotactic spiral patterns*, which are well modelled by the *Fibonacci lattice*. Points in the Fibonacci lattice are arranged uniformly along a spiral pattern on the sphere linking the two poles (see Fig. 2b). The lattice can very easily be generated from the following formula:

$$\begin{cases} \boldsymbol{r}_n = [\cos(\varphi_n)\sin(\theta_n), \sin(\varphi_n)\sin(\theta_n), \cos(\theta_n)]\,, \\ \text{where} \quad \varphi_n = 2\pi n \left(1 - \frac{2}{1+\sqrt{5}}\right) \quad \& \quad \theta_n = \arccos\left(1 - \frac{2n}{N}\right), \end{cases} \tag{44}$$

with $n = 1, \dots, N$. It can be shown [9] that if the knot set $\Xi_N$ is constructed according to the Fibonacci lattice (44), then the nodal width is *quasi-optimal* and approximately given by $\Theta_{\Xi_N} \simeq 2.728/\sqrt{N}$.

Observe that the nodal width of the Fibonacci lattice tends to zero as the number of knots $N$ grows to infinity. This is a general behaviour of quasi-uniform knot sets [9]. Consequently, the uniform approximation error (43) in Proposition 10 tends to zero as the number of knots tends to infinity. In other words, any element of $\mathscr{H}_{\mathscr{D}^{1/2}}(\mathbb{S}^{d-1})$ can be approximated *arbitrarily well* by $\mathscr{D}$-splines with *quasi-uniform* knot sets – called *quasi-uniform spherical splines*, provided a sufficient number of knots. In light of this discussion, we therefore propose to discretise FPBP problems by restricting their search spaces to subspaces spanned by quasi-uniform $\mathscr{D}$-splines. The following theorem shows that the solutions to FPBP problems restricted this way can be obtained by solving a *discrete* penalised basis pursuit (PBP) problem.

**Theorem 3** *(Canonical discretisation of FPBP problems). Consider the notations and assumptions A1 to A6 of Theorem 2. Consider additionally the following:*

*B7* $\mathscr{D}$ *is* spline-admissible *and* positive-definite*;*
*B8* $\{\psi_1, \dots, \psi_N\} \subset \mathcal{M}_{\mathscr{D}}(\mathbb{S}^{d-1})$ *are* zonal functions *of the form*

$$\psi_n := \psi_{\mathscr{D}}(\langle \boldsymbol{r}, \boldsymbol{r}_n \rangle), \quad \forall n = 1, \dots, N, \tag{45}$$

*where* $\psi_{\mathscr{D}}$ *is the* zonal Green kernel *of* $\mathscr{D}$ *and* $\Xi_N = \{\boldsymbol{r}_1, \ldots, \boldsymbol{r}_N\} \subset \mathbb{S}^{d-1}$ *for some* $N \in \mathbb{N}$*;*

*B9* $\mathfrak{S}_{\mathscr{D}}(\mathbb{S}^{d-1}, \Xi_N) = span\{\psi_1, \ldots, \psi_N\} \subset \mathcal{M}_{\mathscr{D}}(\mathbb{S}^{d-1})$ *is the space of* $\mathscr{D}$*-splines with knot set* $\Xi_N$*;*

*B10* $\Psi : \mathbb{C}^N \to \mathfrak{S}_{\mathscr{D}}(\mathbb{S}^{d-1}, \Xi_N)$ *is a* synthesis operator*, defined as* $\Psi(\boldsymbol{x}) := \sum_{n=1}^N x_n \psi_n$, $\forall \boldsymbol{x} \in \mathbb{C}^N$.

*Then, for each* $\boldsymbol{y} \in \mathbb{C}^L$*, the* restricted FPBP problem

$$\mathcal{V} = \operatorname*{Arg\,min}_{f \in \mathfrak{S}_{\mathscr{D}}(\mathbb{S}^{d-1}, \Xi_N)} \{F(\boldsymbol{y}, \boldsymbol{\Phi}(f)) \quad + \quad \lambda \|\mathscr{D} f\|_{TV}\} \tag{46}$$

*and the following* discrete PBP problem

$$\mathfrak{U} = \operatorname*{Arg\,min}_{\boldsymbol{x} \in \mathbb{C}^N} \{F(\boldsymbol{y}, \boldsymbol{G}\boldsymbol{x}) \quad + \quad \lambda \|\boldsymbol{x}\|_1\} \tag{47}$$

*are* equivalent*, in the sense that their solution sets are in* bijection *with one another:*

$$\mathcal{V} = \Psi(\mathfrak{U}) \qquad \textit{and} \qquad \mathfrak{U} = \boldsymbol{\Psi}^\dagger(\mathcal{V}), \tag{48}$$

*where* $\boldsymbol{\Psi}^\dagger : \mathfrak{S}_{\mathscr{D}}(\mathbb{S}^{d-1}, \Xi_N) \to \mathbb{C}^N$ *is the* Moore-Penrose pseudo-inverse *of* $\Psi$*. Moreover, the matrix* $\boldsymbol{G} := \boldsymbol{\Phi}\Psi \in \mathbb{C}^{L \times N}$ *in* (47) *is given by* $G_{ln} := \langle \psi_{\mathscr{D}}(\langle \cdot, \boldsymbol{r}_n \rangle) | \varphi_l \rangle$, $l = 1, \ldots, L$, $n = 1, \ldots, N$*, which simplifies to* $G_{ln} = (\psi_{\mathscr{D}} * \varphi_l)(\boldsymbol{r}_n)$ *when the sampling functionals* $\{\varphi_l, \ l = 1, \ldots, L\}$ *are in* $\mathscr{L}^2(\mathbb{S}^{d-1})$.

**Proof.** The proof of Theorem 3 is available in Appendix C. □

**Remark 10** *(Canonical discretisation scheme).* Notice that the discretisation scheme chosen in Theorem 3 is *canonical* w.r.t. the gTV norm induced by the pseudo-differential operator $\mathscr{D}$. Indeed, it conveniently transforms the gTV norm $\|\mathscr{D}\cdot\|_{TV}$ into a discrete $\ell_1$ norm in (47). As detailed in the proof, this is because the basis functions $\{\psi_{\mathscr{D}}(\langle \cdot, \boldsymbol{r}_n \rangle), \ n = 1, \ldots, N\}$ used in the discretisation are Green functions of the operator $\mathscr{D}$. Had the basis functions been chosen differently, such simplifications would not have been possible, hence making the discrete optimisation problem (47) considerably more difficult to solve in practice.

**Remark 11** *(Choice of N).* The bound in (43) can be used in practice to set $N$. Indeed, one can choose $N$ such that the relative approximation error falls below an acceptable accuracy threshold for any $h \in \mathscr{H}_{\mathscr{D}^{1/2}}(\mathbb{S}^{d-1})$, hence allowing us to approximate the solutions of the FPBP problem with controlled error.

**Remark 12** *(Practical implementation).* Theorem 3 provides us with a simple two-step procedure for computing a practical solution to the restricted FPBP problem (46):

1. Minimise (47) using one of the algorithms described in Section 5.2 and obtain a solution $\boldsymbol{u} \in \mathfrak{U}$.
2. Using the synthesis operator $\Psi : \mathbb{C}^N \to \mathfrak{S}(\mathbb{S}^{d-1}, \Xi_N)$ and the fact that $\mathcal{V} = \Psi(\mathfrak{U})$, map $\boldsymbol{u}$ into a solution $f = \Psi(\boldsymbol{u}) \in \mathcal{V}$ of the restricted FPBP problem (46):

$$f(\boldsymbol{r}) = (\Psi \boldsymbol{u})(\boldsymbol{r}) = \sum_{n=1}^N u_n \psi_{\mathscr{D}}(\langle \boldsymbol{r}, \boldsymbol{r}_n \rangle), \quad \forall \boldsymbol{r} \in \mathbb{S}^{d-1}.$$

The latter step can in this case be interpreted as an *interpolation* on $\mathbb{S}^{d-1}$ of samples $\{u_n, \ n = 1, \ldots, N\} \subset \mathbb{C}$ with sampling locations $\{\boldsymbol{r}_n, \ n = 1, \ldots, N\} \subset \mathbb{S}^{d-1}$ and interpolation kernel $\psi_{\mathscr{D}}$. Since the interpolating functions $\psi_{\mathscr{D}}(\langle \boldsymbol{r}, \boldsymbol{r}_n \rangle)$ are *zonal*, such an interpolation can be carried out *very efficiently* in practice (and even more so when $\psi_{\mathscr{D}}$ has compact support, as investigated in Section 5.3).

**Remark 13** *(Form of the solutions).* Applying [29, Theorem 6] to the discrete PBP problem (47) shows that $\mathfrak{U}$ is convex and compact with $L$-sparse extreme points. From the bijection (48), it implies in turn that $\mathcal{V} = \Psi(\mathfrak{U})$ is the closed convex-hull of extreme points taking the form of *sparse* $\mathscr{D}$-admissible spherical splines with at most $L$ non-zero amplitudes:

$$\forall f \in \delta\mathcal{V}, \quad f = \Psi\boldsymbol{u} \quad \text{with} \quad \|\boldsymbol{u}\|_0 \leq L,$$

where $\|\cdot\|_0$ denotes the "$\ell_0$ norm", counting the number of non-zero elements in a vector. Solutions of the restricted FPBP problem (46) behave hence very similarly as the ones of the unrestricted FPBP problem (35) considered in Theorem 2.

## 5.2. Optimisation algorithms

In this section, we propose to solve the discrete PBP problem (47) by means of provably convergent *fully-split proximal iterative methods* [30], which only involve simple matrix-vector multiplications and *proximal steps*. We treat the most general case where the cost function $F$ is proximable but not necessarily differentiable with the *primal-dual splitting method (PDS)* introduced by Condat in his seminal work [31]. In the simpler (yet prevailing in practice) case where $F$ is also *differentiable* and with *$\beta$-Lipschitz continuous* derivative, we leverage an optimal first-order method called *accelerated proximal gradient descent (APGD)* [32,30], with faster convergence rate than the PDS method. For the sake of simplicity and without loss of generality, we consider the *real case* only, where $\boldsymbol{x} \in \mathbb{R}^N$ and $\boldsymbol{y} \in \mathbb{R}^L$ – i.e. the coefficients and data vector are assumed real. The complex case, less common in practice, can be handled similarly by identifying $\mathbb{C}^N$ and $\mathbb{C}^L$ with $\mathbb{R}^{2N}$ and $\mathbb{R}^{2L}$ respectively and reformulating the optimisation problem (47) in terms of real operations only, using the techniques described in [47, Section 7.8].

### 5.2.1. Primal-dual splitting method for proximable cost functionals

Consider the real PBP problem:

$$\text{Find} \quad \hat{\boldsymbol{x}} \in \operatorname*{Arg\,min}_{\boldsymbol{x}\in\mathbb{R}^N} \left\{F(\boldsymbol{y}, \boldsymbol{G}\boldsymbol{x}) \quad + \quad \lambda\|\boldsymbol{x}\|_1\right\}, \tag{49}$$

where $F(\cdot;\boldsymbol{y}) : \mathbb{R}^L \to \mathbb{R}\cup\{\infty\}$ is a *proper, lower semi-continuous and convex* cost functional whose *proximity operator*

$$\mathbf{prox}_f(\boldsymbol{z}) := \operatorname*{arg\,min}_{\boldsymbol{x}\in\mathbb{C}^N} \left\{ f(\boldsymbol{x}) \quad + \quad \frac{1}{2}\|\boldsymbol{x}-\boldsymbol{z}\|_2^2 \right\}, \qquad \forall \boldsymbol{z} \in \mathbb{C}^N, \tag{50}$$

admits a closed-form formula or can be *efficiently* computed with *high accuracy.* Closed-form expressions for the proximity operators of commonly-used data-fidelity terms are available in [41, Chapter 7, Section 5] and are summarised in Table 2. To solve (49), we apply the *Chambolle-Pock primal-dual splitting algorithm* of [48, Algorithm 1] to our particular case and obtain Algorithm 1. As explained in [48,31], Algorithm 1 solves jointly the *primal* problem (49) and its associated *dual* using an equivalent *saddle-point formulation* of (49)

$$(\hat{\boldsymbol{x}}, \hat{\boldsymbol{z}}) \in \operatorname{Arg} \min_{\mathbb{R}^N} \sup_{\boldsymbol{z}\in\Delta(F^*_{\boldsymbol{y}})} \left\{\langle \boldsymbol{G}\boldsymbol{x}, \boldsymbol{z}\rangle_{\mathbb{R}^L} \;-\; F^*_{\boldsymbol{y}}(\boldsymbol{z}) \;+\; \lambda\|\boldsymbol{x}\|_1\right\}, \tag{51}$$

where $F^*_{\boldsymbol{y}} : \mathbb{R}^L \to \mathbb{R}\cup\{\infty\}$ is the *convex conjugate* of $F$, with *domain* $\Delta(F^*_{\boldsymbol{y}}) \subset \mathbb{R}^L$ and defined as [49]

$$F^*_{\boldsymbol{y}}(\boldsymbol{z}) := \sup_{\boldsymbol{x}\in\mathbb{R}^L} \langle \boldsymbol{z}, \boldsymbol{u}\rangle_{\mathbb{R}^L} - F(\boldsymbol{y}, \boldsymbol{u}), \quad \forall \boldsymbol{u} \in \mathbb{R}^L. \tag{52}$$

**Algorithm 1** The primal-dual splitting algorithm for solving the generic PBP problem (49).

1: **procedure** PDS_PBP($\boldsymbol{y}, \boldsymbol{G}, \lambda, \varepsilon, \tau, \sigma, \boldsymbol{x}_0, \boldsymbol{z}_0$)
2: $\quad n := 0,\ \delta := 0$
3: $\quad$ **while** $\delta = 0$ **do**
4: $\quad\quad n \leftarrow n + 1$
5: $\quad\quad \boldsymbol{x}_n := \mathbf{soft}_{\lambda\tau}\left(\boldsymbol{x}_{n-1} - \tau \boldsymbol{G}^T \boldsymbol{z}_{n-1}\right)$
6: $\quad\quad \boldsymbol{z}_n := \mathbf{prox}_{\sigma F^*_{\boldsymbol{y}}}\left(\boldsymbol{z}_{n-1} + \sigma \boldsymbol{G}\left[2\boldsymbol{x}_n - \boldsymbol{x}_{n-1}\right]\right)$
7: $\quad\quad$ **if** $\|\boldsymbol{x}_n - \boldsymbol{x}_{n-1}\|_2 \leq \varepsilon \|\boldsymbol{x}_{n-1}\|_2$ **then** $\delta \leftarrow 1$
8: $\quad$ **return** $\boldsymbol{x}_n$ ▷ Approximate solution to (49).

**Algorithm 2** The accelerated proximal gradient descent algorithm for solving smooth PBP problems (49).

1: **procedure** APGD_PBP($\boldsymbol{y}, \boldsymbol{G}, \lambda, \varepsilon, \tau, \mathfrak{d}, \boldsymbol{x}_0$)
2: $\quad n := 0,\ \delta := 0,\ \boldsymbol{z}_0 := \boldsymbol{x}_0$
3: $\quad$ **while** $\delta = 0$ **do**
4: $\quad\quad n \leftarrow n + 1$
5: $\quad\quad \boldsymbol{z}_n := \mathbf{soft}_{\lambda\tau}\left(\boldsymbol{x}_{n-1} - \tau \boldsymbol{\nabla} E_{\boldsymbol{y}}(\boldsymbol{x}_{n-1})\right)$
6: $\quad\quad \boldsymbol{x}_n := \boldsymbol{z}_n + \frac{n-1}{n+\mathfrak{d}}(\boldsymbol{z}_n - \boldsymbol{z}_{n-1})$
7: $\quad\quad$ **if** $\|\boldsymbol{x}_n - \boldsymbol{x}_{n-1}\|_2 \leq \varepsilon \|\boldsymbol{x}_{n-1}\|_2$ **then** $\delta \leftarrow 1$
8: $\quad$ **return** $\boldsymbol{x}_n$ ▷ Approximate solution to (49).

**Table 2**
Common data-fidelity functionals and their associated proximity operators.

| Name | $F(\boldsymbol{y}, \boldsymbol{z}),\ \boldsymbol{y}, \boldsymbol{z} \in \mathbb{R}^L$ | $\mathbf{prox}_{\tau F(\boldsymbol{y}, \cdot)}(\boldsymbol{z}),\ \tau > 0,\ \boldsymbol{z} \in \mathbb{R}^L$ | Useful for |
|---|---|---|---|
| **Exact match** | $\iota(\boldsymbol{z} - \boldsymbol{y})$ | $\boldsymbol{y}$ | Noiseless data, interpolation. |
| **$\ell_1$-norm** | $\|\boldsymbol{z} - \boldsymbol{y}\|_1$ | $\mathbf{soft}_\tau(\boldsymbol{z} - \boldsymbol{y}) + \boldsymbol{y}$ | Strong outliers and heavy-tailed noise distributions. |
| **$\ell_2$-ball** | $\iota_{\mathcal{B}_{2,\epsilon}}(\boldsymbol{z} - \boldsymbol{y}),\ \epsilon > 0$ | $\epsilon \frac{\boldsymbol{z} - \boldsymbol{y}}{\|\boldsymbol{z} - \boldsymbol{y}\|} + \boldsymbol{y}$ | Gaussian noise with known noise level. |
| **Generalised KL-divergence** | $\sum_{i=1}^L y_i \log\left(\frac{y_i}{z_i}\right) - y_i + z_i\ \forall \boldsymbol{z}, \boldsymbol{y} \in \mathbb{R}^L_+$ | $\frac{1}{2}(\boldsymbol{z} - \tau + \sqrt{(\boldsymbol{z} - \tau)^2 + 4\boldsymbol{y}\tau})$ | Count data with Poisson noise. |

The proximal operator of the convex conjugate (52) is given, for every $\sigma > 0$, by *Moreau's identity* [30]

$$\mathbf{prox}_{\sigma F^*_{\boldsymbol{y}}}(\boldsymbol{z}) = \boldsymbol{z} - \sigma \mathbf{prox}_{F(\boldsymbol{y}, \cdot)/\sigma}(\boldsymbol{z}/\sigma),\ \forall \boldsymbol{z} \in \mathbb{R}^L. \tag{53}$$

The saddle-point problem (51) allows the complicated composite term $F(\boldsymbol{y}, \boldsymbol{G}\boldsymbol{x})$ in (49) to be split as a sum of two simpler terms $F^*_{\boldsymbol{y}}(\boldsymbol{z}) + \langle \boldsymbol{G}\boldsymbol{x}, \boldsymbol{z} \rangle$, which are easier to optimise. Estimates $(\hat{\boldsymbol{x}}, \hat{\boldsymbol{z}})$ of the primal and dual variables are then obtained by iterating rows 5 and 6 of Algorithm 1 until convergence, starting from an initial guess $(\boldsymbol{x}_0, \boldsymbol{z}_0) \in \mathbb{R}^N \times \mathbb{R}^L$. Notice that the update equations 5 and 6 are not too computationally intensive since they involve only proximal operators – namely the *soft-thresholding operator* [30, Chapter 6] and (53) – and matrix-vector multiplications between $\boldsymbol{G}$, its transpose $\boldsymbol{G}^T$ and the primal/dual variables. Regarding the stopping criterion, we chose to stop Algorithm 1 when the relative improvement in the primal variable $\boldsymbol{x}_n$ falls under a certain pre-determined threshold $\varepsilon > 0$. This stopping criterion, which monitors improvement of the primal variable, is motivated by the fact that we are in this context only interested in solving the primal problem (49). As reported in [31], the convergence of the iterates $(\boldsymbol{x}_n, \boldsymbol{z}_n)$ towards a solution $(\hat{\boldsymbol{x}}, \hat{\boldsymbol{z}})$ of (51) as $n$ tends to infinity depends on the choice of the hyper-parameters $\tau > 0$ and $\sigma > 0$, which control respectively the amount of improvement in the primal and dual variable at each iteration. To ensure convergence, one must have [31, Theorem 3.3]: $\sigma\tau\|\boldsymbol{G}\|_2^2 \leq 1$. In practice, the speed of convergence is improved by choosing $\sigma$ and $\tau$ as large as possible and relatively well-balanced. Consequently, we chose in our implementation of Algorithm 1 $\sigma\tau = 1/\|\boldsymbol{G}\|_2^2$ and $\tau = \sigma$, yielding: $\sigma = \tau = 1/\|\boldsymbol{G}\|_2$.

### 5.2.2. *Accelerated proximal gradient descent method for smooth cost functionals*

In this section, we consider the special case of *smooth* PBP problems (49), where the cost functional $E_{\boldsymbol{y}}(\boldsymbol{x}) := F(\boldsymbol{y}, \boldsymbol{G}\boldsymbol{x}),\ \forall \boldsymbol{x} \in \mathbb{R}^N$, is assumed *differentiable* and with *$\beta$-Lipschitz gradient* $\boldsymbol{\nabla} E_{\boldsymbol{y}} : \mathbb{R}^N \to \mathbb{R}^N$:

$$\|\boldsymbol{\nabla} E_{\boldsymbol{y}}(\boldsymbol{x}) - \boldsymbol{\nabla} E_{\boldsymbol{y}}(\boldsymbol{z})\|_2 \leq \beta \|\boldsymbol{x} - \boldsymbol{z}\|_2, \qquad \forall \boldsymbol{x}, \boldsymbol{z} \in \mathbb{R}^N,$$

for some constant $\beta \geq 0$. While trivially in the scope of Algorithm 1, such problems are however more efficiently optimised by means of Algorithm 2, a specific instance of the generic *accelerated proximal gradient* *descent (APGD)* method [50,51]. Indeed, it has been shown in [50] that, with $0 < \tau \leq 1/\beta$ and $\mathfrak{d} > 2$, APGD achieves the following optimal convergence rates:

$$\lim_{n\to\infty} n^2 \left|\mathcal{J}(\boldsymbol{x}_n) - \mathcal{J}(\boldsymbol{x}^\star)\right| = 0 \qquad \& \qquad \lim_{n\to\infty} n^2\|\boldsymbol{x}_n - \boldsymbol{x}_{n-1}\|_2^2 = 0,$$

for some minimiser $\boldsymbol{x}^\star \in \operatorname{Arg\,min}_{\boldsymbol{x}\in\mathbb{R}^N} \; \{\mathcal{J}(\boldsymbol{x}) := E_{\boldsymbol{y}}(\boldsymbol{x}) + \lambda\|\boldsymbol{x}\|_1\} \neq \emptyset$. In other words, the iterates $\{\boldsymbol{x}_n\}_{n\in\mathbb{N}}$ of APGD form a norm-convergent sequence, minimising the objective functional at a rate $o(1/n^2)$. In our practical implementation of Algorithm 2, we chose the step size $\tau$ as large as possible $\tau = 1/\beta$ and set $\mathfrak{d}$ to the value $\mathfrak{d} = 75$. The latter choice was motivated by the results reported in [52,53], which show significant practical acceleration for values of $\mathfrak{d}$ in the range $[50, 100]$. If performance is key, additional practical speedups can be achieved by adaptively restarting Algorithm 2 [53,54].

### 5.3. *Practical pseudo-differential operators*

One key insight of Theorem 2 is that the solutions of the FPBP problem (35) are $\mathscr{D}$-splines. As such, they inherit all their analytical properties from the zonal Green kernel associated to the pseudo-differential operator $\mathscr{D}$ used in the gTV regularisation term. For practical purposes, it is hence important to choose the pseudo-differential operator in agreement with the desired analytical properties of the solutions. In Example 1, we have for example introduced Sobolev operators $\mathscr{D}_\beta := [\mathrm{Id} - \Delta_{\mathbb{S}^{d-1}}]^\beta$, whose associated zonal Green kernels reproduce, for $\beta > (d-1)/2$, the *Sobolev spaces* [41, Lemma 5.5]

$$\mathscr{H}^\beta(\mathbb{S}^{d-1}) = \left\{ f \in \mathscr{S}'(\mathbb{S}^{d-1}) : \sum_{n\in\mathbb{N}} (1 + n(n+d-2))^\beta \sum_{m=1}^{N_d(n)} |\hat{f}_n^m|^2 < +\infty \right\}.$$

The latter are nested RKHSs, $\mathscr{H}^\gamma(\mathbb{S}^{d-1}) \subset \mathscr{H}^\beta(\mathbb{S}^{d-1}) \subset \mathscr{L}^2(\mathbb{S}^{d-1})$, $\forall \gamma \geq \beta > (d-1)/2$, containing functions with $\beta$-increasing degrees of smoothness. For example, a function $f \in \mathscr{H}^\beta(\mathbb{S}^{d-1})$ for $\beta \in \mathbb{N}$ is differentiable up to order $\beta$, with all its derivatives up to that order square-integrable. Sobolev operators seem hence particularly well-suited to enforce a certain degree of smoothness in the solutions of FPBP problems. Unfortunately, the Sobolev zonal Green kernel, given in (17), does not admit a convenient closed-form expression, hence making Sobolev operators – and consequently Sobolev splines – very cumbersome to work with in practice.

In this section, we hence design two pseudo-differential operators, namely the *Matérn* and *Wendland* operators, with convenient properties for practical purposes. Their zonal Green kernels admit indeed a simple closed-form expression and are well-localised in space. They have moreover a trivial nullspace (hence falling into the scope of Theorem 2) and their Fourier symbols are closely related to those of Sobolev operators, making them suitable for enforcing a certain degree of regularity in the solutions of FPBP problems. Unlike the standard differential operators from Example 1 which are defined via their Fourier symbols, the Matérn and Wendland operators are implicitly defined from their zonal Green kernel, obtained by restricting scaled *radial kernels* to the hypersphere [55,56]:

$$\psi_\beta^\epsilon(\langle \boldsymbol{r}, \boldsymbol{s}\rangle) := \frac{1}{\epsilon^n}\Psi_\beta\left(\epsilon^{-1}\|\boldsymbol{r}-\boldsymbol{s}\|_{\mathbb{R}^d}\right) = \frac{1}{\epsilon^n}\Psi_\beta\left(\epsilon^{-1}\sqrt{2 - 2\langle \boldsymbol{r}, \boldsymbol{s}\rangle}\right), \quad \forall(\boldsymbol{r}, \boldsymbol{s}) \in \mathbb{S}^{d-1}\times\mathbb{S}^{d-1}, \tag{54}$$

where $0 < \epsilon \leq 1$ is a *scale parameter* and $\Psi_\beta : \mathbb{R}_+ \to \mathbb{R},\ \beta > (d-1)/2$, is such that the kernel $\Psi_\beta(\|\boldsymbol{x} - \boldsymbol{y}\|_{\mathbb{R}^d})$, $\boldsymbol{x}, \boldsymbol{y} \in \mathbb{R}^d$, reproduces the Euclidean Sobolev space $\mathscr{H}^{\beta+1/2}(\mathbb{R}^d)$ w.r.t. a certain inner product – see [56, Section 2] for more details, and Sections 5.3.1 and 5.3.2 for examples of such radial basis functions. The resulting kernels (54) are zonal by construction, and can be expanded as [21, Chapter 3]:

$$\psi_\beta^\epsilon(\langle \boldsymbol{r}, \boldsymbol{s}\rangle) = \sum_{n\in\mathbb{N}} \frac{N_d(n)}{\mathfrak{a}_d\left(\hat{\psi}_\beta^\epsilon[n]\right)^{-1}} P_{n,d}(\langle \boldsymbol{r}, \boldsymbol{s}\rangle), \qquad \forall (\boldsymbol{r}, \boldsymbol{s}) \in \mathbb{S}^{d-1}\times\mathbb{S}^{d-1}, \tag{55}$$

where $\{\hat{\psi}_\beta^\epsilon[n]\}_{n\in\mathbb{N}}$ are the *Fourier-Legendre coefficients* of $\psi_\beta^\epsilon : [-1, 1] \to \mathbb{R}$. The latter can moreover be shown to verify [56, Lemma 2.1]:

$$c_1(1+\epsilon n)^{-2\beta} \leq \hat{\psi}_\beta^\epsilon[n] \leq c_2(1+\epsilon n)^{-2\beta}, \quad \forall n \geq 0, \tag{56}$$

for some positive constants $c_1, c_2$. Using (56), it is easy to show that $\hat{\psi}_\beta^\epsilon[n] > 0$, $\forall n \in \mathbb{N}$, and $\hat{\psi}_\beta^\epsilon[n] = \Theta(n^{-p})$ with $p = 2\beta > d-1$. From (55) and Proposition 3, we can hence identify the kernel (54) with the zonal Green kernel of the spline-admissible[32] pseudo-differential operator $\mathscr{D}_\beta^\epsilon$ with Fourier symbol $\{(\hat{\psi}_\beta^\epsilon[n])^{-1},\ n \in \mathbb{N}\}$. Still thanks to (56), it is moreover possible to show [55,56] that, for a given $\beta > (d-1)/2$, the norms $\|f\|_{\mathscr{D}_\beta^\epsilon} = \sum_{n\in\mathbb{N}} \hat{\psi}_\beta^\epsilon[n]^{-1} \sum_{m=1}^{N_d(n)} |\hat{f}_n^m|^2$, are all *equivalent* to the canonical *Sobolev norm*. The native RKHS

$$\mathscr{N}_{\mathscr{D}_\beta^\epsilon} = \left\{ f \in \mathscr{S}'(\mathbb{S}^{d-1}) : \sum_{n\in\mathbb{N}} \frac{1}{\hat{\psi}_\beta^\epsilon[n]} \sum_{m=1}^{N_d(n)} |\hat{f}_n^m|^2 < +\infty \right\},$$

contains therefore exactly the same elements as the Sobolev space $\mathscr{H}^\beta(\mathbb{S}^{d-1})$: $\mathscr{N}_{\mathscr{D}_\beta^\epsilon} = \mathscr{H}^\beta(\mathbb{S}^{d-1})$, $\forall \beta > (d-1)/2,\ \epsilon \in ]0, 1]$. In conclusion, the zonal Green kernels $\psi_\beta^\epsilon$ in (54) reproduce the Sobolev space $\mathscr{H}^\beta(\mathbb{S}^{d-1})$ when the latter is equipped with the inner product: $\langle h, g\rangle_{\mathscr{D}_\beta^\epsilon} = \sum_{n=0}^{+\infty} \hat{\psi}_\beta^\epsilon[n]^{-1} \left[\sum_{m=1}^{N_d(n)} \hat{h}_n^m \overline{\hat{g}_n^m}\right]$, and can hence be used as a replacement for the Sobolev zonal Green kernel to build practical spherical splines. In the subsequent sections, we give examples of functions $\Psi_\beta$ in (54), yielding respectively the *Matérn* and *Wendland* kernels, depicted in Fig. 3.

**Remark 14** *(About the scale parameter $\epsilon$).* For a fixed $\beta > (d-1)/2$, we have seen that the kernels $\psi_\beta^\epsilon$ for $0 < \epsilon \leq 1$ all reproduce the Sobolev space $\mathscr{H}^\beta(\mathbb{S}^{d-1})$. As such, one could question the relevancy of this parameter in the construction (54) proposed in [55,56]. Such doubts are however dispelled when considering approximation errors made by projecting functions from $\mathscr{H}^\beta(\mathbb{S}^{d-1})$ into specific spline spaces $\mathfrak{S}_{\mathscr{D}_\beta^\epsilon}(\mathbb{S}^{d-1}, \Xi_M)$ with fixed knot sets $\Xi_M \subset \mathbb{S}^{d-1}$. Indeed, it was shown in [55,56] that the approximation error is proportional to the quantity $(\Theta_{\Xi_M}/\epsilon)^\beta$, where $\Theta_{\Xi_M} > 0$ is the nodal width of $\Xi_M$ defined in (42) page 24. As such, choosing $\epsilon$ at least as large as the nodal width $\Theta_M$ helps in reducing the approximation error.

### *5.3.1. Matérn operators*

Matérn operators are obtained by choosing $\Psi_\beta$ in (54) as

$$\Psi_\beta(r) := \frac{\Gamma(\nu(\beta))}{2\Gamma(\nu(\beta)+1)} S_{\nu(\beta)}(r), \qquad \forall r \in \mathbb{R}_+. \tag{57}$$

In (57), $\nu(\beta) = \beta - (d-1)/2$, $\beta > (d-1)/2$, $\Gamma$ is the *Gamma function* and $S_\nu : \mathbb{R}_+ \to \mathbb{R}$ denotes the *Matérn function*, defined as [57, Chapter 4, p. 84]

$$S_\nu(r) := \frac{2^{1-\nu}}{\Gamma(\nu)} \left(\sqrt{2\nu} r\right)^\nu K_\nu\left(\sqrt{2\nu} r\right), \quad \forall r > 0,$$

where $K_\nu$ is the *modified Bessel function of the second kind* [58, Section 9.6]. It can be shown that the radial kernel (57) reproduces indeed the Sobolev spaces $\mathscr{H}^{\beta+1/2}(\mathbb{R}^d)$ [59, Theorem 6.13]. For half integers $\nu = p + 1/2$ with $p \in \mathbb{N}$, it is possible to write the Matérn function as the product of an exponential and a polynomial of order $p$ [57, Eq. 4.16], [60, Eq. 10.2.15]:

$$S_{p+1/2}(r) = \exp\left(-\sqrt{2p+1}r\right) \frac{p!}{(2p)!} \sum_{i=0}^{p} \frac{(p+i)!}{i!(p-i)!} \left(\sqrt{8p+4}r\right)^{p-i}, \qquad \forall r \in \mathbb{R}_+.$$

[32] Since $p = 2\beta > d-1$, we have indeed from Proposition 4 that $\mathscr{D}_\beta^\epsilon$ is spline-admissible.

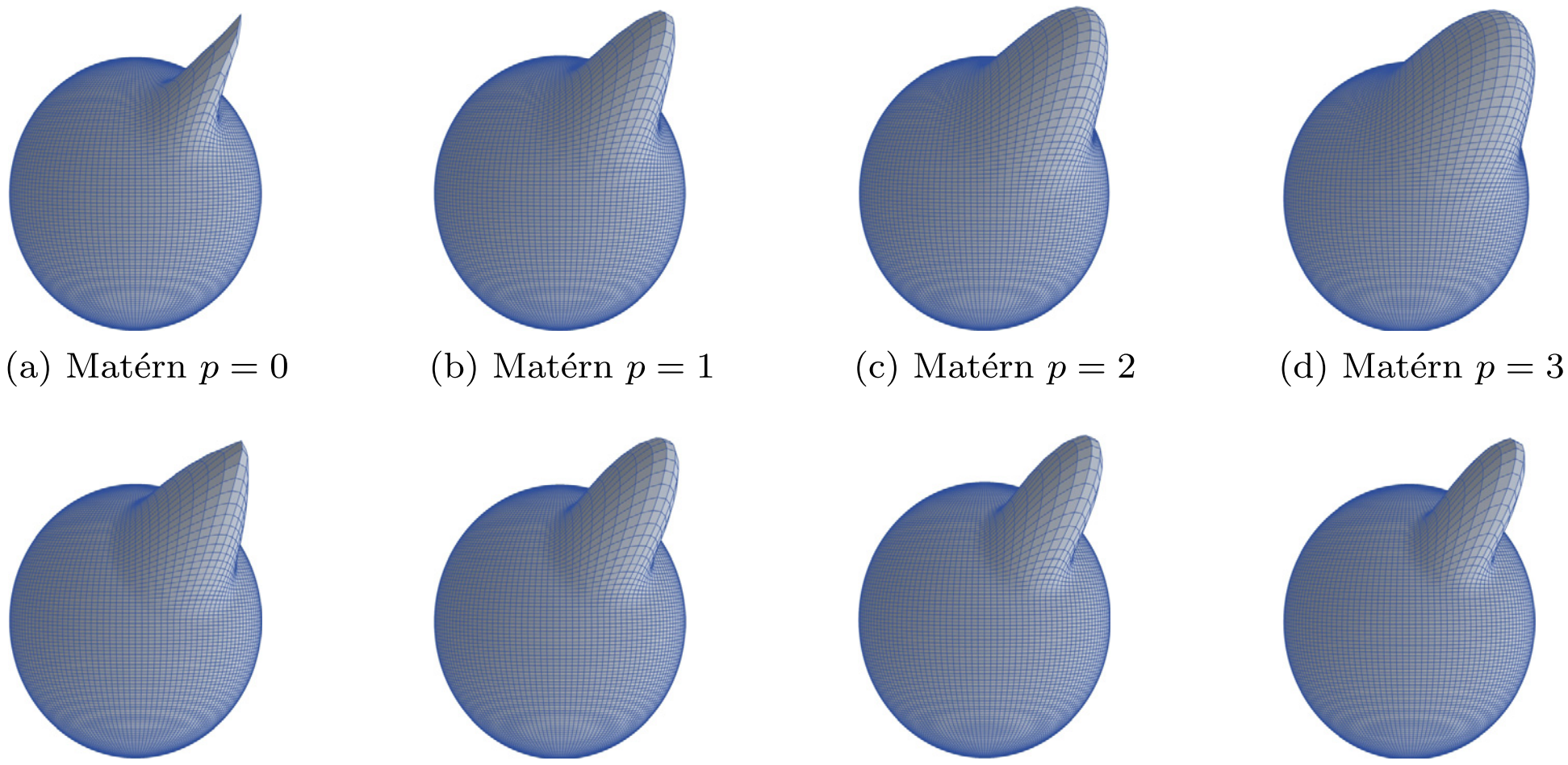

(a) Matérn $p = 0$ (b) Matérn $p = 1$ (c) Matérn $p = 2$ (d) Matérn $p = 3$

(e) Wendland $k = 0$ (f) Wendland $k = 1$ (g) Wendland $k = 2$ (h) Wendland $k = 3$

**Fig. 3.** Traces of Matérn and Wendland Green kernels $\psi_\beta^\epsilon(\langle \cdot, \boldsymbol{s}\rangle)$ (Figs. 3a to 3d and Figs. 3e to 3h respectively) for $\beta = p + 3/2$, $\epsilon = 0.1$, $p = 0, 1, 2, 3$ and focus direction $\boldsymbol{s} = (1, -1, 1)/\sqrt{3}$.

Matérn functions for $\nu \in \{1/2, 3/2, 5/2, 7/2\}$ are listed in [41, Figure 8.1]. In the limit $\nu \to \infty$, the Matérn function converges[33] towards the *Gaussian function* [57, Chapter 4, p. 84]. The zonal Green kernels $\psi_\epsilon^\beta$ of Matérn operators are hence bell-shaped functions, sharply decaying away from their centre[34] (see Fig. 3).

### *5.3.2. Wendland operators*

Wendland operators are obtained by choosing $\Psi_\beta$ in (54) as

$$\Psi_\beta(r) := \phi_{d,k(\beta)}(r), \qquad \forall r \in \mathbb{R}_+,$$

where $k(\beta) = \beta - d/2 \in \mathbb{N}$ and $\phi_{d,k} : \mathbb{R}_+ \to \mathbb{R}$, $k \in \mathbb{N}$ are the *Wendland's functions.* The latter are constructed by repeatedly applying and integral operator $\mathcal{I}$ to *Askey's truncated power functions* $\phi_l$:

$$\phi_{d,k}(r) := (\mathcal{I}^k \phi_l)(r), \quad k \in \mathbb{N}, l := \lfloor d/2 \rfloor + k + 1, \quad \phi_l(r) := (1-r)_+^l, \quad a_+ := \max(a, 0),$$

where $\mathcal{I}$ is given by: $(\mathcal{I}\phi)(r) = \int_r^{+\infty} t\phi(t)\, dt$, $r \geq 0$. It can be shown [14] that the Wendland's functions can be represented as:

$$\phi_{d,k}(r) = (1-r)_+^{l+k} p_{k,l}(r), \quad r \geq 0,$$

where $p_{k,l}$ is a polynomial of degree $k$ whose coefficients depend on $l$. These functions are *compactly supported* piecewise polynomials with support $[0, 1]$ which yield positive definite radial kernels in $\mathbb{R}^d$ with minimal degree and prescribed smoothness [14,61]. They have been introduced by Wendland [62] in the context of high-dimensional approximation. For $d \geq 3$, Wenland's radial kernels $\phi_{d,k(\beta)}(\|\boldsymbol{x}-\boldsymbol{y}\|_{\mathbb{R}^d})$, $\boldsymbol{x},\boldsymbol{y} \in \mathbb{R}^d$ were moreover proven [14,61] to reproduce Sobolev spaces of the form $\mathscr{H}^{\beta+1/2}(\mathbb{R}^d)$. In the case $d=2$, similar results can be obtained via a generalisation of the Wendland's functions called the *missing Wendland's functions* [14]. Examples of Wendland zonal Green kernel and their associated Sobolev spaces are provided in [41, Figure 8.2] for $d=3$.

[33] For practical purposes, $\nu \geq 7/2$ yield Matérn functions almost indistinguishable from the Gaussian function [57, Chapter 4, p. 84].

[34] Using the Gaussian approximation, the function $\psi_\epsilon^\beta(\cos(\theta))$ can be shown to have approximate support size $2\arcsin(5\epsilon/2)$.

#### 5.3.3. Sparse Gram matrices

In many experimental setups, the rapid decays of the Matérn and Wendland Green kernels cause the Gram matrices $\boldsymbol{G}$ in Theorem 3 to be *sparse*, allowing them to be conveniently implemented as such in Algorithms 1 and 2. For example, consider Theorem 3 in the context of the pseudo-differential operator associated to the Wendland kernel $\psi^{\epsilon}_{\beta}$ and sampling functionals $\{\varphi_l = \delta_{\boldsymbol{\rho}_l}, l = 1,\dots,L\}$ with sampling directions $\{\boldsymbol{\rho}_l,\ l = 1,\dots,L\} \subset \mathbb{S}^{d-1}$. Then, the entries of the Gram matrix $\boldsymbol{G} \in \mathbb{R}^{L\times N}$ are given by $G_{ln} = \psi^{\epsilon}_{\beta}(\langle \boldsymbol{\rho}_l, \boldsymbol{r}_n\rangle)$, $l = 1,\dots,L$, $n = 1,\dots,N$. It is then easy to see that, for $\epsilon$ small enough and the point sets $\{\boldsymbol{\rho}_l,\ l = 1,\dots,L\}$ and $\{\boldsymbol{r}_n,\ n = 1,\dots,N\}$ reasonably well distributed over $\mathbb{S}^{d-1}$, most of the entries of $\boldsymbol{G}$ are null. This behaviour extends to many spatially-localised measurement processes such as local averages or filtrations.

## 6. Test cases

In this section, we test our FPBP spherical approximation framework on two real-life spherical approximation problems originating from the fields of meteorology and forestry respectively. In both applications, various sampling functionals and data-fidelity functionals are investigated, demonstrating the versatility and genericity of our approximation framework. Interactive versions of the spherical maps provided in this section are moreover available at the links

https://matthieumeo.github.io/temperature_anomalies.html, & https://matthieumeo.github.io/fire_density.html

respectively.

### 6.1. Sea surface temperature anomalies

In this example, we propose to establish a global map of sea surface temperature anomalies for the month of January 2011. Such maps are used in environmental sciences to monitor global climate change as well as manage the population of marine species and ecosystems particularly sensitive to fluctuations in water temperature. The data consists in 6745 simulated anomalies sampled at various points across the globe by drifting floats of the ARGO fleet [33,34], and corrupted by Gaussian white noise. The map is obtained by solving the FPBP problem (33) with the indicator function of an $\ell_2$-ball as cost functionals (legitimised by the Gaussian noise assumption). The latter being non smooth but proximable, we solve the discrete PBP problem resulting from the canonical discretisation of (33) by means of the PDS Algorithm 1. For comparison purposes, we moreover consider discretising the FPBP problem by means of the non canonical domain discretisation scheme described in [41, Section 2 of Chapter 6]. Finally, we show the benefits of gTV regularisation w.r.t. the generalised Tikhonov (gTikhonov) regularisation discussed in [41, Chapter 5].

#### 6.1.1. Background

Sea surface temperature is usually defined as the temperature of the one millimetre upper layer of the oceans, reflecting the thermal energy stored in the latter. Sea surface temperatures departing from long-term averages (typically 12 years) are called *temperature anomalies*. While some anomalies are transient and simply due to ocean circulation patterns (such as *El Niño* and *La Niña*), other persist over many years, and can hence be potential indicators of global climate changes [63]. Sea surface temperature anomalies are also very important in the monitoring and management of marine ecosystems particularly sensitive to water temperature fluctuations. For example, above-average sea water temperatures can result in *coral bleaching*, a phenomenon is suspected to be responsible of the disappearance of between 29 and 50% of the Great Barrier Reef in 2016 [64]. Similarly, high water temperatures are contributing factors to *harmful algal blooms*, which lead to oxygen depletion in natural waters, with disastrous consequences on marine life [65].

#### *6.1.2. Data description*

For this experiment, we obtained simulated sea surface temperature anomalies by sampling at 6745 locations the global map of sea surface temperature anomalies produced by *NASA's Aqua* satellite [66] in January 2011 [67]. These anomalies, which serve here as a ground truth, were derived by comparing the sea surface temperatures recorded in January 2011 by NASA's Aqua satellite to the 12-year-averaged historical data for the same month collected by the *Pathfinder* satellite [68] between 1985 and 1997. The resulting map is depicted in Fig. 4a. The 6745 sampling locations were chosen as the positions of all floats from the *Argo fleet* [33] during the month of January 2011, obtained at this link [69] and curated by the authors of [34]. Argo is an international program, initiated in the early 2000's, that uses 4000 drifting floats to monitor temperature, salinity and currents in the Earth's oceans. The samples were further polluted by Gaussian white noise with PSNR 10 dB. The resulting samples are plotted in Fig. 4b.

#### *6.1.3. Data model*

Let $f : \mathbb{S}^2 \to \mathbb{R}$ denote the sea surface temperature anomaly function defined at every location on the globe (modelled as the unit sphere $\mathbb{S}^2$). Since temperatures typically have *smooth variations* at the surface of the Earth, we assume $f$ to be an element of some Sobolev space $\mathscr{H}^\beta(\mathbb{S}^2)$, with $\beta > 1$. The $L = 6745$ measurements $\{y_1, \dots, y_L\} \subset \mathbb{R}$ correspond here to noisy anomaly records collected by the Argo floats across the globe at locations

$$\{\boldsymbol{p}_1, \dots, \boldsymbol{p}_L\} \subset \mathbb{S}^2.$$

Assuming a Gaussian white noise model, the float records can moreover be modelled as realisations of independent Gaussian random variables $\{Y_1, \dots, Y_L\}$, centred around the true temperature anomalies obtained by *ideal spatial sampling* of $f$:

$$Y_i \overset{\text{ind}}{\sim} \mathcal{N}(f(\boldsymbol{p}_i), \sigma^2), \tag{58}$$

where $\mathcal{N}$ denotes the Gaussian distribution and $\sigma^2 > 0$ is the (unknown) noise variance, assumed uniform. Note that we have

$$\mathbb{E}[Y_i] = f(\boldsymbol{p}_i) = \langle \delta_{\boldsymbol{p}_i} | f \rangle, \qquad i = 1, \dots, L,$$

which fits well in our generic data model (23) page 17, if we choose the sampling functionals as Dirac measures $\delta_{\boldsymbol{p}_i}$. Note moreover that spatial sampling is indeed well-defined for $f$ since for $\beta > 1$ the Sobolev space $\mathscr{H}^\beta(\mathbb{S}^2)$ is an RKHS.

#### *6.1.4. Methods*

*gTV regularisation.* We consider first recovering $f$ by means of the following FPBP problem:

$$f^\star \in \underset{f \in \mathcal{M}_{\mathscr{D}^\epsilon_{2.5}}(\mathbb{S}^2)}{\operatorname{Arg\,min}} \left\{ \iota_{\mathcal{B}_{2,\rho}} (\boldsymbol{y} - \boldsymbol{\Phi}(f)) \quad + \quad \|\mathscr{D}^\epsilon_{2.5} f\|_{TV} \right\}, \tag{59}$$

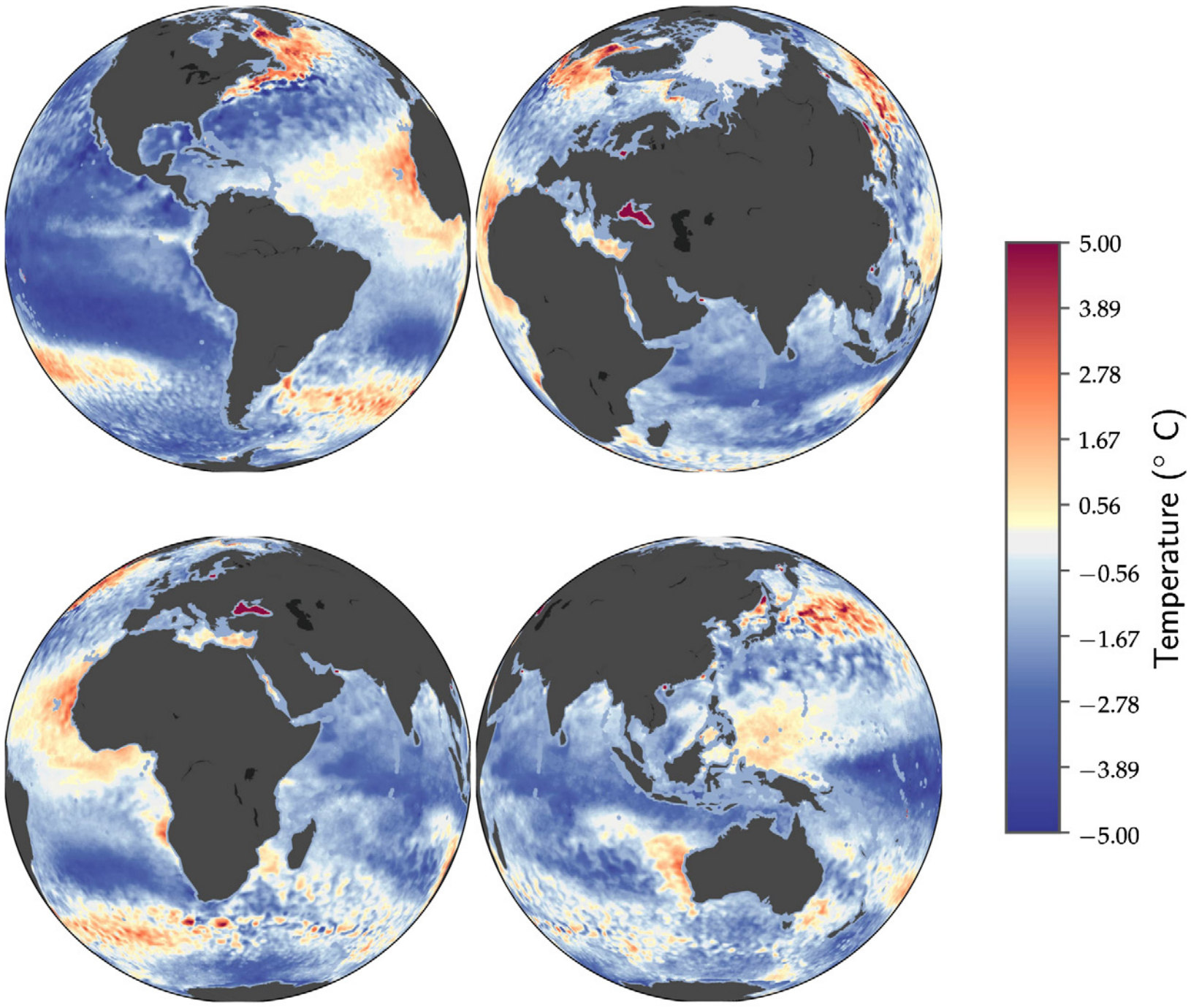

(a) Global map of sea surface temperature anomalies in January 2011 produced from *NASA's Aqua* satellite data.

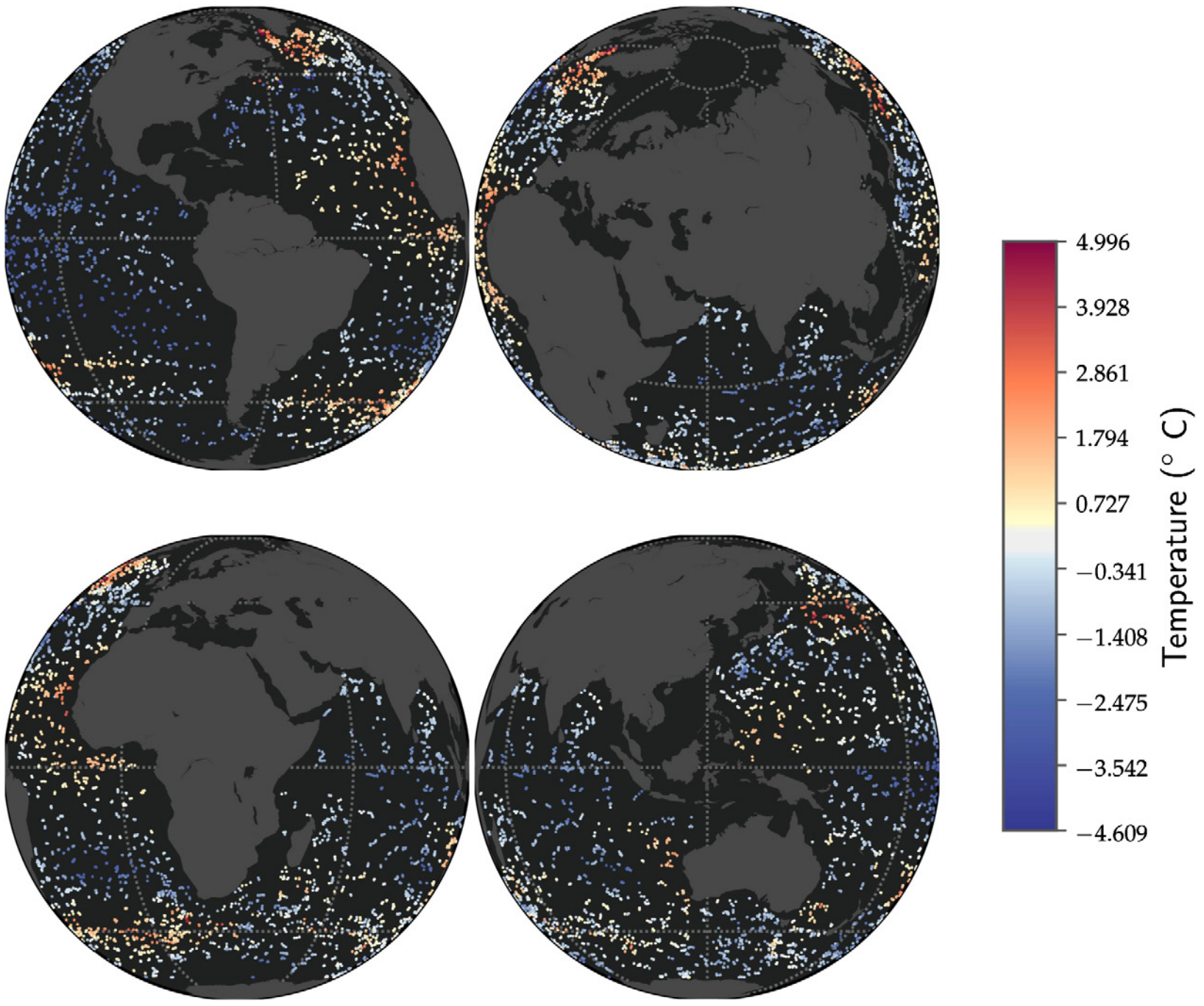

(b) Simulated anomalies recorded by Argo floats in January 2011. Float locations are marked by dots colored according to the recorded anomaly (red=warmer temperatures, blue=cooler temperatures).

**Fig. 4.** The data for the experiments in Section 6.1 consists in 6745 anomalies sampled from a global sea surface temperature map produced from *NASA's Aqua* satellite data in January 2011. The sample locations correspond to the locations of the Argo drifting floats during that month. (For interpretation of the colours in the figures, the reader is referred to the web version of this article.)

where:

- $\iota_{\mathcal{B}_{2,\rho}} : \mathbb{R}^L \to \{0\} \cup \{+\infty\}$ is the indicator function of the $\ell_2$-ball with radius $\rho = 0.5\% \times \|\boldsymbol{y}\|_2$. As explained in [41, Section 5.3 of Chapter 7], the indicator function in (63) defines, under the Gaussian noise model (58), a confidence region containing the true samples $\mathbb{E}[\boldsymbol{y}]$ with probability $1-\alpha$, for some $0 < \alpha < 1$ dependent on $\rho$.
- $\mathscr{D}^\epsilon_{2.5} : \mathscr{S}'(\mathbb{S}^2) \to \mathscr{S}'(\mathbb{S}^2)$ is the pseudo-differential operator associated to the Matérn zonal Green kernel with fixed scale $\epsilon \simeq 0.017$ – corresponding to an angular resolution[35] of approximately 4°:

$$\psi^\epsilon_{2.5}(\langle \boldsymbol{r}, \boldsymbol{s}\rangle) = S_{3/2}\left(\frac{\sqrt{2-2\langle \boldsymbol{r}, \boldsymbol{s}\rangle}}{\epsilon}\right) = \left[1 + \frac{\sqrt{2-2\langle \boldsymbol{r}, \boldsymbol{s}\rangle}}{\epsilon}\right] \exp\left(-\frac{\sqrt{2-2\langle \boldsymbol{r}, \boldsymbol{s}\rangle}}{\epsilon}\right), \quad \forall (\boldsymbol{r}, \boldsymbol{s}) \in \mathbb{S}^2 \times \mathbb{S}^2. \tag{60}$$

- $\boldsymbol{\Phi} : \mathcal{M}_{\mathscr{D}^\epsilon_{2.5}}(\mathbb{S}^2) \to \mathbb{R}^L$ is the *sampling operator* given by

$$\boldsymbol{\Phi}(f) = [\langle f|\delta_{\boldsymbol{p}_1}\rangle, \ldots, \langle f|\delta_{\boldsymbol{p}_L}\rangle], \quad \forall f \in \mathcal{M}_{\mathscr{D}^\epsilon_{2.5}}(\mathbb{S}^2).$$

  Note that $\boldsymbol{\Phi}$ is well-defined over $\mathcal{M}^\epsilon_{\mathscr{D}_{2.5}}(\mathbb{S}^2)$ since the Matérn kernel $\psi^\epsilon_{2.5}$ has continuous traces, and hence from Proposition 6, the predual $\mathscr{C}^\epsilon_{\mathscr{D}_{2.5}}(\mathbb{S}^2)$ contains all Dirac measures.

Note that since the cost functional in (59) is an indicator function, multiplying the gTV regularisation term $\|\mathscr{D}^\epsilon_{2.5} f\|_{TV}$ by a regularisation parameter $\lambda > 0$ would not change the solution set. Indeed, we trivially have:

$$(59) = \operatorname*{Arg\,min}_{f \in \mathcal{M}_{\mathscr{D}^\epsilon_{2.5}}} \{\|\mathscr{D}^\epsilon_{2.5} f\|_{TV} : \ \|\boldsymbol{y} - \boldsymbol{\Phi}(f)\|_2 \leq \rho\} = \operatorname*{Arg\,min}_{f \in \mathcal{M}_{\mathscr{D}^\epsilon_{2.5}}} \{\lambda\|\mathscr{D}^\epsilon_{2.5} f\|_{TV} : \ \|\boldsymbol{y} - \boldsymbol{\Phi}(f)\|_2 \leq \rho\}, \quad \forall \lambda > 0.$$

The trade-off between data-fidelity and sparsity of the generalised derivative of $f^\star$ is here controlled by the radius $\rho > 0$ of the $\ell_2$-ball: choosing $\rho$ small tends to favour solutions $f^\star$ compliant with the data $\boldsymbol{y}$ but with potentially non sparse generalised derivatives $\mathscr{D}^\epsilon_{2.5} f^\star$, while choosing $\rho$ large tends to favour solutions $f^\star$ with sparse generalised derivatives $\mathscr{D}^\epsilon_{2.5} f^\star$ but potentially non compliant with the data $\boldsymbol{y}$.

We consider discretising (59) by means of two discretisation schemes: the canonical quasi-uniform spline-based scheme described in Section 5.1 and the non canonical domain discretisation scheme described in [41, Section 2 of Chapter 6], which amounts to pixelating the sphere by means of the Fibonacci equidistributed spherical point set.

- **Canonical spline-based discretisation:** From the discussion in Section 5.1, solutions of the optimisation problem (59) can be approximated as quasi-uniform Matérn splines:

$$f^\star(\boldsymbol{r}) = \sum_{n=1}^{N} x^\star_n \, \psi^\epsilon_{2.5}(\langle \boldsymbol{r}, \boldsymbol{r}_n\rangle), \quad \forall \boldsymbol{r} \in \mathbb{S}^2,$$

  where $N = 7386$, $\Xi_N = \{\boldsymbol{r}_n, \, n = 1, \ldots, N\} \subset \mathbb{S}^2$ is a Fibonacci lattice (see Example 3) and $\boldsymbol{x}^\star = [x^\star_1, \ldots, x^\star_N] \in \mathbb{R}^N$ is some solution to the discrete optimisation problem:

$$\boldsymbol{x}^\star \in \operatorname*{Arg\,min}_{\boldsymbol{x} \in \mathbb{R}^N} \left\{\iota_{\mathcal{B}_{2,\rho}}(\boldsymbol{y} - \boldsymbol{G}\boldsymbol{x}) \quad + \quad \|\boldsymbol{x}\|_1\right\}. \tag{61}$$

[35] The angular resolution is measured here as the full width at half maximum (FWHM) of the Matérn kernel.

The matrix $\boldsymbol{G} \in \mathbb{R}^{L\times N}$ is moreover given by $G_{ln} = \psi^{\epsilon}_{2.5}(\langle \boldsymbol{p}_l, \boldsymbol{r}_n\rangle)$, $\forall (l,n) \in [\![1,L]\!] \times [\![1,N]\!]$. We solve (61) using the PDS Algorithm 1. Since the Matérn kernel is spatially localised, the matrix $\boldsymbol{G}$ is in practice sparse and is implemented as such in Algorithm 1 for computational and storage efficiency.

- **Non canonical domain discretisation:** For comparison purposes, we also discretising recovering (59) by means of the domain discretisation schemes described in [41, Section 2 of Chapter 6]. To this end, we consider the restriction $\boldsymbol{f} \in \mathbb{R}^N$ of $f$ to the discrete set of directions $\Xi_N = \{\boldsymbol{r}_n,\, n = 1,\dots,N\} \subset \mathbb{S}^2$, where $\Xi_N$ is the same Fibonacci lattice as before. This can be interpreted as pixelating the spherical domain. Due to its appealing conceptual simplicity, this is the prevailing approach in most applications. Unfortunately, such a discretisation *necessarily* incurs some (difficultly assessable) *approximation error*, sometimes even when the number of points composing the tessellation graph tends to infinity [70,12]. We recover $\boldsymbol{f}$ via the discrete domain counterparts of (59), given in this case by:

$$\boldsymbol{f}^\star \in \operatorname*{Arg\,min}_{\boldsymbol{f}\in\mathbb{R}^N} \left\{ \iota_{\mathcal{B}_{2,\rho}}(\boldsymbol{y} - \boldsymbol{J}\boldsymbol{f}) \quad + \quad \|\boldsymbol{D}\boldsymbol{f}\|_1 \right\}. \tag{62}$$

The entries of the *sensing matrix* $\boldsymbol{J} \in \mathbb{R}^{L\times N}$ are defined as $J_{ij} = \delta_{n_i j}$, $\forall i = 1,\dots,L$, $j = 1,\dots,N$, where $\delta_{ij}$ is the *Kronecker delta* and $n_i = \arg\min_{n=1,\dots,N} \|\boldsymbol{p}_i - \boldsymbol{r}_n\|_2$, $i = 1,\dots,L$. Roughly speaking, $\boldsymbol{J}$ corresponds to a discrete sampling matrix on the lattice $\Xi_N$, where the off-lattice sampling locations $\{\boldsymbol{p}_i, i = 1,\dots,L\}$ have been mapped to their closest neighbour in $\Xi_N$. The discrete pseudo-differential operator $\boldsymbol{D} \in \mathbb{R}^{N\times N}$ finally is chosen as the discrete Sobolev operator $\boldsymbol{D} = (\boldsymbol{I}_N + \boldsymbol{L})^2$, where $\boldsymbol{L}$ is the Laplacian of the spherical tessellation graph associated to the point set $\Xi_N$ (see [41, Section 2 of Chapter 6]). Note that $\boldsymbol{D} = (\boldsymbol{I}_N + \boldsymbol{L})^{2.5}$ would have been a more canonical choice since $\mathscr{D}^{\epsilon}_{2.5}$ in (59) is equivalent to the Sobolev operator $(\mathrm{Id} - \Delta_{\mathbb{S}^2})^{2.5}$ (see discussion in Section 5.3). However, computing $\boldsymbol{D} = (\boldsymbol{I}_N + \boldsymbol{L})^{2.5}$ requires computing the eigenvalue decomposition of the matrix $\boldsymbol{D} \in \mathbb{R}^{N\times N}$, which is often impossible in practice due to the size of the latter. Moreover, such a choice of discrete pseudo-differential operator would make [41, Algorithm 7.9] used to solve (62) much more computationally and memory intensive since the latter could no longer perform matrix-vector products involving $\boldsymbol{D}$ with the matrix free [41, Algorithm 7.11]. Indeed, this algorithm was designed for discrete pseudo-differential operators taking the form of polynomials of $\boldsymbol{L}$, and $\boldsymbol{D} = (\boldsymbol{I}_N + \boldsymbol{L})^{2.5}$ is *not* a polynomial in $\boldsymbol{L}$.

*gTikhonov regularisation.* For comparison purposes, we finally consider recovering $f$ by means of the following functional penalised Tikhonov (FPT) problem (see [41, Chapter 5] for details):

$$f^\star = \operatorname*{arg\,min}_{f\in\mathscr{H}_{\mathscr{D}^{\epsilon}_{2.5}}(\mathbb{S}^2)} \left\{ \iota_{\mathcal{B}_{2,\rho}}(\boldsymbol{y} - \boldsymbol{\Phi}(f)) \quad + \quad \|\mathscr{D}^{\epsilon}_{2.5} f\|_2^2 \right\}, \tag{63}$$

where $\boldsymbol{\Phi} : \mathscr{H}_{\mathscr{D}^{\epsilon}_{2.5}}(\mathbb{S}^2) \to \mathbb{R}^L$ is this time given by

$$\boldsymbol{\Phi}(f) = [\langle \delta_{\boldsymbol{p}_1} | f\rangle, \dots, \langle \delta_{\boldsymbol{p}_L} | f\rangle], \quad \forall f \in \mathscr{H}_{\mathscr{D}^{\epsilon}_{2.5}}(\mathbb{S}^2),$$

and $\mathscr{H}_{\mathscr{D}^{\epsilon}_{2.5}}(\mathbb{S}^2)$ is the Hilbert space defined in (28) (see [41, Chapter 5, Section 2.1] for a detailed analysis on this space). From [41, Theorem 5.3], the solution to the optimisation problem (63) is unique and given by:

$$f^\star(\boldsymbol{r}) = \sum_{l=1}^{L} x_l^\star \, \psi^{\epsilon}_{2.5} * \psi^{\epsilon}_{2.5}(\langle \boldsymbol{r}, \boldsymbol{p}_l\rangle), \quad \forall \boldsymbol{r} \in \mathbb{S}^2,$$

where $*$ denotes the spherical convolution operator[36] (see Definition 3) and $\hat{\boldsymbol{x}} = [\hat{x}_1, \ldots, \hat{x}_N] \in \mathbb{R}^N$ is the unique solution to the discrete optimisation problem (see [8, Section V]):

$$\boldsymbol{x}^\star = \arg\min_{\boldsymbol{x} \in \mathbb{R}^L} \left\{ \iota_{\mathcal{B}_{2,\rho}} (\boldsymbol{y} - \boldsymbol{H}\boldsymbol{x}) \quad + \quad \boldsymbol{x}^T \boldsymbol{H} \boldsymbol{x} \right\}. \tag{64}$$

Entries of the matrix $\boldsymbol{H} \in \mathbb{R}^{L \times L}$ are moreover given by $H_{lk} = \psi^\epsilon_{2.5} * \psi^\epsilon_{2.5}(\langle \boldsymbol{p}_l, \boldsymbol{p}_k \rangle)$, $\forall l, k \in [\![1, L]\!]$. We solve (64) using [41, Algorithm 7.3]. Again, since the Matérn kernel is spatially localised, the matrix $\boldsymbol{H}$ is in practice *sparse* and is implemented as such in the iterations of the numerical solver for computational and storage efficiency.

#### *6.1.5. Results*

The various estimates of the sea surface temperature anomaly function obtained by solving (59), (63) and (62) are provided in Figs. 5a, 5b and 6 respectively. The smoothing induced by the gTikhonov regularisation is clearly visible in Fig. 6. In contrast, the continuous and discrete gTV estimates in Figs. 5a and 5b capture far more of the fine fluctuations in the true anomaly map: see for example the eastern coast and southern tip of Africa, as well as the regions surrounding Greenland, Japan or the Indian ocean. This time, both estimates exhibit much more similar features. The discrete gTV estimate in Fig. 5b however appears rougher than the continuous gTV estimate in Fig. 5a due to the clearly visible pixelisation artefacts.

### *6.2. Wildfires*

In this example, we propose to establish global density maps of wildfires across the globe for the year 2016. Wildfire density maps allow scientists to better understand atmospheric chemistry and its impact on climate. The data used consists of fire counts recorded by NASA's Aqua and Terra satellites. The resolution of the raw data is moreover deliberately reduced by a factor of 3 by binning the counts in patches of angular size $\simeq 1.5^\circ \times 1.5^\circ$. The goal of this resolution reduction is to show that the lost resolution can be successfully recovered by spline-based approximation. Two density maps are obtained by solving with Algorithms 1 and 2 respectively two FPBP problems: one with a least-squares cost functional, and one with a KL-divergence cost functional – ideally suited for count data with Poisson-like distribution [41, Chapter 7].

#### *6.2.1. Background*

Home to 80% of terrestrial species, forests host most of Earth's biodiversity and contribute largely to its preservation [71]. Indeed, tree canopies play a crucial role in the regulation of the water cycle, creation of litter and exchange of energy between the ground and the atmosphere, which are all contributing factors to the overall good health of an ecosystem [1]. Changes to forest habitats can lead to the extinction of endangered species and disrupt the entire food chain equilibrium. But forests are more than animal shelters: they also protect humans from natural hazards such as floods or droughts [71]. In addition, forests represent natural and cost-efficient solutions for mitigating climate-change, and can provide 30% of the solution for keeping global warming below 2 °C [71].

In order to enlighten policy-makers and hopefully stem the current environmental crisis, it is hence crucial to monitor deforestation and understand its causes, such as agricultural conversion, commodity production, urbanisation, illegal logging or fires. Fires deserve perhaps a special attention as they contribute largely to the overall greenhouse gas emissions, with an estimated contribution of 30% to the net annual increase in the concentration of atmospheric carbon dioxide [72,73].

[36] Spherical convolutions between compactly supported kernels can be implemented efficiently with graph signal processing techniques, see [41, Chapter 8].

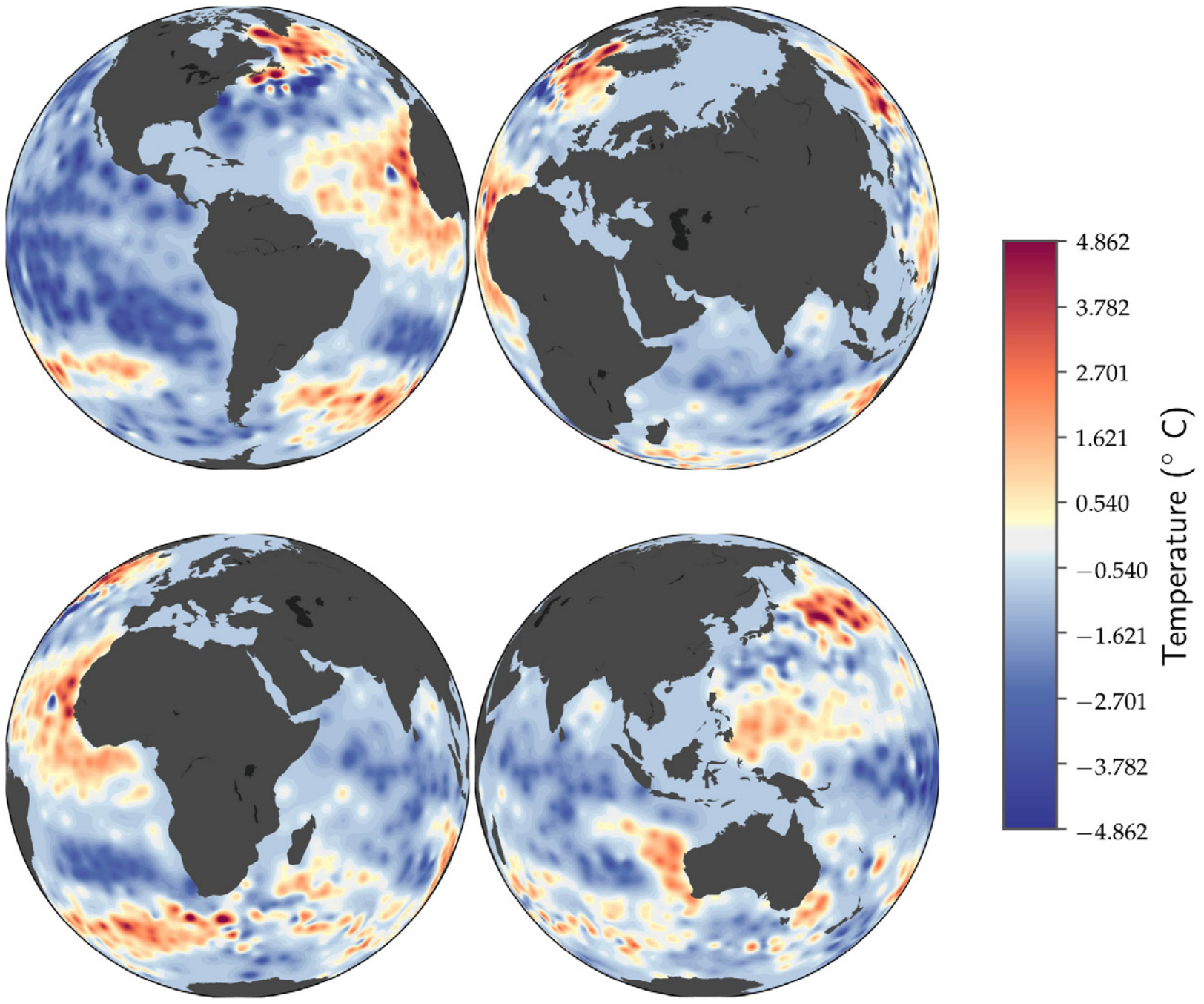

(a) Sea surface temperature anomaly function obtained by solving the FPBP problem (59) with gTV regularisation.

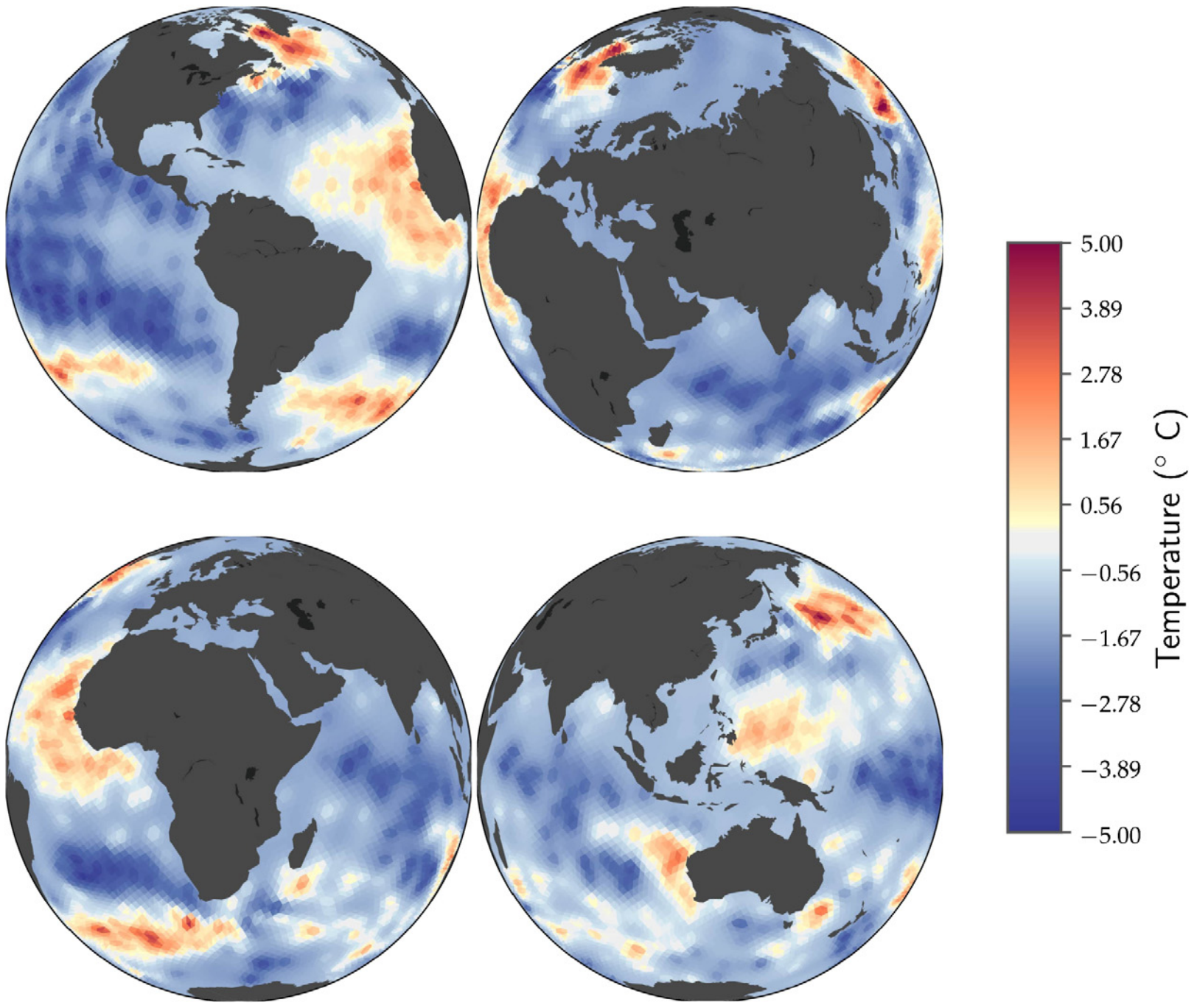

(b) Sea surface temperature anomaly function obtained by solving the discrete problem (62), with discrete gTV regularisation.

**Fig. 5.** Figs. 5a and 5b are the sea surface temperature anomaly functions obtained by solving the optimisation problems (59) and (62) respectively.

### *6.2.2. Data description*

For this experiment, we worked with the *Fire* [74] data product provided by NASA for the year 2016. The dataset, available at [72], provides monthly counts of active fires at a resolution of 0.1 degrees square. The counts are estimated from multispectral images captured by the *MODIS* aboard NASA's *Terra* and *Aqua* satellites [75]. For a better visual appreciation of the results, we aggregated the data from every month of 2016 (see Fig. 7a) and binned it to a lower resolution of approximately 1.5 degrees square (this corresponds to a reduction of resolution by a factor 3). We further corrupted the binned data with Poisson noise, a common noise model for count data. The processed data is displayed in Fig. 7b.

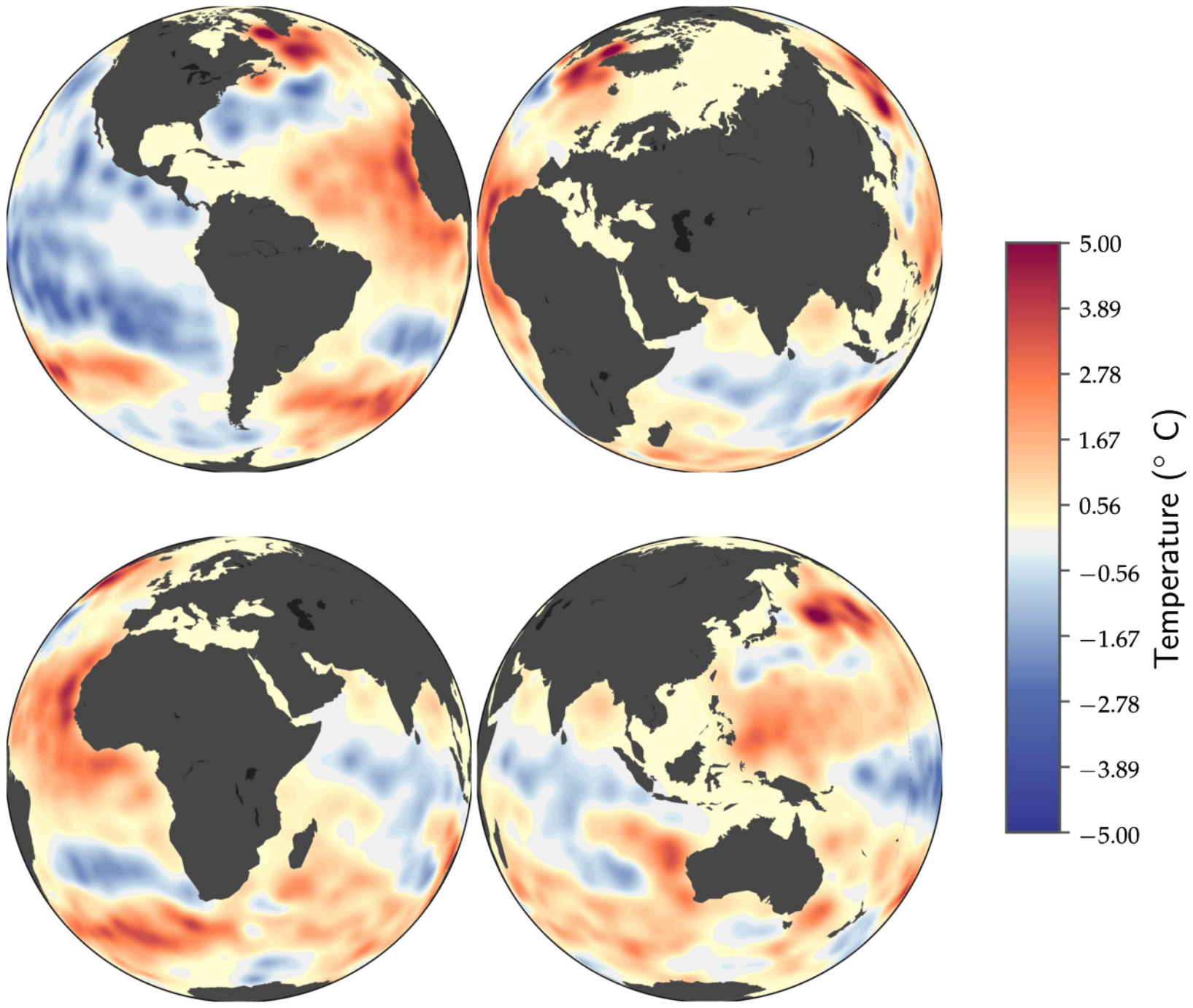

**Fig. 6.** Sea surface temperature anomaly function obtained by solving the FPT problem (63) with gTikhonov regularisation.

### *6.2.3. Data model*

In this experiment, one wishes to estimate the spatial density of fires $f$ at the surface of the Earth, using counts $\{y_1, \ldots, y_L\} \subset \mathbb{N}$ from non-overlapping equal-angle patches $\{B_1, \ldots, B_L\} \subset \mathbb{S}^2$ tiling the sphere. Since the data at hand consists of counts, a Poisson noise model is a suitable choice. This can be achieved by modelling fire locations as random occurrences of some *spatial Poisson point process* [76,77], often used in spatial statistics to model random spatial scattering of objects [78,79]. The distribution of a Poisson point process is entirely determined by its *intensity measure* $\Lambda \in \mathcal{M}(\mathbb{S}^2)$, which counts the expected number of objects (in this case fires) observed in a given region of the sphere. The sought spatial density of the Poisson process is then given – assuming it exists – by the *Radon-Nikodym derivative* [80] $f : \mathbb{S}^2 \to \mathbb{R}$ of $\Lambda$, also called *density* or *intensity* function of the point process. Similarly as in Section 6.1, we assume $f$ to be an element of the RKHS $\mathscr{H}^\beta(\mathbb{S}^2)$ for some $\beta > 1$. With such a formalism, the reported counts can then be seen as realisations of $L = 24000$ independent Poisson random variables $\{Y_1, \ldots, Y_L\}$:

$$Y_i \overset{\text{ind}}{\sim} \text{Poisson}\left(\lambda_i\right), \qquad i = 1, \ldots, L, \tag{65}$$

with rates $\lambda_i > 0$ given by:

$$\lambda_i = \int_{\mathbb{S}^2} f(\boldsymbol{r})\chi_{B_i}(\boldsymbol{r})d\boldsymbol{r} = \int_{B_i} f(\boldsymbol{r})d\boldsymbol{r}, \qquad i = 1, \ldots, L, \tag{66}$$

and where $\chi_{B_i} \in \mathscr{L}^2(\mathbb{S}^2)$ are the *characteristic functions* of the surveyed patches $\{B_i,\, i=1,\ldots,L\} \subset \mathbb{S}^2$. We can reinterpret the rates in (66) as generalised samples of $f$:

$$Y_i \overset{\text{ind}}{\sim} \text{Poisson}\left(\langle f, \chi_{B_i}\rangle_{\mathbb{S}^2}\right), \qquad i=1,\ldots,L,$$

where $\mathbb{E}[Y_i] = \langle f, \chi_{B_i}\rangle_{\mathbb{S}^2}$ , $i=1,\ldots,L$, hence yielding a data model falling into the scope of the generalised sampling framework (23).

#### 6.2.4. Methods

*KL-divergence cost function.* Since the data consists of counts, we consider recovering $f$ by means of the following FPBP problem:

$$f^\star \in \underset{f\in\mathcal{M}_{\mathscr{D}^\epsilon_{3,1}}(\mathbb{S}^2)}{\text{Arg min}} \left\{ D_{KL}(\boldsymbol{y}||\boldsymbol{\Phi}(f)) \quad + \quad \lambda\|\mathscr{D}^\epsilon_{3,1} f\|_{TV}\right\}, \tag{67}$$

where:

- $D_{KL}$ denotes the *generalised Kullback-Leibler divergence* defined in [41, Section 5.5 of Chapter 7]. As explained there, the KL-divergence cost function can be shown to be proportional to the negative log-likelihood of the data $\boldsymbol{y} = [y_1,\ldots,y_L] \in \mathbb{R}^L$ under the Poisson data model (65).
- $\boldsymbol{\Phi} : \mathcal{M}_{\mathscr{D}^\epsilon_{3,1}}(\mathbb{S}^2) \to \mathbb{R}^L$ is the sampling operator given by:

$$\boldsymbol{\Phi}(f) = [\langle f|\chi_{B_1}\rangle,\ldots,\langle f|\chi_{B_L}\rangle], \quad \forall f \in \mathcal{M}_{\mathscr{D}^\epsilon_{3,1}}(\mathbb{S}^2).$$

- $\mathscr{D}^\epsilon_{3,1} : \mathscr{S}'(\mathbb{S}^2) \to \mathscr{S}'(\mathbb{S}^2)$ is the pseudo-differential operator associated to the Wendland zonal Green kernel with a scale $\epsilon \simeq 0.026$ – corresponding again to an angular resolution[37] of approximately 1°:

$$\psi^\epsilon_{3,1}(\langle \boldsymbol{r},\boldsymbol{s}\rangle) = \phi_{3,1}\left(\frac{\sqrt{2-2\langle \boldsymbol{r},\boldsymbol{s}\rangle}}{\epsilon}\right) = \left(1-\frac{\sqrt{2-2\langle \boldsymbol{r},\boldsymbol{s}\rangle}}{\epsilon}\right)^4_+ \left(1+4\frac{\sqrt{2-2\langle \boldsymbol{r},\boldsymbol{s}\rangle}}{\epsilon}\right), \forall(\boldsymbol{r},\boldsymbol{s}) \in \mathbb{S}^2\times\mathbb{S}^2.$$

Note that since the sampling functions are square-integrable and $\mathscr{D}^\epsilon_{3,1}$ is invertible and with spectral growth order $p = 2.5 > (d-1)/2 = 1$, we can use Proposition 8 to show that $\{\chi_{B_i},\, i=1,\ldots,L\} \subset \mathscr{C}_{\mathscr{D}^\epsilon_{3,1}}(\mathbb{S}^2)$ and hence the sampling operator $\boldsymbol{\Phi}$ is indeed well defined.

Again, we consider approximating the solutions of (67) by means of quasi-uniform Wendland splines:

$$f^\star(\boldsymbol{r}) = \sum_{n=1}^N x^\star_n\, \psi^\epsilon_{3,1}(\langle \boldsymbol{r},\boldsymbol{r}_n\rangle), \quad \forall \boldsymbol{r} \in \mathbb{S}^2,$$

where $N = 210216$, $\Xi_N = \{\boldsymbol{r}_n,\, n=1,\ldots,N\} \subset \mathbb{S}^2$ is a Fibonacci lattice (see Example 3) and $\boldsymbol{x}^\star = [x^\star_1,\ldots,x^\star_N] \in \mathbb{R}^N$ is some solution to the discrete optimisation problem:

$$\boldsymbol{x}^\star \in \underset{\boldsymbol{x}\in\mathbb{R}^N}{\text{Arg min}} \left\{ D_{KL}(\boldsymbol{y}||\boldsymbol{G}\boldsymbol{x}) \quad + \quad \lambda\|\boldsymbol{x}\|_1\right\}. \tag{68}$$

From Theorem 3, the matrix $\boldsymbol{G} \in \mathbb{R}^{L\times N}$ is moreover given by $G_{ln} = (\psi^\epsilon_{3,1} * \chi_{B_l})(\boldsymbol{r}_n)$, $\forall(l,n) \in [\![1,L]\!]\times[\![1,N]\!]$. We solve (68) using Algorithm 1, once again leveraging the sparse structure of $\boldsymbol{G}$.

[37] The angular resolution is measured here as the full width at half maximum (FWHM) of the Wendland kernel.

*Choice of the regularisation parameter.* The regularisation parameter $\lambda$ in (67) and (68) was set to $\lambda \simeq 134 \times 10^3$. This value maximises the perceptual quality of the solution $f^\star$ over a collection of candidate values inspected manually.

*Quadratic cost function.* For sufficiently large rates, the Poisson distribution can be well approximated by a Gaussian distribution. This motivates the use of a least-squares data functional in (67), yielding:

$$f^\star \in \underset{f \in \mathcal{M}_{\mathscr{D}^\epsilon_{3,1}}(\mathbb{S}^2)}{\operatorname{Arg\,min}} \left\{ \|\boldsymbol{y} - \boldsymbol{\Phi}(f)\|_2^2 \quad + \quad \lambda \|\mathscr{D}^\epsilon_{3,1} f\|_{TV} \right\}. \tag{69}$$

The discrete optimisation problem (68) then becomes:

$$\boldsymbol{x}^\star \in \underset{\boldsymbol{x} \in \mathbb{R}^N}{\operatorname{Arg\,min}} \left\{ \|\boldsymbol{y} - \boldsymbol{G}\boldsymbol{x}\|_2^2 \quad + \quad \lambda \|\boldsymbol{x}\|_1 \right\}, \tag{70}$$

which can be solved efficiently using Algorithm 2.

*Choice of the regularisation parameter.* The value of $\lambda$ in (69) and (70) was set to $\lambda \simeq 29 \times 10^9$, which maximises the perceptual quality of the solution $f^\star$ over a collection of candidate values inspected manually. Note that this value is 6 orders of magnitude larger than the value of regularisation parameter considered in (67). This does not mean however that the regularisation term has more "weight" in (69) than in (67): the penalty strength must be assessed *relatively* to the natural scales of the cost functionals involved in (67) and (69). In this particular example, the scales of the KL-divergence and quadratic cost functions differ by 6 to 9 orders of magnitude: we have for example $D_{KL}(\boldsymbol{y}||\boldsymbol{G}\boldsymbol{1}) \simeq 5 \times 10^{13}$ while $\|\boldsymbol{y} - \boldsymbol{G}\boldsymbol{1}\|_2^2 \simeq 1.64 \times 10^{22}$. Therefore, the different values of $\lambda$ used in (69) and (67) yield comparable relative penalty strengths.

*6.2.5. Results*

The fire density maps estimated with both recovery strategies (67) and (69) are provided in Figs. 8a and 8b respectively. We observe that the recovered estimates (with KL-divergence and least-squares cost functions respectively) have a much higher resolution than the original corrupted binned counts in Fig. 7b, recovering almost the natural resolution of the unprocessed data in Fig. 7a. We observe however that the KL-divergence cost function seems to better recover regions with low intensity count, such as the Arabian Peninsula or Australia. In contrast, the least-squares data-fidelity functional has a tendency of yielding sparser density estimates, where all low intensity count regions are set to zero. This behaviour was already observed in image restoration under Poisson noise [81].

## 7. Conclusion

We have proposed a functional penalised basis pursuit approximation framework for functional inverse problems on the hypersphere. Unlike ad-hoc discrete methods traditionally favoured by practitioners, such problems present the advantage of being directly formulated in the continuous spherical domain, the natural domain for analogue spherical signals encountered in nature. In Theorem 2 we showed that, if regularised by means of gTV regularisation, functional inverse problems admitted *finite dimensional* solutions, which could hence be estimated in practice despite being defined over a *continuous* domain. More precisely, we showed in that the solutions were convex combinations of spherical $\mathscr{D}$-splines with *sparse* innovations, i.e. less than available measurements. This result inspired in Section 5.1 a canonical search space discretisation scheme with controlled approximation error (see Theorem 3). In Section 5.2, we moreover proposed algorithmic solutions adapted from the primal-dual splitting method and APGD to solve the discrete optimisation problems resulting from the canonical discretisation scheme. The proposed algorithms were shown

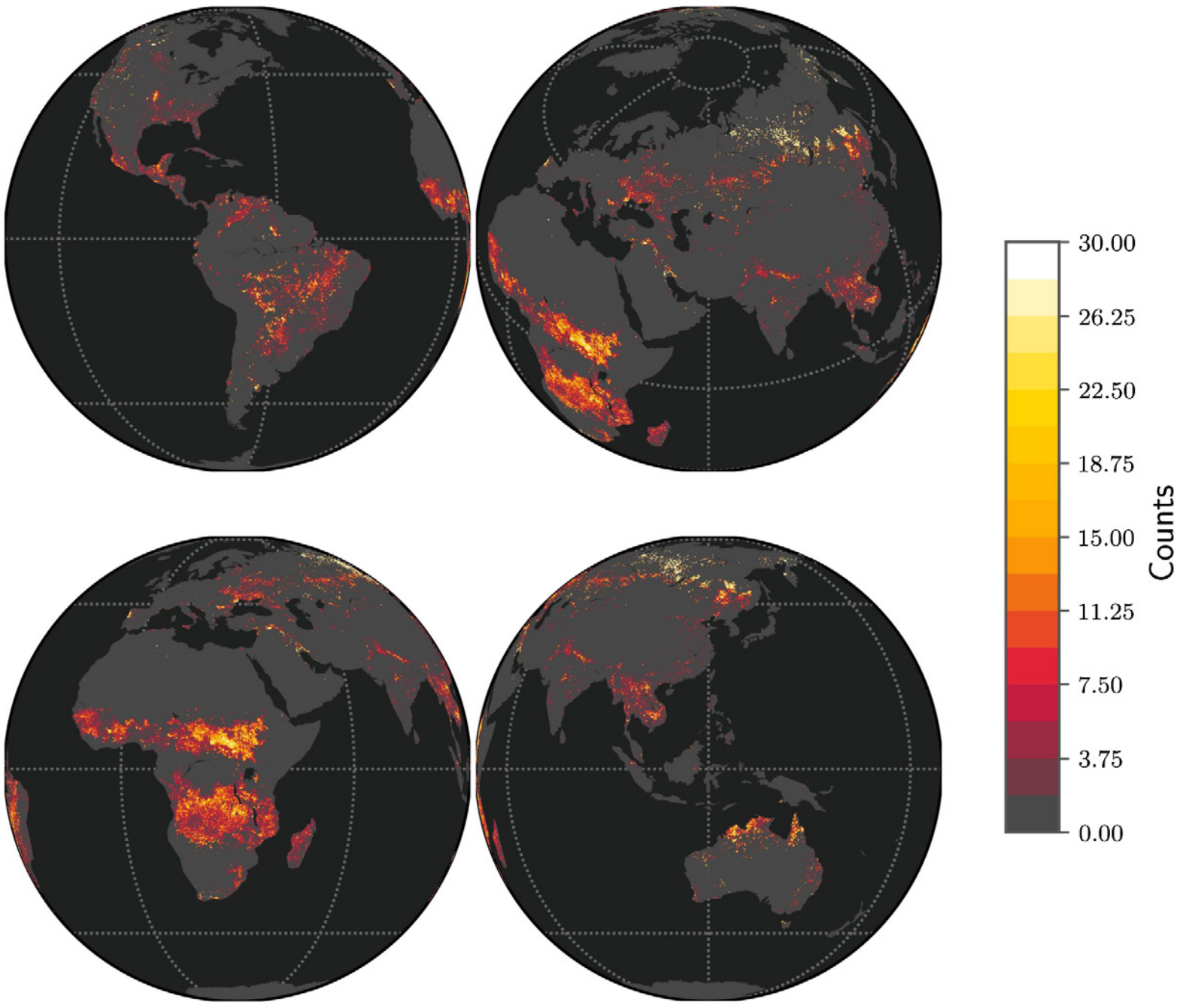

(a) Aggregated fire counts at full resolution (0.1 degree).

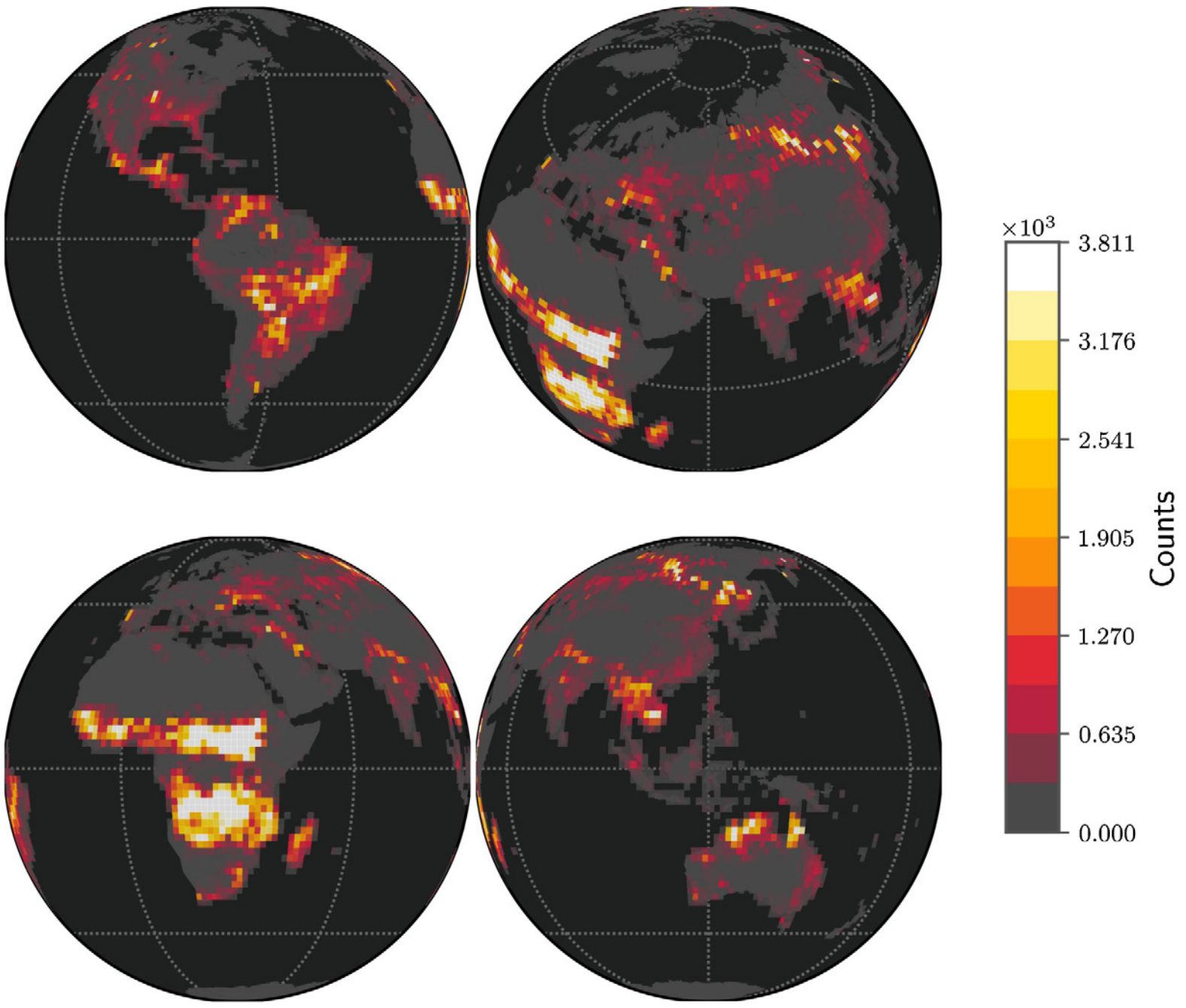

(b) Aggregated fire counts with reduced resolution (1.5 degree) and additional Poisson corruption.

**Fig. 7.** Aggregated fire counts for the year 2016 produced from MODIS data, a sensor aboard the *NASA's Terra/Aqua* satellites.

to be *computationally efficient*, provably *convergent* and compatible with *most common cost functionals* – including non-differentiable ones, such as the KL-divergence often used in the context of Poisson noise. In Section 5.3 finally, we introduced the last ingredient to our spherical approximation framework, namely the Wendland and Matérn splines, particularly convenient for practical purposes. We concluded in Sec-

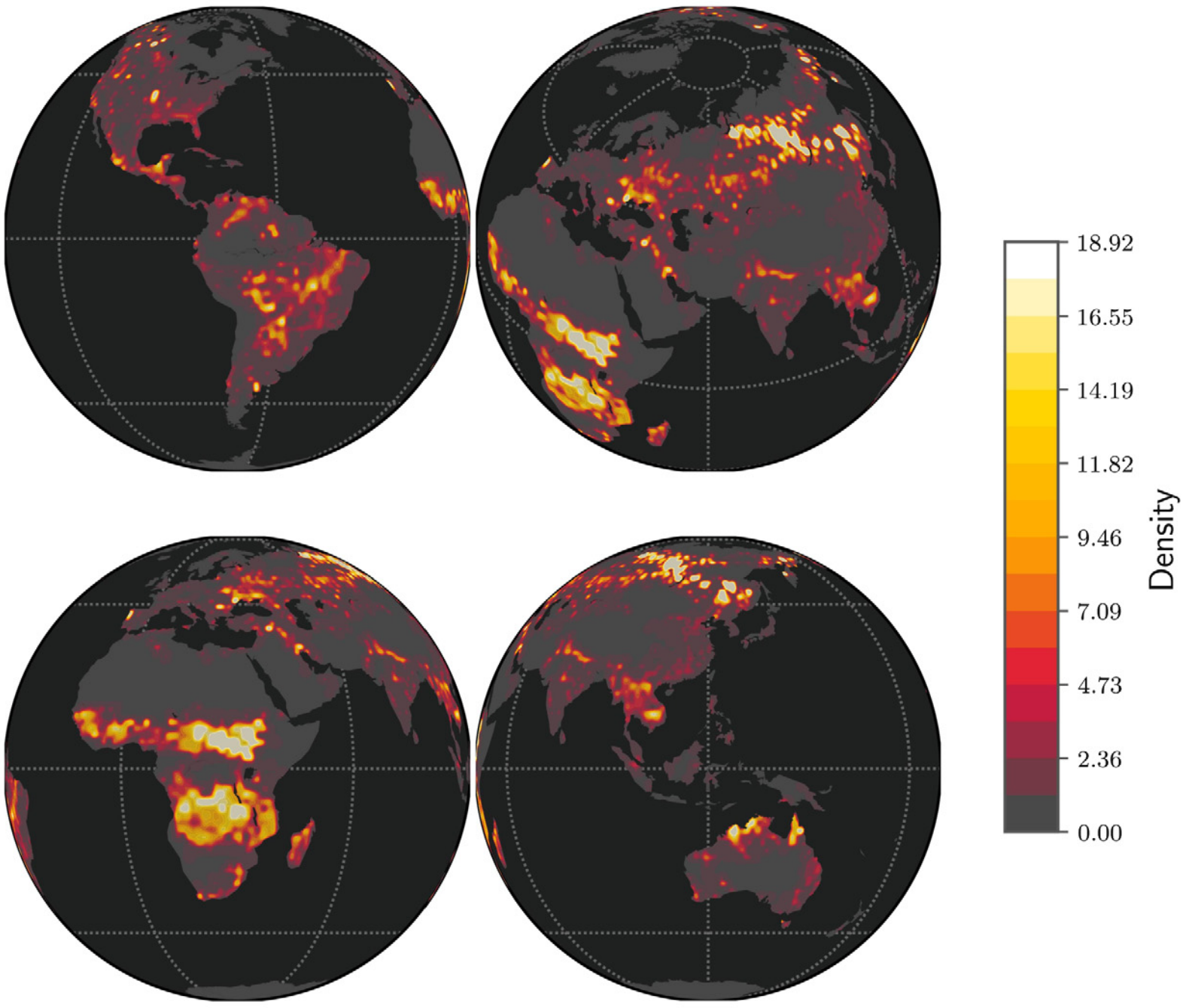

(a) Sparse spline approximation of the fire density function with KL-divergence data-fidelity.

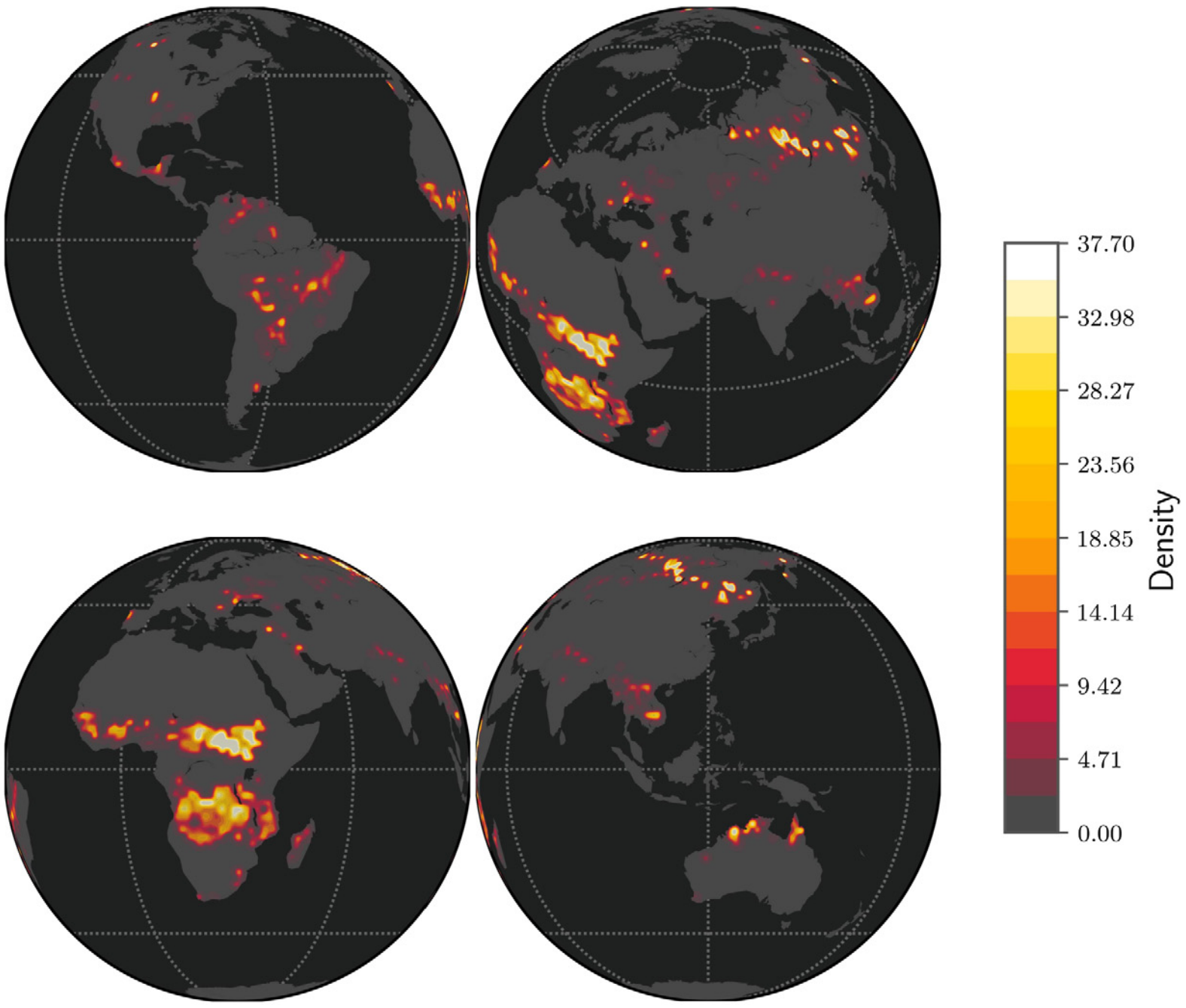

(b) Sparse spline approximation of the fire density function with least-squares data-fidelity.

**Fig. 8.** Figs. 8a and 8b are the fire density functions obtained by solving the optimisation problems (68) and (70) respectively.

tion 6 and tested our continuous-domain spherical approximation framework and novel algorithms on two datasets from the fields of meteorology and forestry. The sampling functionals, cost functions and regularisation strategies considered in each case were diverse, showing the versatility of both the theoretical framework and algorithmic solutions. In the meteorology example, we moreover illustrated the superiority of continuous-domain vs. discrete-domain recovery, both in terms of accuracy and resolution. This superiority was partially explained by the fact that continuous-domain methods could, unlike their discrete-domain counterparts, directly process the irregular spatial samples without resorting to *ad-hoc* gridding steps.

In conclusion, we hope that the contributions of this paper will spark interest among the community of practitioners, and give rise to the development of new reconstruction algorithms for spherical approximation problems. As future work, we plan on extending the framework to random, vector-valued or time-varying spherical fields. We also would like to investigate the sliding Frank-Wolfe algorithm proposed in [82] for solving FPBP problems without resorting to the spline-based discretisation scheme proposed in this paper.

**Acknowledgments**

The author is indebted to Michael Unser, Julien Fageot, and thank them for the insightful discussions, as well as their invaluable feedback on earlier versions of this manuscript. Their expert advice has largely contributed to improving the overall quality of this manuscript. The author is also grateful to Sepand Kashani, Victor Panaretos, Paul Hurley, and Martin Vetterli for their useful feedback, proofreading and numerous pieces of advice.

## Appendix A. Proofs of Section 4

**Proof of Proposition 2.** First, notice that since $|\hat{D}_n| = \Theta(n^p)$, we have in particular $|\hat{D}_n| = \Omega(n^p)$ and hence $\mathscr{D}^\dagger h \in \mathscr{S}(\mathbb{S}^{d-1})$ for all $h \in \mathscr{S}(\mathbb{S}^{d-1})$ (using similar arguments as for the proof of point 4 of [41, Proposition 4.1]). The compositions $\mathscr{D}\mathscr{D}^\dagger$ and $\mathscr{D}^\dagger\mathscr{D}$ are hence well-defined. Finally, we have from (9) and (12) that, for all $h \in \mathscr{S}(\mathbb{S}^{d-1})$

$$\mathscr{D}^\dagger\mathscr{D}\mathscr{D}^\dagger h = \sum_{n\notin\mathfrak{K}_\mathscr{D}} \frac{\hat{D}_n}{\hat{D}_n^2}\left[\sum_{m=1}^{N_d(n)} \hat{h}_n^m Y_n^m\right] = \sum_{n\notin\mathfrak{K}_\mathscr{D}} \frac{1}{\hat{D}_n}\left[\sum_{m=1}^{N_d(n)} \hat{h}_n^m Y_n^m\right] \mathscr{D}^\dagger h,$$

$$\mathscr{D}\mathscr{D}^\dagger\mathscr{D} h = \sum_{n\notin\mathfrak{K}_\mathscr{D}} \frac{\hat{D}_n^2}{\hat{D}_n}\left[\sum_{m=1}^{N_d(n)} \hat{h}_n^m Y_n^m\right] = \sum_{n\in\mathbb{N}} \hat{D}_n\left[\sum_{m=1}^{N_d(n)} \hat{h}_n^m Y_n^m\right] = \mathscr{D}h,$$

which shows that $\mathscr{D}^\dagger$ is a generalised inverse of $\mathscr{D}$. Moreover, we have that

$$\mathscr{D}\mathscr{D}^\dagger h = \sum_{n\notin\mathfrak{K}_\mathscr{D}} \frac{\hat{D}_n}{\hat{D}_n}\left[\sum_{m=1}^{N_d(n)} \hat{h}_n^m Y_n^m\right] = \sum_{n\notin\mathfrak{K}_\mathscr{D}} \left[\sum_{m=1}^{N_d(n)} \hat{h}_n^m Y_n^m\right] = \mathscr{D}\mathscr{D}^\dagger h.$$

Since both $\mathscr{D}$ and $\mathscr{D}^\dagger$ are self-adjoint, we have $(\mathscr{D}\mathscr{D}^\dagger)^* = \mathscr{D}^\dagger\mathscr{D} = \mathscr{D}\mathscr{D}^\dagger$ and $(\mathscr{D}^\dagger\mathscr{D})^* = \mathscr{D}\mathscr{D}^\dagger = \mathscr{D}^\dagger\mathscr{D}$, which shows that $\mathscr{D}^\dagger$ is actually the Moore-Penrose pseudo-inverse of $\mathscr{D}$ and concludes the proof. □

**Proof of Proposition 3.** Let $\boldsymbol{s} \in \mathbb{S}^{d-1}$. From the generalised spherical harmonic transform applied to $\Psi_{\boldsymbol{s}}^\mathscr{D}$ we have:

$$\Psi_{\boldsymbol{s}}^\mathscr{D} = \sum_{n\in\mathbb{N}}\sum_{m=1}^{N_d(n)} \langle \Psi_{\boldsymbol{s}}^\mathscr{D} | Y_n^m \rangle Y_n^m = \sum_{n\in\mathbb{N}\setminus\mathfrak{K}_\mathscr{D}}\sum_{m=1}^{N_d(n)} \frac{1}{\hat{D}_n} Y_n^m(\boldsymbol{s}) Y_n^m,$$

since from (14), we have $\langle \Psi_{\boldsymbol{s}}^\mathscr{D} | Y_n^m \rangle = \langle \mathscr{D}^\dagger \delta_{\boldsymbol{s}} | Y_n^m \rangle = \langle \delta_{\boldsymbol{s}} | \mathscr{D}^\dagger Y_n^m \rangle = \hat{D}_n^{-1} Y_n^m(\boldsymbol{s})$ if $n \notin \mathfrak{K}_\mathscr{D}$ and zero otherwise. We have hence, from the bilinearity of the Schwartz duality product and the addition theorem [21, Theorem 5.11]:

$$
\begin{aligned}
\langle \Psi_{\boldsymbol{s}}^{\mathscr{D}} | \varphi \rangle &= \sum_{n\in\mathbb{N}\setminus\mathfrak{K}_{\mathscr{D}}} \frac{1}{\hat{D}_n} \sum_{m=1}^{N_d(n)} Y_n^m(\boldsymbol{s}) \langle Y_n^m | \varphi \rangle = \sum_{n\in\mathbb{N}\setminus\mathfrak{K}_{\mathscr{D}}} \frac{1}{\hat{D}_n} \sum_{m=1}^{N_d(n)} Y_n^m(\boldsymbol{s}) \left\langle \varphi, Y_n^m \right\rangle_{\mathbb{S}^{d-1}} \\
&= \sum_{n\in\mathbb{N}\setminus\mathfrak{K}_{\mathscr{D}}} \frac{1}{\hat{D}_n} \sum_{m=1}^{N_d(n)} Y_n^m(\boldsymbol{s}) \int_{\mathbb{S}^{d-1}} \varphi(\boldsymbol{r}) \overline{Y_n^m}(\boldsymbol{r})\, d\boldsymbol{r} = \int_{\mathbb{S}^{d-1}} \varphi(\boldsymbol{r}) \left[ \sum_{n\in\mathbb{N}\setminus\mathfrak{K}_{\mathscr{D}}} \frac{1}{\hat{D}_n} \sum_{m=1}^{N_d(n)} Y_n^m(\boldsymbol{s}) \overline{Y_n^m}(\boldsymbol{r}) \right] d\boldsymbol{r} \\
&= \int_{\mathbb{S}^{d-1}} \varphi(\boldsymbol{r}) \left[ \sum_{n\in\mathbb{N}\setminus\mathfrak{K}_{\mathscr{D}}} \frac{N_d(n)}{\mathfrak{a}_d \hat{D}_n} P_{n,d}(\langle \boldsymbol{r}, \boldsymbol{s} \rangle) \right] d\boldsymbol{r} = \int_{\mathbb{S}^{d-1}} \varphi(\boldsymbol{r}) \psi_{\mathscr{D}}(\langle \boldsymbol{r}, \boldsymbol{s} \rangle)\, d\boldsymbol{r}, \qquad \forall \varphi \in \mathscr{S}(\mathbb{S}^{d-1}). \quad \square
\end{aligned}
$$

**Proof of Proposition 4.** Let $\boldsymbol{s} \in \mathbb{S}^{d-1}$ be fixed but arbitrary. We show that, under the assumptions of Proposition 4, the series

$$
\psi_{\mathscr{D}}(\langle \boldsymbol{r}, \boldsymbol{s} \rangle) = \sum_{n\in\mathbb{N}\setminus\mathfrak{K}_{\mathscr{D}}} \frac{N_d(n)}{\mathfrak{a}_d \hat{D}_n} P_{n,d}(\langle \boldsymbol{r}, \boldsymbol{s} \rangle), \quad \boldsymbol{r} \in \mathbb{S}^{d-1}, \tag{A.1}
$$

converges uniformly (w.r.t. the variable $\boldsymbol{r}$). Since every summand is continuous, we can then conclude that the limit $\psi_{\mathscr{D}}(\langle \boldsymbol{r}, \boldsymbol{s} \rangle)$ is continuous (see [21, Theorem 2.14]) and hence in particular pointwise defined – i.e. $\mathscr{D}$ is spline-admissible. To show that (A.1) is uniformly convergent, we consider its remainder for some $N > \max(\mathfrak{K}_{\mathscr{D}})$. Then, from the addition theorem [21, Theorem 5.11] and the Cauchy-Schwarz inequality we get, for each $\boldsymbol{r} \in \mathbb{S}^{d-1}$:

$$
\begin{aligned}
\left| \sum_{n=N}^{+\infty} \frac{N_d(n)}{\mathfrak{a}_d \hat{D}_n} P_{n,d}(\langle \boldsymbol{r}, \boldsymbol{s} \rangle) \right| &= \left| \sum_{n=N}^{+\infty} \frac{1}{\hat{D}_n} \sum_{m=1}^{N_d(n)} Y_n^m(\boldsymbol{s}) \overline{Y_n^m}(\boldsymbol{r}) \right| = \left| \sum_{n=N}^{+\infty} \sum_{m=1}^{N_d(n)} \left( \frac{Y_n^m(\boldsymbol{s})}{\operatorname{sgn}(\hat{D}_n)\sqrt{|\hat{D}_n|}} \right) \overline{\left( \frac{Y_n^m(\boldsymbol{r})}{\sqrt{|\hat{D}_n|}} \right)} \right| \\
&\leq \left| \sum_{n=N}^{+\infty} \frac{\sum_{m=1}^{N_d(n)} |Y_n^m(\boldsymbol{s})|^2}{|\hat{D}_n|} \right| \left| \sum_{n=N}^{+\infty} \frac{\sum_{m=1}^{N_d(n)} |Y_n^m(\boldsymbol{r})|^2}{|\hat{D}_n|} \right| = \left| \sum_{n=N}^{+\infty} \frac{N_d(n)}{\mathfrak{a}_d |\hat{D}_n|} \right|^2.
\end{aligned}
$$

Moreover, since $|\hat{D}_n| = \Theta(n^p)$ we have from (6) $N_d(n)|\hat{D}_n|^{-1} = \mathcal{O}(n^{d-2-p})$. Since $p > d-1 \Rightarrow d-2-p < -1$, the series $\sum_{n\in\mathbb{N}\setminus\mathfrak{K}_{\mathscr{D}}} \frac{N_d(n)}{\mathfrak{a}_d|\hat{D}_n|}$ is convergent and hence its remainder tends to zero. Therefore

$$
\left| \psi_{\mathscr{D}}(\langle \boldsymbol{r}, \boldsymbol{s} \rangle) - \sum_{\substack{n=0 \\ n\notin\mathfrak{K}_{\mathscr{D}}}}^{N-1} \frac{N_d(n)}{\mathfrak{a}_d \hat{D}_n} P_{n,d}(\langle \boldsymbol{r}, \boldsymbol{s} \rangle) \right| = \left| \sum_{n=N}^{+\infty} \frac{N_d(n)}{\mathfrak{a}_d \hat{D}_n} P_{n,d}(\langle \boldsymbol{r}, \boldsymbol{s} \rangle) \right| \leq \left| \sum_{n=N}^{+\infty} \frac{N_d(n)}{\mathfrak{a}_d |\hat{D}_n|} \right|^2 \overset{N\to+\infty}{\longrightarrow} 0.
$$

Moreover, since the bound is independent of $\boldsymbol{r}$ (and $\boldsymbol{s}$ as a matter of fact) the convergence is uniform, which achieves the proof. $\square$

**Proof of Proposition 5.** Consider the function $\mathfrak{s}' : \mathbb{S}^{d-1} \to \mathbb{C}$ defined as

$$
\mathfrak{s}'(\boldsymbol{r}) = \sum_{i=1}^{M} \alpha_i \psi_{\mathscr{D}}(\langle \boldsymbol{r}, \boldsymbol{r}_i \rangle) \quad + \quad \sum_{n\in\mathfrak{K}_{\mathscr{D}}} \sum_{m=1}^{N_d(n)} \langle \mathfrak{s} | Y_n^m \rangle Y_n^m(\boldsymbol{r}), \quad \boldsymbol{r} \in \mathbb{S}^{d-1}.
$$

Notice that since $\mathscr{D}$ is spline-admissible, the functions $\psi_{\mathscr{D}}(\langle \cdot, \boldsymbol{r}_i \rangle)$ are ordinary functions which are hence bounded on $\mathbb{S}^{d-1}$ and hence $\psi_{\mathscr{D}}(\langle \cdot, \boldsymbol{r}_i \rangle) \in \mathscr{L}^2(\mathbb{S}^{d-1})$, which implies in turn that $\mathfrak{s}' \in \mathscr{L}^2(\mathbb{S}^{d-1})$. We can hence interpret $\mathfrak{s}'$ as an element of $\mathscr{S}'(\mathbb{S}^{d-1})$ with pointwise definition: $\langle \mathfrak{s}' | \varphi \rangle := \langle \varphi, \mathfrak{s}' \rangle_{\mathbb{S}^{d-1}}$, $\forall \varphi \in \mathscr{S}(\mathbb{S}^{d-1})$. We now show that $\mathfrak{s} = \mathfrak{s}'$, i.e. $\langle \mathfrak{s} | \varphi \rangle = \langle \mathfrak{s}' | \varphi \rangle$, $\forall \varphi \in \mathscr{S}(\mathbb{S}^{d-1})$. First, we write $\mathscr{S}(\mathbb{S}^{d-1}) = \mathcal{R}(\mathscr{D}) \oplus \mathcal{N}(\mathscr{D})$ such that every element $\varphi$ of $\mathscr{S}(\mathbb{S}^{d-1})$ can be written as $\varphi = \mathscr{D}h + \eta$, with $(h, \eta) \in \mathcal{R}(\mathscr{D}) \times \mathcal{N}(\mathscr{D})$. We have hence

$$\begin{aligned}\langle \mathfrak{s}|\varphi\rangle &= \langle \mathfrak{s}|\mathscr{D}h\rangle + \langle \mathfrak{s}|\eta\rangle = \langle \mathscr{D}\mathfrak{s}|h\rangle + \sum_{n\in\mathfrak{K}_{\mathscr{D}}}\sum_{m=1}^{N_d(n)} \hat{\eta}_n^m \langle \mathfrak{s}|Y_n^m\rangle \\ &= \left\langle \sum_{i=1}^{N} \alpha_i \delta_{\boldsymbol{r}_i} \middle| h \right\rangle + \sum_{n\in\mathfrak{K}_{\mathscr{D}}}\sum_{m=1}^{N_d(n)} \hat{\eta}_n^m \langle \mathfrak{s}|Y_n^m\rangle = \sum_{i=1}^{N}\alpha_i h(\boldsymbol{r}_i) + \sum_{n\in\mathfrak{K}_{\mathscr{D}}}\sum_{m=1}^{N_d(n)} \langle \mathfrak{s}|Y_n^m\rangle \hat{\eta}_n^m.\end{aligned}$$

Similarly we have from Proposition 3

$$\begin{aligned}\langle \mathfrak{s}'|\varphi\rangle = \langle \varphi, \mathfrak{s}'\rangle_{\mathbb{S}^{d-1}} &= \sum_{i=1}^{M} \alpha_i \langle \varphi, \psi_{\mathscr{D}}(\langle\cdot, \boldsymbol{r}_i\rangle)\rangle_{\mathbb{S}^{d-1}} + \sum_{n\in\mathfrak{K}_{\mathscr{D}}}\sum_{m=1}^{N_d(n)} \langle \mathfrak{s}|Y_n^m\rangle \langle \varphi, Y_n^m\rangle_{\mathbb{S}^{d-1}} \\ &= \sum_{i=1}^{M} \alpha_i \langle \Psi_{\boldsymbol{s}}^{\mathscr{D}}|\varphi\rangle + \sum_{n\in\mathfrak{K}_{\mathscr{D}}}\sum_{m=1}^{N_d(n)} \langle \mathfrak{s}|Y_n^m\rangle \langle \eta, Y_n^m\rangle_{\mathbb{S}^{d-1}} \\ &= \sum_{i=1}^{M} \alpha_i \langle \delta_{\boldsymbol{r}_i}|\mathscr{D}^\dagger\mathscr{D}h\rangle + \sum_{n\in\mathfrak{K}_{\mathscr{D}}}\sum_{m=1}^{N_d(n)} \langle \mathfrak{s}|Y_n^m\rangle \hat{\eta}_n^m \\ &= \sum_{i=1}^{M} \alpha_i \langle \delta_{\boldsymbol{r}_i}|h\rangle + \sum_{n\in\mathfrak{K}_{\mathscr{D}}}\sum_{m=1}^{N_d(n)} \langle \mathfrak{s}|Y_n^m\rangle \hat{\eta}_n^m \\ &= \sum_{i=1}^{M} \alpha_i h(\boldsymbol{r}_i) + \sum_{n\in\mathfrak{K}_{\mathscr{D}}}\sum_{m=1}^{N_d(n)} \langle \mathfrak{s}|Y_n^m\rangle \hat{\eta}_n^m,\end{aligned}$$

and hence we have indeed $\langle \mathfrak{s}|\varphi\rangle = \langle \mathfrak{s}'|\varphi\rangle \ \forall \varphi \in \mathscr{S}(\mathbb{S}^{d-1})$ as claimed. Equation (21) then follows trivially from the fact that the summations involving $\mathfrak{K}_{\mathscr{D}}$ vanish when $\mathfrak{K}_{\mathscr{D}} = \emptyset$. Finally, (22) follows from the definition of $\mathfrak{S}_{\mathscr{D}}(\mathbb{S}^{d-1}, \Xi_M)$ (see Definition 7) and the fact that when $\mathfrak{K}_{\mathscr{D}} = \emptyset$ the spline coefficients are unconstrained (see Remark 5). □

## Appendix B. Proofs of Section 3

**Proof of Proposition 6.** Notice first that $\mathscr{D}$ maps isometrically $(\mathscr{C}(\mathbb{S}^{d-1}), \|\cdot\|_\infty)$ onto $(\mathscr{C}_{\mathscr{D}}(\mathbb{S}^{d-1}), \|\cdot\|_{\mathscr{D},\infty})$. Indeed, every element $h$ of $\mathscr{C}_{\mathscr{D}}(\mathbb{S}^{d-1})$ can be uniquely written as $h = \mathscr{D}\eta$, $\eta \in \mathscr{C}(\mathbb{S}^{d-1})$, with $\|h\|_{\mathscr{D},\infty} = \|\mathscr{D}^{-1}h\|_\infty = \|\eta\|_\infty$. We have hence the isometries

$$(\mathscr{C}(\mathbb{S}^{d-1}), \|\cdot\|_\infty) \cong (\mathscr{C}_{\mathscr{D}}(\mathbb{S}^{d-1}), \|\cdot\|_{\mathscr{D},\infty}) \text{ and } (\mathcal{M}_{\mathscr{D}}(\mathbb{S}^{d-1}), \|\cdot\|_{\mathscr{D},TV}) \cong (\mathcal{M}(\mathbb{S}^{d-1}), \|\cdot\|_{TV}).$$

Moreover, we have from the *Riesz-Markov representation theorem* the duality pair $(\mathscr{C}(\mathbb{S}^{d-1}), \|\cdot\|_\infty)' \cong (\mathcal{M}(\mathbb{S}^{d-1}), \|\cdot\|_{TV})$, which yields

$$(\mathscr{C}_{\mathscr{D}}(\mathbb{S}^{d-1}), \|\cdot\|_{\mathscr{D},\infty})' \cong (\mathscr{C}(\mathbb{S}^{d-1}), \|\cdot\|_\infty)' \cong (\mathcal{M}(\mathbb{S}^{d-1}), \|\cdot\|_{TV}) \cong (\mathcal{M}_{\mathscr{D}}(\mathbb{S}^{d-1}), \|\cdot\|_{\mathscr{D},TV}),$$

and hence $(\mathscr{C}_{\mathscr{D}}(\mathbb{S}^{d-1}), \|\cdot\|_{\mathscr{D},\infty})' \cong (\mathcal{M}_{\mathscr{D}}(\mathbb{S}^{d-1}), \|\cdot\|_{\mathscr{D},TV})$ as claimed. □

**Proof of Theorem 2.** We show Theorem 2 as a corollary of Lemma 1 below, pertaining to penalised convex optimisation problems in non strictly convex abstract Banach spaces. The latter leverages the *Krein-Milman* *theorem* [83, p. 75] as well as results from [28,8,27] to characterise the solution set $\mathcal{V}$ of (26) as the weak* closed convex hull of extreme points with bounded df. Our result generalises [28, Theorem 6] to the case of non strictly convex cost functionals $F$, often encountered in practice.

**Lemma 1** *(Extreme point representer theorem). Consider the following assumptions:*

*C1* $(\mathscr{B}, \|\cdot\|_{\mathscr{B}})$ *is a* Banach space, *with* topological dual $(\mathscr{B}', |||\cdot|||)$;
*C2* $span\{\varphi_i, i = 1, \dots, L\} \subset \mathscr{B}$, *with the* $\varphi_i$ *being* linearly independent;
*C3* $\mathbf{\Phi} : \mathscr{B}' \to \mathbb{C}^L$ *is a* sampling *operator, defined as in* (24);
*C4* $F : \mathbb{C}^L \times \mathbb{C}^L \to \mathbb{R}_+ \cup \{+\infty\}$ *is a cost functional as in* (27);
*C5* $\Lambda : \mathbb{R}_+ \to \mathbb{R}_+$ *is some* arbitrary strictly increasing convex *function.*

*Then, for any* $\boldsymbol{y} \in \mathbb{C}^L$*, the solution set of the optimisation problem*

$$\mathcal{V} = \operatorname*{Arg\,min}_{f\in\mathscr{B}'} \left\{F\left(\boldsymbol{y}, \mathbf{\Phi}(f)\right) \;+\; \Lambda\left(|||f|||\right)\right\}, \tag{B.1}$$

*is* non-empty *and the* weak* closed convex hull *of its* extreme points. *The latter are moreover necessarily of the form:*

$$f^\star = \sum_{m=1}^{M} \alpha_m e_m, \tag{B.2}$$

*where* $1 \le M \le L$, $\{\alpha_1, \dots, \alpha_M\} \subset \mathbb{C}$ *and* $e_m$ *are* extreme points *of the closed unit regularisation ball* $\mathcal{B} := \{f \in \mathscr{B}' : \Lambda(|||f|||) \le 1\}$.

**Proof.** Using the exact same arguments as in part i) of [28, Theorem 5] (which remain valid under the assumptions of Lemma 1), one can show that the functional $f \mapsto F\left(\boldsymbol{y}, \mathbf{\Phi}(f)\right) + \Lambda\left(|||f|||\right)$ is *proper, weak* lower semi-continuous, convex and coercive* on $\mathscr{B}'$. From [8, Proposition 8] the solution set $\mathcal{V}$ is hence *non-empty, convex* and *weak* compact*. Since $\mathscr{B}'$ equipped with the weak* topology is *locally convex and Hausdorff*, we can moreover invoke the *Krein-Milman theorem* [83, p. 75] to conclude that $\mathcal{V}$ is the *weak* closed convex hull of its extreme points*. In particular it has extreme points. Let $f_e \in \mathcal{V}$ be an arbitrary extreme point of $\mathcal{V}$ and let $\boldsymbol{z}_e := \mathbf{\Phi}(f_e) \in \mathbb{C}^L$. Then $f_e$ is also in the solution set of the generalised interpolation problem

$$\mathcal{V}_e = \operatorname*{Arg\,min}_{f\in\mathscr{B}'} \left\{\Lambda\left(|||f|||\right) \quad \text{s.t.} \quad \mathbf{\Phi}(f) = \boldsymbol{z}_e\right\}. \tag{B.3}$$

Using [27, Theorem 3.1] (with $j = 0$) we can moreover show that extreme points of $\mathcal{V}_e$ are of the form (B.2). Since $\mathcal{V}_e \subset \mathcal{V}$ and $f_e \in \mathcal{V}_e$, $f_e$ is also an extreme point of $\mathcal{V}_e$ and hence is indeed of the form (B.2). This shows that every extreme point of $\mathcal{V}$ is of the form (B.2), which achieves the proof. □

We now apply Lemma 1 to the FPBP problem (35) in Theorem 2. To this end, we set $(\mathscr{B}, \|\cdot\|_{\mathscr{B}}) = (\mathscr{C}_{\mathscr{D}}(\mathbb{S}^{d-1}), \|\mathscr{D}^{-1}\cdot\|_\infty)$, $(\mathscr{B}', |||\cdot|||) = (\mathcal{M}_{\mathscr{D}}(\mathbb{S}^{d-1}), \|\mathscr{D}\cdot\|_{TV})$ and $\Lambda(t) = \lambda t$. The assumptions of Lemma 1 are then indeed verified since $(\mathscr{C}_{\mathscr{D}}(\mathbb{S}^{d-1}), \|\mathscr{D}^{-1}\cdot\|_\infty)$ and $(\mathcal{M}_{\mathscr{D}}(\mathbb{S}^{d-1}), \|\mathscr{D}\cdot\|_{TV})$ form a duality pair of Banach spaces, and $\Lambda$ is convex and strictly increasing. We get hence from Theorem 1 that the solution set to the FPBP problem (35) is *nonempty* and the *weak* closed convex hull* of its *extreme points*. The latter are moreover necessarily of the form:

$$f^\star = \sum_{i=1}^{M} \beta_i e_i, \tag{B.4}$$

where $1 \leq M \leq L$, $\{\beta_1, \dots, \beta_M\} \subset \mathbb{C}$ and $e_i \in \mathscr{S}'(\mathbb{S}^{d-1})$ are *extreme points* of the closed regularisation ball $\mathcal{B}_{gTV,1/\lambda} = \{f \in \mathcal{M}_{\mathscr{D}}(\mathbb{S}^{d-1}) : \|\mathscr{D}f\|_{TV} \leq 1/\lambda\}$. We now compute the extreme points of $\mathcal{B}_{gTV,1/\lambda}$. We denote by $\delta\mathcal{V}$ the extreme points of an arbitrary convex set $\mathcal{V}$. We define moreover the gTV unit ball on $\mathcal{M}_{\mathscr{D}}(\mathbb{S}^{d-1})$ $\mathcal{B}_{gTV} = \{f \in \mathcal{M}_{\mathscr{D}}(\mathbb{S}^{d-1}) : \|\mathscr{D}f\|_{TV} \leq 1\}$, as well as the TV unit ball on the space of $\mathbb{C}$-valued regular Borel measures $\mathcal{M}(\mathbb{S}^{d-1})$: $\mathcal{B}_{TV} = \{f \in \mathcal{M}(\mathbb{S}^{d-1}) : \|f\|_{TV} \leq 1\}$. First, we trivially have $\delta\mathcal{B}_{gTV,1/\lambda} = \lambda^{-1}\delta\mathcal{B}_{gTV}$. Second, we have shown in Section 4.2 that $\mathscr{D}^{-1}$ is an *isometric isomorphism* from $\mathcal{M}(\mathbb{S}^{d-1})$ to $\mathcal{M}_{\mathscr{D}}(\mathbb{S}^{d-1})$, which yields $\delta\mathcal{B}_{gTV} = \mathscr{D}^{-1}(\delta\mathcal{B}_{TV})$ [27]. Finally, it is well known [28, Section 3.5] that extremes points of the total variation unit ball of complex regular Borel measures are of the form $z\delta_{\boldsymbol{r}}$ with $\boldsymbol{r} \in \mathbb{S}^{d-1}$ and $z \in \mathbb{C}$, $|z| = 1$. In conclusion, we hence get

$$\delta\mathcal{B}_{gTV,1/\lambda} = \lambda^{-1}\delta\mathcal{B}_{gTV} = \lambda^{-1}\mathscr{D}^{-1}(\delta\mathcal{B}_{TV}) = \left\{z\lambda^{-1}\mathscr{D}^{-1}\delta_{\boldsymbol{r}},\ \boldsymbol{r} \in \mathbb{S}^{d-1}, |z| = 1\right\}. \tag{B.5}$$

Plugging (B.5) into (B.4) allows us to write any extreme point of the solution set (35) as

$$f^\star = \sum_{i=1}^{M} \frac{\beta_i z_i}{\lambda}\mathscr{D}^{-1}\delta_{\boldsymbol{r}_i} = \sum_{i=1}^{M} \alpha_i \mathscr{D}^{-1}\delta_{\boldsymbol{r}_i},$$

for some constants $\{\alpha_1, \dots, \alpha_M\} \subset \mathbb{C}$ and directions $\{\boldsymbol{r}_1, \dots, \boldsymbol{r}_M\} \subset \mathbb{S}^{d-1}$, and where $1 \leq M \leq L$. This provides us with the first equality in (36). The second equality in (36) follows trivially from the definition of Green functions (see Definition 5). The last equality in (36) is obtained by considering the gSHT of $\mathscr{D}^{-1}\delta_{\boldsymbol{r}}$, $\boldsymbol{r} \in \mathbb{S}^{d-1}$. We have indeed, $\langle \mathscr{D}^{-1}\delta_{\boldsymbol{r}} | Y_n^m \rangle = \langle \delta_{\boldsymbol{r}} | \mathscr{D}^{-1} Y_n^m \rangle = Y_n^m(\boldsymbol{r})/\hat{D}_n$, $\forall n \in \mathbb{N}$, $m = 1, \dots, N_d(n)$, and hence $\mathscr{D}^{-1}\delta_{\boldsymbol{r}} = \sum_{n\in\mathbb{N}}\sum_{m=1}^{N_d(n)} \hat{D}_n^{-1} Y_n^m(\boldsymbol{r}) Y_n^m$. Finally, we have from Proposition 3 that

$$f^\star = \sum_{i=1}^{M} \alpha_i \Psi_{\boldsymbol{r}_i}^{\mathscr{D}} = \sum_{i=1}^{M} \alpha_i \psi_{\mathscr{D}}(\langle \cdot, \boldsymbol{r}_i \rangle),$$

where $\psi_{\mathscr{D}}$ is the zonal Green kernel of $\mathscr{D}$ and where the equality is in the sense of (15). Moreover, when $\mathscr{D}$ is spline-admissible, all traces $\{\psi_{\mathscr{D}}(\langle \cdot, \boldsymbol{r} \rangle),\ \boldsymbol{r} \in \mathbb{S}^{d-1}\}$ of the zonal Green kernel are ordinary functions and hence $f^\star$ is also an ordinary function. This shows (37) and achieves the proof of Theorem 2. □

**Proof of Proposition 7.** According to Proposition 6, the Dirac measures $\{\delta_{\boldsymbol{r}},\ \boldsymbol{r} \in \mathbb{S}^{d-1}\}$ belong to $\mathscr{C}_{\mathscr{D}}(\mathbb{S}^{d-1})$ i.f.f. there exists, for every $\boldsymbol{r} \in \mathbb{S}^{d-1}$, a function $\Psi_{\boldsymbol{r}} \in \mathscr{C}(\mathbb{S}^{d-1})$ such that: $\mathscr{D}\Psi_{\boldsymbol{r}} = \delta_{\boldsymbol{r}}$. The functions $\Psi_{\boldsymbol{r}}$ satisfying the above equation are actually the Green functions of $\mathscr{D}$, which, from Proposition 3, can be identified with traces $\psi_{\mathscr{D}}(\langle \cdot, \boldsymbol{r} \rangle)$ of the zonal Green kernel of $\mathscr{D}$. Finally, we have shown in Proposition 4 that, for a pseudo-differential operator with spectral growth order $p > d - 1$ $\{\psi_{\mathscr{D}}(\langle \cdot, \boldsymbol{r} \rangle),\ \boldsymbol{r} \in \mathbb{S}^{d-1}\} \subset \mathscr{C}(\mathbb{S}^{d-1})$, and hence $\Psi_{\boldsymbol{r}}$ can indeed be identified with a continuous function and consequently all Dirac measures belong to the predual $\mathscr{C}_{\mathscr{D}}(\mathbb{S}^{d-1})$. □

**Proof of Proposition 8.** From Proposition 6, a function $f \in \mathscr{L}^2(\mathbb{S}^{d-1})$ is in $\mathscr{C}_{\mathscr{D}}(\mathbb{S}^{d-1})$ if there exists $\eta \in \mathscr{C}(\mathbb{S}^{d-1})$ s.t. $f = \mathscr{D}\eta$. Since $\mathscr{D}$ is assumed invertible, this is equivalent to requiring that $\mathscr{D}^{-1}f \in \mathscr{C}(\mathbb{S}^{d-1})$, which is guaranteed if the series of functions

$$(\mathscr{D}^{-1}f)(\boldsymbol{r}) = \sum_{n\in\mathbb{N}} \frac{1}{\hat{D}_n} \sum_{m=1}^{N_d(n)} \hat{f}_n^m Y_n^m(\boldsymbol{r}), \quad \boldsymbol{r} \in \mathbb{S}^{d-1}, \tag{B.6}$$

converges *uniformly* (see [21, Theorem 2.14]). To show that (B.6) is uniformly convergent, we consider its remainder for some $N \in \mathbb{N}$. Then, from the addition theorem [21, Theorem 5.11] and the Cauchy-Schwarz inequality we get, for each $\boldsymbol{r} \in \mathbb{S}^{d-1}$:

$$\left|\sum_{n=N}^{+\infty}\frac{1}{\hat{D}_n}\sum_{m=1}^{N_d(n)}\hat{f}_n^m Y_n^m(\boldsymbol{r})\right| \leq \left|\sum_{n=N}^{+\infty}\frac{\sum_{m=1}^{N_d(n)}|Y_n^m(\boldsymbol{r})|^2}{|\hat{D}_n|^2}\right|\left|\sum_{n=N}^{+\infty}\sum_{m=1}^{N_d(n)}|\hat{f}_n^m|^2\right| = \left|\sum_{n=N}^{+\infty}\frac{N_d(n)}{\mathfrak{a}_d|\hat{D}_n|^2}\right|\left|\sum_{n=N}^{+\infty}\sum_{m=1}^{N_d(n)}|\hat{f}_n^m|^2\right|.$$

Since $f \in \mathscr{L}^2(\mathbb{S}^{d-1})$ we have trivially $\lim_{N\to+\infty}\left|\sum_{n=N}^{+\infty}\sum_{m=1}^{N_d(n)}|\hat{f}_n^m|^2\right| = 0$. Moreover, since $|\hat{D}_n| = \Theta(n^p)$ we have from (6) $N_d(n)|\hat{D}_n|^{-2} = \mathcal{O}(n^{d-2-2p})$. Since $p > (d-1)/2 \Rightarrow d-2-2p < -1$ we have hence that the series $\sum_{n\in\mathbb{N}}\frac{N_d(n)}{\mathfrak{a}_d|\hat{D}_n|^2}$ is convergent and hence is remainder tends to zero. We have hence

$$\left|(\mathscr{D}^{-1}f)(\boldsymbol{r}) - \sum_{n=0}^{N-1}\frac{1}{\hat{D}_n}\sum_{m=1}^{N_d(n)}\hat{f}_n^m Y_n^m(\boldsymbol{r})\right| = \left|\sum_{n=N}^{+\infty}\frac{1}{\hat{D}_n}\sum_{m=1}^{N_d(n)}\hat{f}_n^m Y_n^m(\boldsymbol{r})\right|$$
$$\leq \left|\sum_{n=N}^{+\infty}\frac{N_d(n)}{\mathfrak{a}_d|\hat{D}_n|^2}\right|\left|\sum_{n=N}^{+\infty}\sum_{m=1}^{N_d(n)}|\hat{f}_n^m|^2\right| \stackrel{N\to+\infty}{\longrightarrow} 0.$$

Moreover, since the upper bound is independent on $\boldsymbol{r}$ the convergence is uniform, which achieves the proof. □

## Appendix C. Proofs of Section 5

**Proof of Proposition 9.** It is easy to see that the Hermitian square-root of $\mathscr{D}$ has Fourier symbol $\{\sqrt{\hat{D}_n},\ n \in \mathbb{N}\}$. The latter is moreover a pseudo-differential operator since the Fourier coefficients $\hat{D}_n$ are all positive (from the assumption of positive-definiteness of $\mathscr{D}$) and hence $\{\sqrt{\hat{D}_n},\ n \in \mathbb{N}\} \subset \mathbb{R}_+$. The rest of the assumptions of Definition 4 trivially follow from $\mathscr{D}$ being a pseudo-differential operator. Moreover, from the assumption $p > d-1$ we get that the spectral growth order of $\mathscr{D}^{1/2}$, equal to $p/2$, is strictly larger than $(d-1)/2$. We can hence apply [41, Lemma 5.5] to conclude that the generalised Sobolev space $\mathscr{H}_{\mathscr{D}^{1/2}}(\mathbb{S}^{d-1})$ is an RKHS, containing all Dirac measures in its dual. Moreover, using the same arguments as in the proof of [41, Theorem 5.3], it is possible to show that the Riesz map $R_{\mathscr{H}_{\mathscr{D}^{1/2}}} : \mathscr{H}'_{\mathscr{D}^{1/2}}(\mathbb{S}^{d-1}) \to \mathscr{H}_{\mathscr{D}^{1/2}}(\mathbb{S}^{d-1})$ is $\mathscr{D}^{-1}$. Therefore, $\mathscr{D}$-splines[38] are all contained in $\mathscr{H}_{\mathscr{D}^{1/2}}(\mathbb{S}^{d-1})$ as images[39] by the Riesz map $\mathscr{D}^{-1}$ of elements of the dual, namely linear combinations of Dirac measures. □

**Proof of Proposition 10.** Proposition 9 tells us that $\mathscr{H}_{\mathscr{D}^{1/2}}(\mathbb{S}^{d-1})$ is an RKHS. Therefore, any element $h \in \mathscr{H}_{\mathscr{D}^{1/2}}(\mathbb{S}^{d-1})$ is an ordinary function, and we have from Proposition 3

$$\langle h, \psi_{\mathscr{D}}(\langle\cdot,\boldsymbol{r}\rangle)\rangle_{\mathscr{D}^{1/2}} = \left\langle \mathscr{D}^{1/2}h, \mathscr{D}^{1/2}\psi_{\mathscr{D}}(\langle\cdot,\boldsymbol{r}\rangle)\right\rangle_{\mathbb{S}^{d-1}} = \left\langle \mathscr{D}^{1/2}\Psi_{\boldsymbol{r}}^{\mathscr{D}}\middle|\mathscr{D}^{1/2}h\right\rangle = \left\langle \mathscr{D}\Psi_{\boldsymbol{r}}^{\mathscr{D}}\middle|h\right\rangle = \langle\delta_{\boldsymbol{r}}|h\rangle = h(\boldsymbol{r}),$$

for all $\boldsymbol{r} \in \mathbb{S}^{d-1}$, which shows that the zonal Green kernel $\psi_{\mathscr{D}}$ is the *reproducing kernel* [84] of $\mathscr{H}_{\mathscr{D}^{1/2}}(\mathbb{S}^{d-1})$. Additionally, since $\mathscr{D}$ is positive-definite, it is in particular invertible, and we get from (22) that

$$\mathfrak{S}_{\mathscr{D}}(\mathbb{S}^{d-1}, \Xi_N) = \operatorname{span}\{\psi_{\mathscr{D}}^n := \psi_{\mathscr{D}}(\langle\cdot,\boldsymbol{r}_n\rangle), \boldsymbol{r}_n \in \Xi_N\}.$$

The positive-definiteness of $\mathscr{D}$ implies moreover (see [41, Remark 5.9] and [21, Theorem 6.27]) that the family of functions $\{\psi_{\mathscr{D}}^n, n = 1,\ldots,N\}$ is *linearly independent* and hence forms a basis for $\mathfrak{S}_{\mathscr{D}}(\mathbb{S}^{d-1}, \Xi_N)$. Consequently, the orthogonal projection of $h$ onto $\mathfrak{S}_{\mathscr{D}}(\mathbb{S}^{d-1}, \Xi_N)$ can be written as

[38] $\mathscr{D}$-splines exist indeed since $p > d-1$ implies that $\mathscr{D}$ is spline-admissible from Proposition 4.
[39] From Definition 7, a $\mathscr{D}$-spline is such that $\mathscr{D}\mathfrak{s} = \sum_{i=1}^N \alpha_i\delta_{\boldsymbol{r}_i}$ which, for $\mathscr{D}$ invertible, is equivalent to $\mathfrak{s} = \mathscr{D}^{-1}(\sum_{i=1}^N \alpha_i\delta_{\boldsymbol{r}_i})$.

$$\mathfrak{s}_N^\perp = \sum_{n=1}^N \langle h, \psi_{\mathscr{D}}^n \rangle_{\mathscr{D}^{1/2}} \widetilde{\psi}_{\mathscr{D}}^n = \sum_{n=1}^N h(\boldsymbol{r}_n)\widetilde{\psi}_{\mathscr{D}}^n,$$

where the second equality follows from the fact that $\psi_{\mathscr{D}}$ reproduces functions in $\mathscr{H}_{\mathscr{D}}^{1/2}(\mathbb{S}^{d-1})$ and $\{\widetilde{\psi}_{\mathscr{D}}^n,\, n = 1, \dots, N\} \subset \mathfrak{S}_{\mathscr{D}}(\mathbb{S}^{d-1}, \Xi_N)$ is the *dual basis* [44, Chapter 2] of $\{\psi_{\mathscr{D}}^n,\, n = 1, \dots, N\}$, verifying the *biorthogonality property* $\left\langle \widetilde{\psi}_{\mathscr{D}}^m, \psi_{\mathscr{D}}^n \right\rangle_{\mathscr{D}^{1/2}} = \delta_{mn}$, $\forall m, n = 1, \dots, N$. We have hence

$$\langle \mathfrak{s}_N^\perp, \psi_{\mathscr{D}}^n \rangle_{\mathscr{D}^{1/2}} = h(\boldsymbol{r}_n) = \langle h, \psi_{\mathscr{D}}^n \rangle_{\mathscr{D}^{1/2}}, \quad n = 1, \dots, N. \tag{C.1}$$

Moreover, [21, Lemma 6.34] tells us that, for spline-admissible pseudo-differential operators with growth order $p > \frac{d+1}{2}$, the zonal Green kernel $\psi_{\mathscr{D}}$ is *uniformly Lipschitz continuous*, i.e. there exists $L_{\mathscr{D}} > 0$ which only depends on the sequence $\{\hat{D}_n\}_{n\in\mathbb{N}}$ such that for any $\boldsymbol{\rho} \in \mathbb{S}^{d-1}$

$$|\psi_{\mathscr{D}}(\langle \boldsymbol{r}, \boldsymbol{\rho} \rangle) - \psi_{\mathscr{D}}(\langle \boldsymbol{s}, \boldsymbol{\rho} \rangle)| \leq L_{\mathscr{D}}^2 \|\boldsymbol{r} - \boldsymbol{s}\|_2, \quad \forall \boldsymbol{r}, \boldsymbol{s} \in \mathbb{S}^{d-1}. \tag{C.2}$$

With these two observations, we are now ready to prove the result. First, we get from (C.1) as well as the Cauchy-Schwarz and triangle inequalities

$$\begin{aligned} \left|h(\boldsymbol{r}) - \mathfrak{s}_N^\perp(\boldsymbol{r})\right| &= |h(\boldsymbol{r}) - h(\boldsymbol{r}_n) + \mathfrak{s}_N^\perp(\boldsymbol{r}_n) - \mathfrak{s}_N^\perp(\boldsymbol{r})| = |\left\langle \psi_{\mathscr{D}}^{\boldsymbol{r}} - \psi_{\mathscr{D}}^{\boldsymbol{r}_n}, h - h_N^\perp \right\rangle_{\mathscr{D}^{1/2}}| \\ &\leq \|\psi_{\mathscr{D}}^{\boldsymbol{r}} - \psi_{\mathscr{D}}^{\boldsymbol{r}_n}\|_{\mathscr{D}^{1/2}} (\|h_N^\perp\|_{\mathscr{D}^{1/2}} + \|h\|_{\mathscr{D}^{1/2}}) \leq 2\|\psi_{\mathscr{D}}^{\boldsymbol{r}} - \psi_{\mathscr{D}}^{\boldsymbol{r}_n}\|_{\mathscr{D}^{1/2}} \|h\|_{\mathscr{D}^{1/2}}. \end{aligned}$$

Second, we obtain from the reproducing property, equation (C.2) and the definition (42) of the nodal width $\Theta_{\Xi_N}$:

$$\begin{aligned} \|\psi_{\mathscr{D}}^{\boldsymbol{r}} - \psi_{\mathscr{D}}^{\boldsymbol{r}_n}\|_{\mathscr{D}^{1/2}}^2 &= \langle \psi_{\mathscr{D}}^{\boldsymbol{r}} - \psi_{\mathscr{D}}^{\boldsymbol{r}_n}, \psi_{\mathscr{D}}^{\boldsymbol{r}} - \psi_{\mathscr{D}}^{\boldsymbol{r}_n} \rangle_{\mathscr{D}^{1/2}} = \psi_{\mathscr{D}}(\langle \boldsymbol{r}, \boldsymbol{r} \rangle) + \psi_{\mathscr{D}}(\langle \boldsymbol{r}_n, \boldsymbol{r}_n \rangle) - \psi_{\mathscr{D}}(\langle \boldsymbol{r}_n, \boldsymbol{r} \rangle) - \psi_{\mathscr{D}}(\langle \boldsymbol{r}, \boldsymbol{r}_n \rangle) \\ &\leq 2L_{\mathscr{D}}^2 \|\boldsymbol{r} - \boldsymbol{r}_n\|_{\mathbb{R}^d} \leq 2L_{\mathscr{D}}^2 \Theta_{\Xi_N}. \end{aligned}$$

In conclusion, this yields: $\sup_{\boldsymbol{r}\in\mathbb{S}^{d-1}} |h(\boldsymbol{r}) - \mathfrak{s}_N^\perp(\boldsymbol{r})| \leq 2^{3/2} L_{\mathscr{D}} \sqrt{\Theta_{\Xi_N}} \|h\|_{\mathscr{D}^{1/2}}$, which achieves the proof. □

**Proof of Theorem 3.** The spline-admissible pseudo-differential operator $\mathscr{D}$ being positive-definite, its Green kernel $\psi_{\mathscr{D}}$ is strictly positive-definite (see [21, Definition 6.25 and Theorem 6.27]) and hence according to [21, Lemma 6.26], the family of functions $\{\psi_n = \psi_{\mathscr{D}}(\langle \cdot, \boldsymbol{r}_n \rangle),\, n = 1, \dots, N\}$ is *linearly independent* for every set $\Xi_N = \{\boldsymbol{r}_1, \dots, \boldsymbol{r}_n\} \subset \mathbb{S}^{d-1}$ of $N$ distinct points. The synthesis operator $\Psi$ defines hence a bijection between $\mathbb{C}^N$ and $\mathfrak{S}_{\mathscr{D}}(\mathbb{S}^{d-1}, \Xi_N) = \operatorname{span}\{\psi_n,\, n = 1, \dots, N\}$. From this isomorphism, we get notably

$$\begin{aligned} \mathcal{V} &= \operatorname*{Arg\,min}_{f\in\mathfrak{S}_{\mathscr{D}}(\mathbb{S}^{d-1}, \Xi_N)} \{F(\boldsymbol{y}, \boldsymbol{\Phi}(f)) \quad + \quad \lambda\|\mathscr{D}f\|_{TV}\} \\ &= \Psi\left(\operatorname*{Arg\,min}_{\boldsymbol{x}\in\mathbb{C}^N} \{F(\boldsymbol{y}, \boldsymbol{\Phi}\Psi(\boldsymbol{x})) \quad + \quad \lambda\|\mathscr{D}\Psi(\boldsymbol{x})\|_{TV}\}\right) \\ &= \Psi\left(\operatorname*{Arg\,min}_{\boldsymbol{x}\in\mathbb{C}^N} \{F(\boldsymbol{y}, \boldsymbol{\Phi}\Psi(\boldsymbol{x})) \quad + \quad \lambda\|\boldsymbol{x}\|_1\}\right), \end{aligned} \tag{C.3}$$

since we have

$$\|\mathscr{D}\Psi(\boldsymbol{x})\|_{TV} = \left\|\sum_{n=1}^N x_n \mathscr{D}\psi_{\mathscr{D}}(\langle \cdot, \boldsymbol{r}_n \rangle)\right\|_{TV} = \left\|\sum_{n=1}^N x_n \delta_{\boldsymbol{r}_n}\right\|_{TV} = \|\boldsymbol{x}\|_1.$$

Notice that the linear operator $\boldsymbol{\Phi}\Psi : \mathbb{C}^N \to \mathbb{C}^L$ is finite-dimensional, and can hence be represented as a matrix. From the bilinearity of the Schwartz duality product, we have indeed

$$(\boldsymbol{\Phi}\Psi\boldsymbol{x})_l = \langle \Psi\boldsymbol{x}|\varphi_l\rangle = \left\langle \sum_{n=1}^N x_n \psi_{\mathscr{D}}(\langle\cdot,\boldsymbol{r}_n\rangle)\middle|\varphi_l\right\rangle = \sum_{n=1}^N x_n \left\langle \psi_{\mathscr{D}}(\langle\cdot,\boldsymbol{r}_n\rangle)\middle|\varphi_l\right\rangle = \sum_{n=1}^N x_n G_{ln}, \quad \forall l = 1,\dots,L.$$

We can hence identify $\boldsymbol{\Phi}\Psi$ with a matrix $\boldsymbol{G} \in \mathbb{C}^{L\times N}$, with entries given by $G_{ln} := \langle \psi_{\mathscr{D}}(\langle\cdot,\boldsymbol{r}_n\rangle)|\varphi_l\rangle$, $l = 1,\dots,L$, $n = 1,\dots,N$. Since $\mathscr{D}$ is spline-admissible, the traces of the zonal Green kernel are ordinary functions and hence in particular square-integrable. When the sampling functionals $\{\varphi_1,\dots,\varphi_L\}$ are in $\mathscr{L}^2(\mathbb{S}^{d-1})$ we can hence obtain a simpler expression for the entries of $\boldsymbol{G}$:

$$G_{ln} = \left\langle \psi_{\mathscr{D}}(\langle\cdot,\boldsymbol{r}_n\rangle)\middle|\varphi_l\right\rangle = \langle \varphi_l, \psi_{\mathscr{D}}(\langle\cdot,\boldsymbol{r}_n\rangle)\rangle_{\mathbb{S}^{d-1}} = (\psi_{\mathscr{D}} * \varphi_l)(\boldsymbol{r}_n), \quad l = 1,\dots,L, \quad n = 1,\dots,N.$$

Equation (C.3) finally reduces to

$$\mathcal{V} = \Psi\left(\operatorname*{Arg\,min}_{\boldsymbol{x}\in\mathbb{C}^N}\{F(\boldsymbol{y},\boldsymbol{G}\boldsymbol{x}) \quad + \quad \lambda\|\boldsymbol{x}\|_1\}\right) = \Psi\left(\mathfrak{U}\right),$$

as claimed. From [41, Proposition 6.1] we furthermore get $\mathfrak{U} = \boldsymbol{\Psi}^\dagger(\mathcal{V}) = (\boldsymbol{\Psi}^*\Psi)^{-1}\boldsymbol{\Psi}^*(\mathcal{V})$, which concludes the proof. □

## Appendix D. Supplementary material

Supplementary material related to this article can be found online at https://doi.org/10.1016/j.acha.2020.12.004.